\documentclass[11pt]{article}

\usepackage[margin=1in]{geometry}
\usepackage{amsmath,amsthm,amssymb}
\usepackage{commath,mathtools}
\usepackage{algorithm,algorithmic}
\usepackage{graphicx}
\usepackage{subcaption}
\usepackage{hyperref}
\usepackage{color}
\usepackage{rotating}
\usepackage{lscape}
\usepackage[sort,nocompress]{cite}
\usepackage{verbatim}

\providecommand{\keywords}[1]{\small \textbf{\textit{Keywords---}} #1}

\newtheorem{theorem}{Theorem}

\newtheorem{lemma}[theorem]{Lemma}
\newtheorem*{remark}{Remarks}

\usepackage[pagewise,mathlines]{lineno}

\DeclareMathOperator*{\argmin}{\mathrm{arg\,min}}

\DeclareMathOperator*{\minimize}{\mathrm{minimize}}

\newcommand{\Lcal}{\mathcal{L}}
\newcommand{\Acal}{\mathcal{A}}
\newcommand{\Bcal}{\mathcal{B}}

\newcommand{\Lint}{L_{\text{int}}}
\newcommand{\Lbdry}{L_{\text{bdry}}}

\newcommand{\tj}{\theta_{j}}
\newcommand{\tjp}{\theta_{j+1}}
\newcommand{\btj}{\bar{\theta}_{j}}
\newcommand{\btjp}{\bar{\theta}_{j+1}}
\newcommand{\be}{\bar{e}}
\newcommand{\Gcal}{\mathcal{G}}

\title{Numerical Solution of Inverse Problems by Weak Adversarial Networks}

\author{
Gang Bao
\footnote{School of Mathematical Sciences, Zhejiang University, Hangzhou, Zhejiang, China. Email: \texttt{baog@zju.edu.cn}.}
\and
Xiaojing Ye
\footnote{Department of Mathematics and Statistics, Georgia State University, Atlanta, GA, 30303, USA.
Email: \texttt{xye@gsu.edu}.}
\and
Yaohua Zang
\footnote{School of Mathematical Sciences, Zhejiang University, Hangzhou, Zhejiang, China. Email: \texttt{yhchuang@zju.edu.cn}.}
\and
Haomin Zhou
\footnote{School of Mathematics, Georgia Institute of Technology, Atlanta, GA, 30332, USA.
Email: \texttt{hmzhou@math.gatech.edu}.}
}
\date{}

\begin{document}
\maketitle

\begin{abstract}
In this paper, a weak adversarial network approach is developed to numerically solve a class of inverse problems, including electrical impedance tomography and dynamic electrical impedance tomography problems. The weak formulation of the PDE for the given inverse problem is leveraged, where the solution and the test function are parameterized as deep neural networks. Then, the weak formulation and the boundary conditions induce a minimax problem of a saddle function of the network parameters. As the parameters are alternatively updated, the network gradually approximates the solution of the inverse problem. Theoretical justifications are provided on the convergence of the proposed algorithm. The proposed method is completely mesh-free without any spatial discretization, and is particularly suitable for problems with high dimensionality and low regularity on solutions. Numerical experiments on a variety of test inverse problems demonstrate the promising accuracy and efficiency of this approach.
\end{abstract}


\section{Introduction}
%
Inverse problems (IP) are ubiquitous in a vast number of scientific disciplines, including geophysics \cite{Portal2016Contribution}, signal processing and imaging \cite{bertero1998introduction}, computer vision \cite{paragios2006handbook}, remote sensing and control \cite{yamamoto1995stability}, statistics \cite{kaipio2006statistical}, and machine learning \cite{goodfellow2016deep}.
Let $\Omega$ be an open and bounded set in $\mathbb{R}^d$, then an IP defined on $\Omega$ can be presented in a general form as:
\begin{subequations}
\label{eq:ip}
\begin{align}
\Acal[u,\gamma] = 0,& \quad \mbox{in}\ \Omega \label{eq:ip1a} \\
\Bcal[u,\gamma] = 0,& \quad \mbox{on}\ \partial \Omega \label{eq:ip1b}
\end{align}
\end{subequations}
where $\Acal[u,\gamma]$ specifies a differential equation, in which $u$ is the solution and $\gamma$ the coefficient in the inverse medium problem or the source function in the inverse source problem. Equation $\Acal$ can be an ordinary differential equation (ODE), or a partial differential equation (PDE), or an integro-differential equation (IDE), that $(u,\gamma)$ needs to satisfy (almost) everywhere inside the region $\Omega$.
The boundary value (and initial value if applicable) is given by $\Bcal[u,\gamma]$ on $\partial \Omega$.
Depending on specific applications, partial information of $u$ and/or $\gamma$ may be available in the interior of $\Omega$.
Then IP \eqref{eq:ip} is to find $(u,\gamma)$ that satisfies both \eqref{eq:ip1a} and \eqref{eq:ip1b}.

To instantiate our approach, we mostly use the classical inverse conductivity problem in electrical impedance tomography (EIT) \cite{Cheney1999Electrical,khan2019review} as an example to present our main idea and the derivations in this paper. However, our methodology can be readily applied to other classes of IPs with modifications. An example of dynamic EIT problem will be shown in Section \ref{sec:experiment}.
The goal of EIT is to determine the electrical conductivity distribution $\gamma(x)$ of an unknown medium defined on $\Omega$ based on the potential $u$, the current $-\gamma\partial_{\vec{n}} u$ measurements, and the knowledge of $\gamma$ (and hence $\partial_{\vec{n}} u$) on/near the boundary $\partial \Omega$ of the domain $\Omega$:
\begin{subequations}
\label{eq:eit}
\begin{align}
-\nabla \cdot ( \gamma \nabla u) - f = 0, & \quad \mbox{in}\ \Omega \label{eq:ip2a} \\
u - u_b = 0,\ \gamma - \gamma_b = 0,\ \partial_{\vec{n}}u - u_n = 0, & \quad \mbox{on}\ \partial \Omega \label{eq:ip2b}
\end{align}
\end{subequations}
where $u_b$ is the measured voltage, $\gamma_b$ is the conductivity near the surface of the object and $u_n\triangleq \nabla u\cdot\vec{n}$ with $\vec{n}$ being the outer normal of $\partial \Omega$.
Note that our approach is not to estimate the Dirichlet-to-Neumann (DtN) map associated with the conductivity function as in classical methods specific to the EIT problem \cite{Curtis1991The,Francini2000Recovering,klibanov2019convexification}.
Instead, our goal is to directly solve a general class of IPs \eqref{eq:ip} numerically using the given data, with the EIT problem \eqref{eq:eit} as a prototype example without exploiting its special structure (e.g., the DtN map).
To make our presentation concise and focused, we only consider IPs with $\Acal[u,\gamma]$ characterized by PDEs in \eqref{eq:ip1a}, and assume that the given IP is well-defined and admits at least one (weak) solution.

Our approach is to train deep neural networks that can represent the solution $(u,\gamma)$ of a given IP, with substantial improvement over classical numerical methods especially for problems with high dimensionality.
More specifically, we leverage the weak formulation of the PDE \eqref{eq:ip1a} and convert the IP into an operator norm minimization problem of $u$ and $\gamma$.
Then we parameterize both $u$, the unknown coefficient $\gamma$, and the test function $\varphi$ as deep neural networks $u_\theta$, $\gamma_\theta$, and $\varphi_{\eta}$ respectively, with network parameters $(\theta,\eta)$, and form a minimax problem of a saddle function of the parameters $(\theta,\eta)$.
Finally, we apply the stochastic gradient descent method to alternately update the network parameters so that $(u_{\theta},\gamma_{\theta})$ gradually approximates the solution of the IP.
The parameterization of $(u,\gamma)$ using deep neural networks requires no discretization of the spatial and temporal domain, and hence is completely mesh free.
This is a promising alternative compared to the classical finite difference method (FDM) and finite element methods (FEM) which suffer the issue of the so-called curse of dimensionality, a term first used in \cite{bellman1966dynamic}.
Moreover, our approach combines the training of the weak solution (primal network) $(u,\gamma)$ and the test function (adversarial network) $\varphi$ governed by the weak formulation of the PDE, which requires less regularity of the solution $(u,\gamma)$ and can be more advantageous in many real-world applications when the solution has singularities.

The remainder of this paper is organized as follows. We first review the recent work on deep learning based solutions to forward and inverse problems in Section \ref{sec:related_work}.
In Section \ref{sec:iwan}, we provide the detailed derivation of our method and a series of theoretical results to support the validity of the proposed approach.
We discuss several implementation techniques that can improve practical performance and conduct a series of numerical experiments to demonstrate the effectiveness of the proposed approach in Section \ref{sec:experiment}.
Section \ref{sec:conclusion} concludes this paper with some general remarks.

\section{Related Work}
\label{sec:related_work}
The past few years have witnessed an emerging trend of using deep learning based methods to solve forward and inverse problems.
These methods can be roughly classified into two categories.
The first category includes methods that approximate the solution of a given problem based on supervised learning approaches.
These methods require a large number of input-output pairs through numerical simulation and experiments to train the desired networks.
In this category, deep neural networks are used to generate approximate intermediate results from measurement data for further refinement \cite{michalikova2014image,martin2015nonlinear,fernandez2018towards,tan2018image,yao2019two,rymarczyk2019comparison}, applied to improve the solution of classical numerical methods in the post-processing phase \cite{Jin2017Deep,Kang2017A,antholzer2019deep,Hauptmann2017Model,martin2017post,hamilton2018deep,hamilton2019beltrami,wei2019dominant}, or approximate the mapping from given parameters of an inverse problem to its solution but require spatial discretization and cannot be applied to high-dimensional problems \cite{joshi2019generative,adler2017solving,li2020nett}.

The second category features unsupervised learning methods that directly solve the forward or inverse problem based on the problem formulation rather than additional training data, which can be more advantageous than those in the first category in practice.
For example, feed-forward neural networks are used to parameterize the coefficient functions and trained by minimizing the performance function in \cite{dadvand2006artificial}. In \cite{khoo2018switchnet}, a neural network architecture called SwitchNet is proposed to solve the inverse scattering problem through the mapping between the scatterers and the scattered field.
In \cite{fan2019solving}, a deep learning approach specific to 2D and 3D EIT problems is developed to represent the DtN map by a compact neural network architecture.
The backward stochastic differential equation (BSDE) corresponding to the PDE in a forward problem is parameterized in part by neural networks, such that the solution of the PDE can be obtained by integrating the BSDE for a target point in the domain \cite{han2018solving,e2017deep,beck2019machine}.
In \cite{e2018deep}, the solution of a forward problem is parameterized as a deep neural network, which is trained by minimizing the loss function composed of the energy functional associated with the PDE and a penalty term on the boundary value condition.
Another mesh-free framework, called physics-informed neural networks (PINN), for solving both the forward and inverse problems using deep neural networks based on the strong formulation of PDEs is proposed in \cite{raissi2019physics}, where a constant coefficient function is considered for the inverse problem part.
{Specifically, PINN parameterizes the unknowns of a given PDE using deep neural networks, which are trained by minimizing the loss function formed as the least squares of the violation of the PDE at sampled points in the domain and boundary conditions.}
Some empirical study of PINN is also conducted in \cite{dockhorn2019discussion}.
Solutions to IPs based on PINN with data given in problem domain are also considered in \cite{jo2019deep}, and refinement of solutions using adaptively sampled collocation points is proposed in \cite{anitescu2019artificial}.
In \cite{zang2019weak}, the weak formulation of the PDE is leveraged as the objective function, where the solution of the PDE and the test function are both parameterized as deep neural networks trying to minimize and maximize the objective function, respectively.
In \cite{khodayi2019varnet}, a similar variational form is used where the test function is fixed basis instead of neural networks to be learned.
In \cite{meng2019composite}, three neural networks, one for low-fidelity data and the other two for the linear and nonlinear functions for high-fidelity data, are used by following the PINN approach.
The PINN with a multi-fidelity network structure is also proposed for stochastic PDE cases, where polynomial chaotic expansions are used to express the solutions, i.e., as a linear combination of random basis with coefficient functions to be learned \cite{chen2019learning}.
In \cite{bar2019unsupervised}, the solution of an IP is parameterized by deep neural network and learned by minimizing a cost function that enforces the conditions of IP and additional regularization, where solutions to the PDE are required during the training.

Recently, meta-learning based approaches for forward problems are also considered \cite{chen2019learning,feliu2019meta,lu2019deeponet}.
In \cite{feliu2019meta}, the mapping from the coefficient of a differential operator to the pseudo-differential operator (e.g., the Green function) is learned by leveraging the compressed form of the wavelet transform.
In \cite{lu2019deeponet}, a deep operator network consisting of a branch network and a trunk network is introduced. The network encodes the input function evaluated at a finite number of locations (branch-net) and the locations for the output function (trunk-net) and, the output function is given by the inner product of the two plus a bias.
Learning network width and depth parameters are also considered using Bayesian optimization in \cite{chen2019learning}.

Our approach to the IP follows our earlier work \cite{zang2019weak} for forward problems, which differs from the aforementioned existing methods in the use of the weak formulation of PDEs.
The weak formulation is a powerful approach for solving PDEs as it requires less regularity and allows for necessary singularities of the solutions, which is an important feature appreciated in many real-world applications such as imaging and abnormality detections.
From the theoretical point of view, our method employs neural network parameterizations of both the solution (as the primal network) and the test function (as the adversarial network), and performs an adversarial training in a way that the test function critics on the solution network where the PDE is violated, and the solution network corrects itself at those spots until the PDE is satisfied (almost) everywhere in the domain.
{However, as inverse problems are often ill-posed and more difficult to solve than forward problems in general, we mostly focus on the inverse problem \eqref{eq:eit} in EIT in this work. Some experimental results on similar problems are also presented in Section \ref{sec:experiment}.}

{
The adversarial training in the present work has a similar flavor as the one used in generative adversarial network \cite{goodfellow2014generative}, where a generator network is aimed at mapping generic random samples (such as those from a given multivariate Gaussian) to ones following the same distribution as the training samples, and a discriminator network is to distinguish these samples produced by the generator network from the true samples.
The generator and adversarial networks act as the two players in a zero-sum game, and are alternately updated by gradient descent and ascent on the objective function respectively to reach an equilibrium.
In particular, a notable variant of GAN, called Wasserstein GAN \cite{arjovsky2017wasserstein}, also has a min-max structure of a primal network (generator) and adversarial network (dual function of optimal transport due to the Wasserstein distance between generated and sample distributions) as our formulation. However, WGAN requires its dual function in the max problem to be 1-Lipschitz, which is very difficult to realize numerically and has generated a series of followup work to overcome the issue \cite{gulrajani2017improved,miyato2018spectral}, Spectral Normalization for Generative Adversarial Networks}.
In contrast, the structure of weak solution versus test function in our work arises naturally from the weak formulation in the PDE theory, which enjoys numerous theoretical justifications and computational benefits for solving IPs for PDEs without imposing restrictive constraint on the adversarial network (test function), as we show in the present work.

{
In contrast to many existing deep learning methods that require a large amount of demonstration data (e.g., coefficient/boundary value and solution pairs) for training, our method follows an unsupervised learning strategy and only needs the formulation of the PDE and boundary conditions in the given IP. In \cite{ulyanov2018deep}, an unsupervised learning study reveals that generic convolutional neural networks (CNN) automatically bias towards smooth signals and can produce results similar to some sophisticated reconstructions in image denoising without any training data. This phenomenon, known as deep image prior (DIP), is further exploited in \cite{heckel2019denoising,dittmer2019regularization}. The most notable difference between DIP and the present work is that, our method is completely mesh-free and does not require any spatial discretization, which is suitable for high-dimensional problems. In DIP and its followup works, on the other hand, the reconstruction network is applied to discretized 2D or 3D images. Moreover, our goal is to use the representation power of deep networks to parameterize the solution of an IP in continuous space, whereas the main interests in DIP are on its intriguing automatic regularization properties.
}

\section{Weak Adversarial Network for Inverse Problems}
\label{sec:iwan}
The proposed weak adversarial network approach for IPs is inspired by the weak formulation of PDEs.
To obtain the weak formulation of the PDE in \eqref{eq:ip1a}, we multiply both sides of \eqref{eq:ip1a} by an arbitrary test function $\varphi \in H^{1}_{0}(\Omega)$ (the Hilbert space of functions with bounded first-order weak derivatives and compactly supported in $\Omega$) and integrate over $\Omega$:
\begin{equation}\label{eq:weak_u}
\langle \Acal[u,\gamma], \varphi \rangle \coloneqq \int_\Omega \Acal[u,\gamma](x)\varphi(x) \dif x = 0.
\end{equation}
One of the main advantages of weak formulation \eqref{eq:weak_u} is that we can subsequently apply integration by parts to transfer certain gradient operator(s) in $\Acal[u,\gamma]$ to $\varphi$, such that the requirement on the regularity of $u$ (and $\gamma$ if applicable) can be reduced.
For example, in the case of inverse conductivity problem \eqref{eq:eit}, the integration by parts and the fact that $\varphi=0$ on $\partial \Omega$ together yield
\begin{equation}\label{eq:weak_u_eit}
\langle \Acal[u,\gamma], \varphi \rangle = \int_\Omega \del[1]{\gamma \nabla u \cdot \nabla \varphi - f\varphi } \dif x = 0,
\end{equation}
where $\gamma \nabla u$ is not necessarily differentiable as in \eqref{eq:eit} in the classical sense anymore (we use $\nabla$ to denote the gradient operator with respect to $x$, and $\nabla_{\theta}$ as the gradient with respect to $\theta$ and so on in this paper).
We call $(u,\gamma) \in H^1(\Omega) \times L^2(\Omega)$ a \textit{weak solution} (or \textit{generalized solution}) of the inverse problem \eqref{eq:ip} if $(u,\gamma)$ satisfies the boundary condition \eqref{eq:ip1b} and \eqref{eq:weak_u} for all $\varphi \in H_0^1(\Omega)$.
Here $L^2(\Omega)$ is the Lebesgue space of square integrable functions on $\Omega$, and $H^1(\Omega) \subset L^2(\Omega)$ is the Hilbert space of functions with bounded first-order weak derivatives.
Note that any classical (strong) solution of \eqref{eq:ip} is also a weak solution.
In this work, we seek for weak solutions of inverse problem \eqref{eq:ip} so that we may be able to provide an answer to the problem even if it does not admit a solution in the classical sense.

Following the work \cite{zang2019weak}, we consider the weak formulation of the PDE $\Acal[u,\gamma]=0$ in \eqref{eq:ip}. To cope with the unknown solution $u$ and parameter $\gamma$ of the PDE in an inverse problem, we parameterize both $u$ and $\gamma$ as deep neural networks, and consider $\Acal[u,\gamma]: H_0^1(\Omega) \to \mathbb{R}$ as a linear functional such that $\Acal[u,\gamma](\varphi) \coloneqq \langle \Acal[u,\gamma], \varphi \rangle$ as defined in \eqref{eq:weak_u}. We define the norm of $\Acal[u,\gamma]$ induced by the $H_1$  norm as
\begin{equation}
\label{eq:norm_op}
 \|\Acal[u,\gamma]\|_{op} \coloneqq \sup_{\varphi\in H^1_0, \varphi\ne0}\frac{\langle\Acal[u,\gamma],\varphi\rangle}{\|\varphi\|_{H^1}},
\end{equation}
{where the $H^1$-norm of $\varphi$ is given by $\|\varphi\|_{H^1(\Omega)}^2=\int_\Omega (|\varphi(x)|^2 + |\nabla \varphi(x)|^2)\dif x$}.
Therefore, $(u,\gamma)$ is a weak solution of \eqref{eq:ip} if and only if $\|\Acal[u,\gamma]\|_{op} = 0$ and $\Bcal[u,\gamma]=0$ on $\partial \Omega$.
As $\|\Acal[u,\gamma] \|_{op}\ge 0$, we know that a weak solution $(u,\gamma)$ to \eqref{eq:ip} thus solves the following problem in observation of \eqref{eq:norm_op}:
\begin{equation}
\label{eq:min_op}
\minimize_{u,\gamma}\ \|\Acal[u,\gamma]\|_{op}^2 = \minimize_{u,\gamma}\ \sup_{\varphi \in H_0^1, \varphi \ne 0} \frac{|\langle\Acal[u,\gamma],\varphi \rangle|^2}{\|\varphi\|_{H^1}^2},
\end{equation}
among all $(u,\gamma) \in H^1(\Omega) \times L^2(\Omega)$, and attains minimal value $0$.
This result is summarized in the following theorem, and the proof is provided in Appendix \ref{pf:thm:min_op}.
\begin{theorem}
\label{thm:min_op}
Suppose that $(u^{*},\gamma^*)$ satisfies the boundary condition $\mathcal{B}[u^*,\gamma^*]=0$, then $(u^*,\gamma^*)$ is a weak solution of \eqref{eq:ip} if and only if  $\|\mathcal{A}[u^*,\gamma^*]\|_{op}=0$.
\end{theorem}

Theorem \ref{thm:min_op} implies that, to find the weak solution of \eqref{eq:ip}, we can instead seek for the optimal solution $(u^*,\gamma^*)$ that satisfies $\Bcal[u^*,\gamma^*]=0$ and meanwhile minimizes \eqref{eq:min_op} by achieving minimum operator norm value $\|\Acal[u^*,\gamma^*]\|_{op}=0$ due to the nonnegativity of the operator norm.
In other words, $(u^*,\gamma^*)$ is a weak solution of the problem \eqref{eq:ip} if and only if both $\|\Acal[u^*,\gamma^*]\|_{op}$ and $\|\Bcal[u^*,\gamma^*]\|_{L^2(\partial\Omega)}$ vanish.
Therefore, we can solve $(u^*,\gamma^*)$ from the following minimization problem which is equivalent to \eqref{eq:ip}:
\begin{equation}\label{eq:obj}
\minimize_{u,\gamma}\ I(u,\gamma) = \|\Acal[u,\gamma]\|_{op}^2 + \beta \|\Bcal[u,\gamma]\|_{L^2(\partial\Omega)}^2,
\end{equation}
and $\beta>0$ is a weight parameter that balances the two terms in the objective function $I(u,\gamma)$.
Note that both terms of the objective function in \eqref{eq:obj} are nonnegative and vanish simultaneously only at a weak solution $(u^*,\gamma^*)$ of \eqref{eq:ip}.

A promising alternative to classical numerical methods for high-dimensional PDEs is the use of deep neural networks since they do not require domain discretization and are completely mesh free.
Deep neural networks are compositions of multiple simple functions (called layers) so that they can approximate rather complicated functions.
Consider a simple multi-layer neural network $u_{\theta}$ as follows:
\begin{equation}\label{eq:u_param}
u_\theta(x) = w_{K}^{\top}\, l_{K-1}\circ \cdots \circ l_{0}(x) + b_{K},
\end{equation}
where the $k$th layer $l_k: \mathbb{R}^{d_k} \to \mathbb{R}^{d_{k+1}}$ is given by $l_k(z) = \sigma_k (W_k z + b_k)$ with weight $W_k\in \mathbb{R}^{d_{k+1} \times d_k}$ and bias $b_k \in \mathbb{R}^{d_{k+1}}$ for $k=0,1,\dots,K-1$, and the network parameters of all layers are collectively denoted by $\theta$ as follows,
\begin{equation}
\label{eq:theta}
\theta \coloneqq (w_K, b_K, W_{K-1}, b_{K-1}, \dots, W_0, b_0).
\end{equation}
Throughout, all vectors in this paper are column vectors by default.
In \eqref{eq:u_param}, $x\in\Omega$ is the input of the network, $d_0=d$ is the problem dimension of \eqref{eq:ip} (also known as the size of input layer), $w_K \in \mathbb{R}^{d_K}$ and $b_K\in\mathbb{R}$ are parameters in the last $K$th layer (also called the output layer).
Typical choices of the nonlinear activation function $\sigma_k$ include sigmoid function $\sigma(z) = (1+e^{-z})^{-1}$, hyperbolic tangent (tanh) function $\sigma(z)=(e^{z} - e^{-z})/(e^{z} + e^{-z})$, and rectified linear unit (ReLU) function $\sigma(z)=\max(0,z)$, which are applied componentwisely.
The training of deep neural networks refers to the process of optimizing $\theta$ using available data or constraints such that the function $u_{\theta}$ can approximate the (unknown) target function.
More details about deep neural networks can be found in \cite{goodfellow2016deep}.

Despite of the simple structures like \eqref{eq:u_param}, deep neural networks are capable to approximate rather complicated continuous function (and its derivatives if needed) uniformly on a compact support $\bar{\Omega}$.
This significant result is known as the universal approximation theorem \cite{hornik1991approximation}.
%
%
%
The expressive power of neural networks ensured by the universal approximation theorem suggests a promising mesh-free parameterization of the weak solution $(u,\gamma)$ of \eqref{eq:ip}.
In what follows, we select sufficiently deep neural network structures of form \eqref{eq:u_param} for both $u$ and $\gamma$.
Specific structures, i.e., layer number $K$ and sizes $\{d_1,\dots,d_{K-1}\}$, used in our numerical experiments will be provided in Section \ref{sec:experiment}.
Note that $u$ and $\gamma$ are two separate networks, but we use a single letter $\theta$ to denote their network parameters rather than $\theta_u$ and $\theta_\gamma$ to simplify notations.
That is, we parameterize $(u,\gamma)$ as deep neural networks $(u_\theta,\gamma_\theta)$, and attempt to find the parameter $\theta$ such that $(u_\theta,\gamma_\theta)$ solves \eqref{eq:obj}.
To this end, the test function $\varphi$ in the weak formulation \eqref{eq:weak_u} is also parameterized as a deep neural network $\varphi_\eta$ in a similar form of \eqref{eq:u_param} and \eqref{eq:theta} with parameter denoted by $\eta$.
With the parameterized $(u_\theta,\gamma_\theta)$ and $\varphi_{\eta}$, we follow the inner product notation in \eqref{eq:weak_u} and define
\begin{equation}
E(\theta,\eta) \coloneqq |\langle \Acal[u_\theta, \gamma_\theta], \varphi_\eta \rangle|^2.
\end{equation}
Instead of normalizing $E(\theta,\eta)$ by $\|\varphi_\eta\|_{H_1}^2$ as in the original definition of (squared) operator norm \eqref{eq:norm_op}, we approximate (up to a constant scaling of) the squared operator norm in \eqref{eq:norm_op} by the following max-type function of $\theta$:
\begin{equation}
\label{eq:Lint}
\Lint(\theta) \coloneqq \max_{|\eta|^2 \le 2B} E(\theta,\eta)
\end{equation}
where $B>0$ is a prescribed bound to constrain the magnitude of network parameter $\eta$.
Here $|\eta|^2 = \sum_{k}(\sum_{ij}[W_k]_{ij}^2 + \sum_{i}[b_k]_i^2)$, and $[M]_{ij}\in\mathbb{R}$ stands for the $(i,j)$th entry of a matrix $M$, and $[v]_i \in \mathbb{R}$ the $i$th component of a vector $v$.
It is worth noting that the bound constraint on the $\ell_2$-norm of $\eta$ in \eqref{eq:Lint} is similar to the weight clipping (equivalent to bound on $\ell_{\infty}$-norm) method used in WGAN \cite{arjovsky2017wasserstein}. However, they serve different purposes: the constraint in \eqref{eq:Lint} is introduced so that the integrals, such as \eqref{eq:weak_u_eit}, are bounded (the actual value of this bound can be arbitrary). In this case, the stochastic gradients obtained by the Monte-Carlo approximations in our numerical implementation have bounded variance, which is needed in the proof of Theorem \ref{thm:converge} below. On the other hand, the weight clipping in WGAN is to ensure the dual function realized by the neural network is in the class of $1$-Lipschitz functions $\mathcal{F}:=\{f: \Omega \to \mathbb{R}: |f(x)-f(y)| \le |x-y|,\ \forall\, x,y \in \Omega\}$. As noted in \cite{arjovsky2017wasserstein}, weight clipping is a simple but not appropriate way to implement the $1$-Lipschitz constraint, and hence there is a series of followup work to tackle this issue, such as \cite{gulrajani2017improved,miyato2018spectral}.
%

Furthermore, we define the loss function associated with the boundary condition \eqref{eq:ip1b} by
\begin{equation}
\label{eq:Lbdry}
\Lbdry(\theta) \coloneqq \|\Bcal[u_\theta,\gamma_\theta]\|_{L^2(\partial \Omega)}^2 = \int_{\partial\Omega}|\Bcal[u_\theta,\gamma_\theta](x)|^2 \dif S(x).
\end{equation}
For instance, if the boundary condition of $(u,\gamma)$ is given in \eqref{eq:ip2b} with known boundary value $(u_b,\gamma_b,u_n)$, then $\Lbdry(\theta)=\int_{\partial\Omega} |u_\theta(x)-u_b(x)|^2 + |\gamma_\theta(x)-\gamma_b(x)|^2 + |\partial_{\vec{n}(x)}u(x) - u_n(x)|^2 \dif S(x)$.
Finally, we define the total loss function $L(\theta)$, and solve the following minimization problem of its optimal $\theta^*$:
\begin{equation}\label{eq:L}
\minimize_{\theta}\ L(\theta),\quad \text{where}\quad L(\theta) \coloneqq \Lint(\theta) + \beta \Lbdry(\theta),
\end{equation}
where we also constrain on the magnitude of the parameter $\theta$ such that $|\theta|^2 \le 2B$ for the same $B$ to simplify notation.
Note that here both $\theta$ and $\eta$ are finite dimensional vectors, and $\Lint(\theta), \Lbdry(\theta), E(\theta,\eta)\in \mathbb{R}_+$, hence it is possible to apply numerical optimization algorithms to find the minimizer of $L(\theta)$.

A standard approach to solving a minimization problem like \eqref{eq:L} is the projected gradient descent method which performs the following iteration:
\begin{equation}
\label{eq:gd}
\theta \leftarrow \Pi(\theta - \tau \nabla_{\theta} L(\theta)),
\end{equation}
where $\Pi(\theta) = \min(\sqrt{2B},|\theta|)\cdot (\theta/|\theta|)$ is the projection of $\theta$ to the ball centered at origin with radius $\sqrt{2B}$, and $\tau>0$ is the step size.
As we can see, the main computation of \eqref{eq:gd} is on the gradient $\nabla_\theta L(\theta) = \nabla_{\theta} \Lint(\theta) + \beta \nabla_{\theta} \Lbdry(\theta)$.
The computation of $\nabla_{\theta} \Lbdry(\theta)$ is straightforward as shown later.
The loss $\Lint(\theta)$, however, is defined as a maximization problem \eqref{eq:Lint}, and we need to write its gradient as a function of $\theta$ first.
To this end, we have the following lemma to compute the gradient $\nabla_{\theta} \Lint(\theta)$, and the proof is provided in Appendix \ref{pf:lem:dLint}.
\begin{lemma}
\label{lem:dLint}
Suppose $\Lint(\theta)$ is defined in \eqref{eq:Lint}. Then the gradient $\nabla_\theta\Lint(\theta)$ at any $\theta$ is given by $\nabla_{\theta} \Lint(\theta) = \partial_{\theta} E(\theta, \eta(\theta))$, where $\eta(\theta)$ is a solution of $\max_{|\eta|^2 \le 2B} E(\theta,\eta)$ for the specified $\theta$.
\end{lemma}
\begin{remark}
Lemma \ref{lem:dLint} suggests that, to obtain $\nabla_\theta \Lint(\theta)$ at any given $\theta$, we can first take the partial derivative of $E$ with respect to $\theta$ with $\eta$ untouched, and then evaluate the partial derivative using $\theta$ and any solution $\eta(\theta)$ of the maximization problem \eqref{eq:Lint}.
\end{remark}

The exact gradients of $\Lint(\theta)$ and $\Lbdry(\theta)$ require integrations of functions parameterized by deep neural networks over $\Omega$ and $\partial\Omega$ in continuous space, which are computationally intractable in practice.
Therefore, we use Monte-Carlo (MC) approximations of these integrals.
To this end, we need the following result on the approximation of integrals using samples, and the proof is provided in Appendix \ref{pf:lem:integral}.
\begin{lemma}
\label{lem:integral}
Suppose $\Omega \subset \mathbb{R}^d$ is bounded, and $\rho$ is a probability density defined on $\Omega$ such that $\rho(x)>0$ for all $x\in\Omega$.
Given a function $\psi \in L^2(\Omega)$, denote $\Psi = \int_{\Omega} \psi(x) \dif x$.
Let $x^{(1)},\dots,x^{(N)}$ be $N$ independent samples drawn from $\rho$.
Consider the following estimator $\hat{\Psi}$ of $\Psi$:
\begin{equation}
\label{eq:Psi}
\hat{\Psi} = \frac{1}{N}\sum_{i=1}^N\frac{\psi(x^{(i)})}{\rho(x^{(i)})}.
\end{equation}
Then the first and second moments of $\hat{\Psi}$ are given by
\begin{equation}
\label{eq:Psi_moment}
\mathbb{E}[\hat{\Psi}]=\Psi\quad\text{and}\quad
\mathbb{E}[\hat{\Psi}^2]= \frac{N-1}{N}\Psi^2 + \frac{1}{N} \int_{\Omega} \frac{\psi(x)^2}{\rho(x)} \dif x .
\end{equation}
Hence the variance of $\hat{\Psi}$ is $N^{-1}\cdot (\int_{\Omega} (\psi^2/\rho) \dif x - (\int_\Omega \psi \dif x)^2)$.
In particular, with the uniform distribution $\rho(x)=1/|\Omega|$, the variance of $\hat{\Psi}=(|\Omega|/N)\cdot \sum_{i}\psi(x^{(i)})$ is $N^{-1}\cdot (|\Omega|\int_{\Omega} \psi^2 \dif x - (\int_{\Omega} \psi \dif x)^2)$.
\end{lemma}
\begin{remark}
We have several remarks regarding Lemma \ref{lem:integral}:
\begin{itemize}
\item
The estimator $\hat{\Psi}$ of the integral $\Psi$ is unbiased.
\item
The variance of $\hat{\Psi}$ shown above decreases at the rate of $O(1/N)$ in the number $N$ of sample collocation points. By H\"{o}lder's inequality and that $\rho$ is a probability density, we know
\begin{equation*}
\abs[2]{\int_\Omega \psi \dif x} \le \int_\Omega |\psi| \dif x = \int_\Omega \frac{|\psi|}{\sqrt{\rho}}\sqrt{\rho} \dif x \le \del[2]{\int_\Omega \frac{|\psi|^2}{\rho} \dif x}^{1/2}\del[2]{\int_\Omega \rho \dif x}^{1/2} = \del[2]{\int_\Omega \frac{|\psi|^2}{\rho} \dif x}^{1/2},
\end{equation*}
which also verifies that $\mathrm{V}(\hat{\Psi})\ge 0$.
More importantly, the equalities hold if $\psi$ does not change sign and $\rho\propto |\psi|$.
Therefore, we can set $\rho$ as close to $|\psi|$ (up to a normalizing constant) as possible to reduce the variance, but meanwhile ensure $\rho$ is easy to sample from and evaluate as required in \eqref{eq:Psi}. This is closely related to the concept of importance sampling.
\item
The result \eqref{eq:Psi} and \eqref{eq:Psi_moment} in Lemma \ref{lem:integral} can be easily extend to the case with unbounded domain $\Omega$, provided that $\psi/\sqrt{\rho}\in L^2(\Omega)$.
\end{itemize}
\end{remark}

Lemma \ref{lem:integral} provides a feasible way to approximate the gradient of $L(\theta)$ for \eqref{eq:gd}.
For instance, to compute $\nabla_\theta \Lbdry(\theta)$, we can take gradient of \eqref{eq:Lbdry} 
with respect to $\theta$, sample $N_b$ collocation points $\{x_b^{(i)}:1\le i\le N_b\}$ on the boundary $\partial \Omega$
and approximate $\nabla_\theta \Lbdry(\theta)$ by
summation of function evaluations at the sample points. If we take $\Bcal[u,\gamma] = (u-u_b,\gamma-\gamma_b,\partial_{\vec{n}}u - u_n)$ and uniformly sample $x_b^{(i)}$, the estimate becomes
\begin{align}
\label{eq:dLbdry_est}
\nabla_{\theta} \Lbdry(\theta)
= &\ 2 \int_{\partial\Omega} \del[2]{(u_\theta-u_b)\nabla_\theta u_\theta + (\gamma_\theta-\gamma_b)\nabla_\theta \gamma_\theta + (\partial_{\vec{n}} u_{\theta} - u_n)\nabla_{\theta}\nabla u \cdot \vec{n} } \dif S(x) \nonumber \\
\approx &\ \frac{2|\partial\Omega|}{N_b}\sum_{i=1}^{N_b}\Big( (u_\theta(x_b^{(i)})-u_b(x_b^{(i)}))\nabla_\theta u_\theta(x_b^{(i)}) + (\gamma_\theta(x_b^{(i)})-\gamma_b(x_b^{(i)}))\nabla_\theta \gamma_\theta(x_b^{(i)}) \\
&\qquad \qquad \qquad + (\partial_{\vec{n}} u_\theta(x_b^{(i)})-u_b(x_b^{(i)}))\nabla_\theta \nabla u_{\theta} x_b^{(i)}\cdot \vec{n}x_b^{(i)}  \Big). \nonumber
\end{align}
Similarly, we can compute the stochastic gradient of  $\nabla_{\theta}\Lint(\theta)$. In the case of taking $\Acal[u,\gamma] = \nabla\cdot(\gamma\nabla u) - f$ in $\Omega$ with $f$ given, and uniformly sampling $N_r$ collocation points $\{x_r^{(i)}:1\le i\le N_r\}$ inside the region $\Omega$, $\nabla_\theta \Lint(\theta)$ can be estimated by
\begin{align}
\nabla_{\theta} \Lint(\theta)
= &\ 2 I(\theta)\int_{\Omega} \del[2]{\nabla_\theta \gamma_\theta( \nabla u_\theta \cdot \nabla \varphi_{\eta(\theta)}) + \gamma_\theta( \nabla_{\theta}\nabla u_\theta \cdot \nabla \varphi_{\eta(\theta)})} \dif S(x) \label{eq:dLint_est} \\
\approx &\ \frac{2|\Omega|\hat{I}(\theta)}{N_r}\sum_{i=1}^{N_r}\del[2]{\nabla_\theta \gamma_\theta(x_r^{(i)})( \nabla u_\theta(x_r^{(i)}) \nabla \varphi_{\eta(\theta)}(x_r^{(i)})) + \gamma_\theta(x_r^{(i)})( \nabla_{\theta}\nabla u_\theta(x_r^{(i)})\nabla \varphi_{\eta(\theta)}(x_r^{(i)}))} \nonumber
\end{align}
where $I(\theta)$ and its estimator $\hat{I}(\theta)$ are given by
\begin{align*}
I(\theta) = \int_{\Omega} \gamma_\theta (\nabla u_\theta \cdot \nabla \varphi_{\eta(\theta)}) \dif x, \quad
\hat{I}(\theta) = \frac{|\Omega|}{2} \sum_{i=1}^{N_r} \gamma_\theta(x_r^{(i)}) (\nabla u_\theta(x_r^{(i)}) \cdot \nabla \varphi_{\eta(\theta)}(x_r^{(i)})),
\end{align*}
and $\eta(\theta)$ is a solution of the maximization problem \eqref{eq:Lint} according to Lemma \ref{lem:dLint}.
All integrals in the gradients can be approximated in a similar way.
These approximated gradients are in fact stochastic gradients, which are unbiased and have bounded variances due to the boundedness of the network parameters.
With these approximations, \eqref{eq:gd} reduces to the stochastic projected gradient descent method, which ensures convergence to a local stationary point of \eqref{eq:L} with proper choice of step sizes.
Since \eqref{eq:L} is constrained, the \textit{gradient mapping}, defined by $\Gcal(\theta)\coloneqq \tau^{-1}[\theta-\Pi(\theta-\tau \nabla_{\theta} L(\theta))]$, is used as the convergence criterion of $\theta$ \cite{ghadimi2016mini,li2018simple,reddi2016proximal}.
Note that the definition of gradient mapping takes the normalization of step size $\tau$ into consideration.
Moreover, without the projection $\Pi$, the gradient mapping reduces to $\Gcal(\theta) = \nabla_{\theta}L(\theta)$, whose magnitude is an evaluation criterion for local stationary points (i.e., $|\nabla_{\theta} L(\theta)|=0$) for unconstrained case.
This result is stated in the following theorem, and the proof is given in Appendix \ref{pf:thm:converge}.
\begin{theorem}
\label{thm:converge}
For any $\varepsilon>0$, let $\{\theta_j\}$ be a sequence of the network parameter in $(u_\theta,\gamma_\theta)$ generated by the gradient descent algorithm \eqref{eq:gd} with integrals in $\nabla_\theta L(\theta)$ approximated by sample averages as in \eqref{eq:Psi} with sample complexities $N_r,N_b=O(\varepsilon^{-1})$ in each iteration, then $\min_{1\le j\le J}\mathbb{E}[|\Gcal(\theta_j)|^2] \le \varepsilon$ after $J=O(\varepsilon^{-1})$ iterations.
\end{theorem}
\begin{remark}
Theorem \ref{thm:converge} establishes the convergence and iteration complexity of \eqref{eq:gd} to the so-called $\varepsilon$-solution of the problem. The result is based on the expected magnitude of the gradient mapping, which is a standard convergence criterion in nonconvex constrained stochastic optimization.
However, this only ensures approximation to a stationary point (not necessarily a local or global minimizer) on expectation.
In theory, one can apply additional global optimization techniques to \eqref{eq:obj} in order to find a global minimizer (possibly only with high probability at best) with substantially higher computational cost.
However, we will not exploit this issue further in this work.
%
\end{remark}

Now we summarize the steps of our algorithm for solving IPs using weak adversarial networks.
To simplify the presentation, we introduce the following notation to indicate the stochastic gradient descent (SGD) procedure for finding a minimizer of a loss function $L(\theta)$:
\begin{equation}
\label{eq:sgd}
\theta^* \leftarrow \text{SGD}(G(\theta), X, \theta_0, \tau, J),
\end{equation}
which means the output $\theta^*$ is the result $\theta_J$ after we execute the (projected) SGD scheme with step size $\tau$ below for $j=0,\dots,J-1$ with initial $\theta_0$:
\begin{equation}
\label{eq:sgd_step}
\theta_{j+1} \leftarrow \Pi(\theta_j - \tau \hat{G}(\theta_j;X)).
\end{equation}
Here $X=\{x^{(i)}:1\le i\le N\}$ is the set of $N$ sampled collocation points, $G(\theta)\coloneqq \nabla_{\theta} L(\theta)$ is the gradient of the loss function $L(\theta)$ to be minimized, and $\hat{G}(\theta;X)$ stands for the stochastic approximation of $G(\theta)$ at any given $\theta$, where the integrals are estimated as in \eqref{eq:Psi} using the sampled collocation points $X$.
{Therefore, each iteration of our algorithm consists of two steps.
In Step 1, we fix $\theta$ and solve the maximization problem with objective function $E(\theta,\eta)$ defined in \eqref{eq:Lint} by applying stochastic gradient ascent for $J_\eta$ steps to obtain an approximate maximizer $\eta$;
In Step 2, we fix this $\eta$, and update $\theta$ by one stochastic gradient descent step using gradient $\nabla_{\theta} L(\theta)=\partial_\theta E(\theta,\eta) + \beta \nabla_{\theta} \Lbdry(\theta)$.
Then we go to Step 1 to start the next iteration.
Hence, our objective function is $E(\theta,\eta) + \beta \Lbdry(\theta)$, for which we seek for the optimal point $(\theta^*,\eta^*)$ via a min-max optimization $\min_{\theta}\max_{\eta} E(\theta,\eta) + \beta \Lbdry(\theta)$.}
This procedure is referred to Inverse Problem Solver using Weak Adversarial Network (IWAN) and summarized in Algorithm \ref{alg:iwan}. The parameter values in our numerical implementations are presented in Section \ref{sec:experiment}.
\begin{algorithm}[t]
\caption{\underline{I}nverse Problem Solver by \underline{W}eak \underline{A}dversarial \underline{N}etwork (IWAN)}
\label{alg:iwan}
\begin{algorithmic}
\STATE \textbf{Input:} {The domain $\Omega$ and data for the Inverse Problem \eqref{eq:ip}.}
\STATE \textbf{Initialize:} {$(u_\theta,\gamma_\theta)$, $\varphi_{\eta}$.}
\FOR{$j = 1, \dots, J$:}
\STATE Sample $X_r=\{x_r^{(i)}: 1\le i\le N_r\}\subset \Omega$ and $X_b=\{x_b^{(i)}: 1\le i\le N_b\}\subset \partial \Omega$.
\STATE $\eta \leftarrow \text{SGD}(-\nabla_\eta E(\theta,\eta), X_r, \eta, \tau_{\eta}, J_{\eta})$.
\STATE $\theta \leftarrow \text{SGD}(\partial_\theta E(\theta,\eta) + \beta \nabla_{\theta} \Lbdry(\theta), (X_r, X_b), \theta, \tau_{\theta}, 1)$.
\ENDFOR
\STATE \textbf{Output:} {$(u_{\theta},\gamma_{\theta})$.}
\end{algorithmic}
\end{algorithm}

\section{Numerical Experiments}
\label{sec:experiment}
\subsection{Implementation Details}
\label{subsec:implementation}
In this subsection, we discuss several implementation details and modifications regarding Algorithm \ref{alg:iwan}.
First, to avoid spending excessive time in solving the inner maximization problem $\max_{\eta} E(\theta,\eta)$ in \eqref{eq:Lint} with a fixed $\theta$, we only apply a few iterations $J_\eta$ to compute $\eta$. Then we switch to update  $\theta$ for one iteration. See the two SGD steps in Algorithm \ref{alg:iwan}.
{This can improve overall efficiency and avoid spending excessive time on the inner maximization problem of $\eta$, especially when $\theta$ is still far from optimal yet.}
In fact, we can employ two separate test functions $\varphi_\eta$ and $\bar{\varphi}_\eta$ (we again use the same $\eta$ for notation simplicity).
In each iteration $j$, we alternately update $(u_\theta,\varphi_\eta,\gamma_\theta,\bar{\varphi}_{\eta})$ in order, each with one or a few SGD steps \eqref{eq:sgd_step}.
We will specify the numbers of steps for these networks for our experiments below.

During the derivations in Section \ref{sec:iwan}, we require bounded network parameters $\theta$ and $\eta$, where the bound $B$ can be arbitrarily large, to ensure finite variances of the integral estimators using samples so that the SGD is guaranteed to converge.
An alternative way to handle the boundedness constraints is to add $|\theta|^2$ and $|\eta|^2$ as regularization terms to the objective function in \eqref{eq:obj}.
{One can also use the operator norm \eqref{eq:norm_op} with denominator replaced by $\|\varphi\|_2^2:=\int_{\Omega} |\varphi|^2 \dif x$ (approximated by MC similarly as in \eqref{eq:Psi}), which is also adopted in our implementation.
This replacement does not cause issue in numerical implementation since the test function $\varphi_{\eta}$ is realized by a network with fixed width/depth and bounded parameters, and hence is guaranteed to be in $H^1$.
}

A test function $\varphi_{\eta}$ is required to vanish on $\partial\Omega$ in the weak formulation \eqref{eq:weak_u}.
One simple technique to ensure this is to precompute a function $\varphi_0\in C(\Omega)$ such that $\varphi_0(x)=0$ if $x\in\partial\Omega$ and $\varphi_0(x)>0$ if $x\in\Omega$ (e.g., a distance function to $\partial \Omega$ would work).
Then we seek for a parameterized network $\varphi_\eta'$ with no constraint on its arbitrary boundary, and set the test function $\varphi_\eta$ to $\varphi_0\varphi_\eta'$ which still takes zero value on $\partial\Omega$.

We implemented our algorithm using TensorFlow \cite{abadi2016tensorflow} (Python version 3.7), a state-of-the-art deep learning package that can efficiently employ GPUs for parallel computing.
The gradients with respect to network parameters ($\theta$ and $\eta$) and input ($x$) are computed by the TensorFlow builtin auto-differentiation module.
During training, we can also substitute the standard SGD optimizer by many of its variants, such as AdaGrad, RMSProp, Adam, Nadam etc.
In our experiments, we use AdaGrad supplied by the TensorFlow package, which appears to provide better performance than other optimizers in most of our tests.
All other parameters, such as the network structures (numbers of layers and neurons), step sizes (also known as the learning rates), number of iterations, will be specified in Section \ref{sec:experiment}.

\subsection{Experiment Setup}
In this section, we conduct a set of numerical experiments to show the practical performance of Algorithm \ref{alg:iwan} in solving inverse problems.
To quantitatively evaluate the accuracy of an approximate solution $\gamma$, we use the \textit{relative error} (in the $L^2$ sense) of $\gamma$ to the ground truth $\gamma^*$, defined by $\|\gamma-\gamma^*\|_2/\|\gamma^*\|_2$, where $\|\gamma\|_2^2 \coloneqq \int_{\Omega} |\gamma(x)|^2 \dif x$.
In practice, we compute $\|\gamma\|_2^2 = (|\Omega|/N)\cdot\sum_{i=1}^N |\gamma(x^{(i)})|^2$ by evaluating $\gamma$ on a fixed set of $N$ mesh grid points $\{x^{(i)}\in\Omega:1\le i\le N\}$ in $\Omega$.
More specifically, we used a regular mesh grid of size $100\times 100$ for $(x_1,x_2)$, and sampled one point $x$ for each of these grid points, i.e., for each grid point $(x_1,x_2)$, randomly draw values of the other coordinates within the domain $\Omega$ such that $N=10^4$. These points were sampled in advance and then used for all comparison algorithms to compute their test relative error. Note that these points are different from those sampled for training in these methods.
%


In all of our experiments, we parameterize each of $(u_\theta,\varphi_\eta,\gamma_\theta,\bar{\varphi}_{\eta})$ as a 9-layer fully connected neural network with 20 neurons per layer as in \eqref{eq:u_param} unless otherwise noted.
We set $\sigma_k$ to \texttt{tanh} for $k=1,2$, \texttt{softplus} for $k=4,6,8$, \texttt{sinc} for $k=3,5,7$ in $u_\theta$, and \texttt{tanh} for $k=1,2, 4, 6$, \texttt{elu} for $k=3, 5$, and \texttt{sigmoid} for $k=7, 8$ in $\gamma_{\theta}$. We use \texttt{elu} in the output layer of $\gamma_{\theta}$.
In parallel, we set $\sigma_k$ to \texttt{tanh} for $k=1,2$ and \texttt{sinc} for $k\ge3$ in $\varphi_\eta$ and $\bar{\varphi}_{\eta}$.
{Unless otherwise noted, we apply one SGD update with step size $\tau_\theta=0.01$ to both of $u_{\theta}$ and $\gamma_\theta$ (each of them performs the $\theta$ update in Algorithm \ref{alg:iwan}), and two ($J_\eta=2$) SGD updates with step size $\tau_{\eta}=0.008$ to both of $\varphi_{\eta}$ and $\bar{\varphi}_{\eta}$ (each of them performs the $\eta$ update in Algorithm \ref{alg:iwan}), following the order of $u_\theta, \varphi_\eta, \gamma_\theta,\bar{\varphi}_{\eta}$ in every iteration $j$ of Algorithm \ref{alg:iwan}.}
We set the weight $\beta=10,000$ for the boundary loss function $\Lbdry(\theta)$ in \eqref{eq:obj}, but also set a weight $\beta'$ to $\Lint(\theta)$ and specify its value in the experiment.
Other parameters will also be specified below.
All the experiments are implemented, trained, and tested in the TensorFlow framework \cite{abadi2016tensorflow} on a machine equipped with Intel 2.3GHz CPU and an Nvidia Tesla P100 GPU and 16GB of graphics card memory.

\subsection{Experimental Results on Inverse Conductivity Problems}
\label{subsec:test}
\begin{figure}
\centering
\begin{subfigure}[b]{.3\textwidth}
\includegraphics[width=\textwidth]{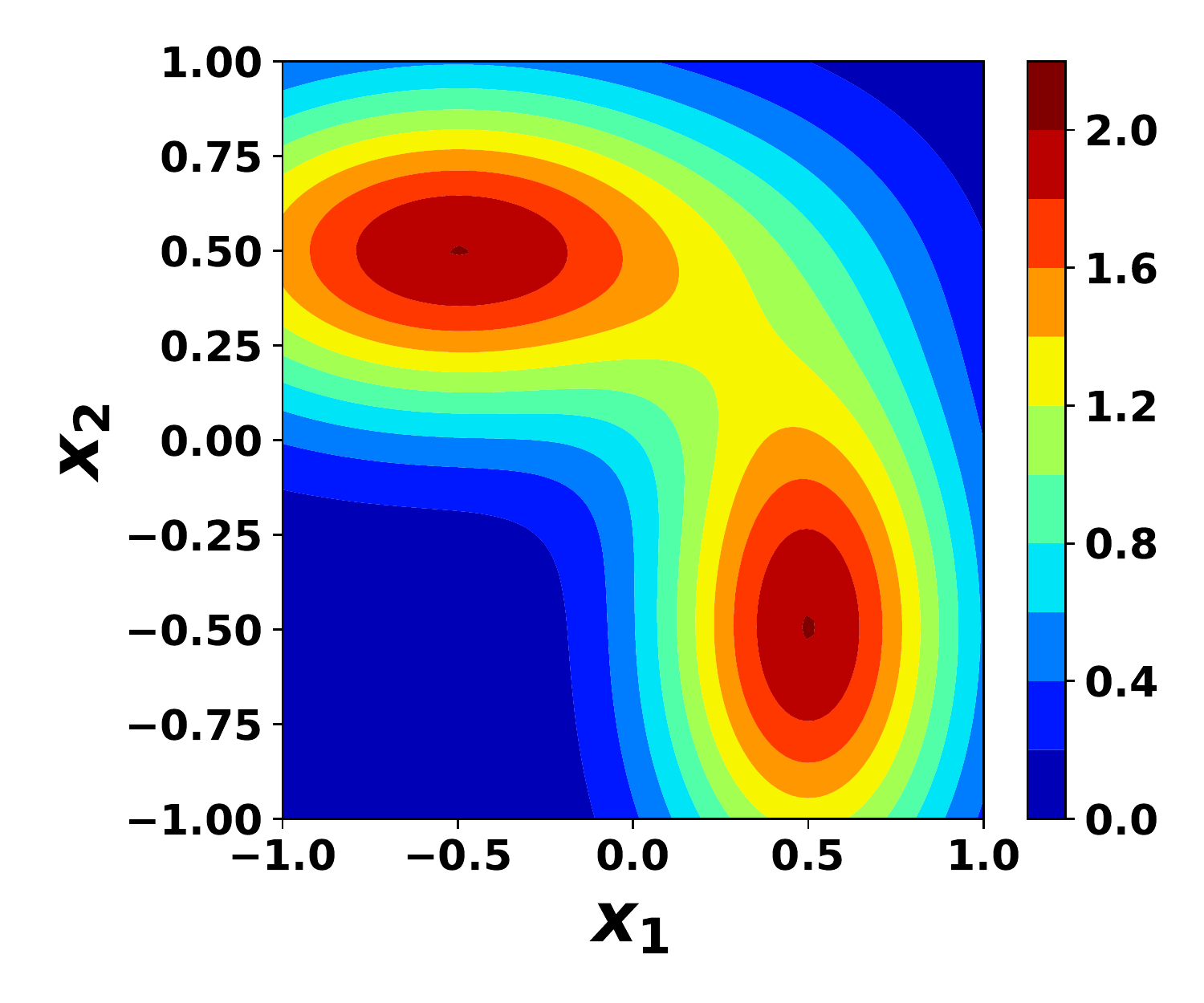}
\caption{Ground truth $\gamma^*$.}
\label{subfig:smooth_gamma}
\end{subfigure}
\begin{subfigure}[b]{.3\textwidth}
\includegraphics[width=\textwidth]{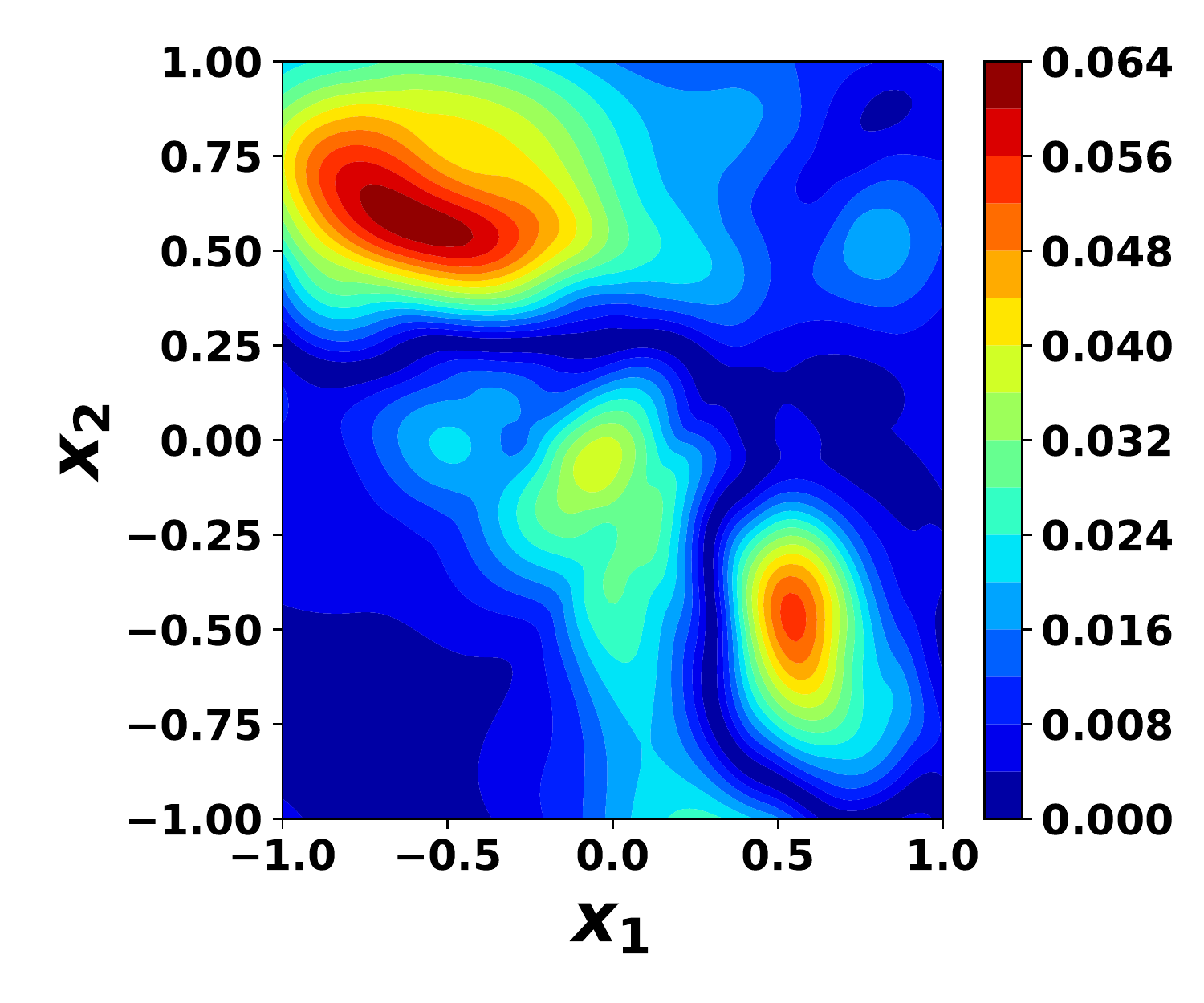}
\caption{Pointwise error $|\gamma^*-\gamma_{\theta}|$}
\label{subfig:smooth_abs_error}
\end{subfigure}
\begin{subfigure}[b]{.3\textwidth}
\includegraphics[width=\textwidth]{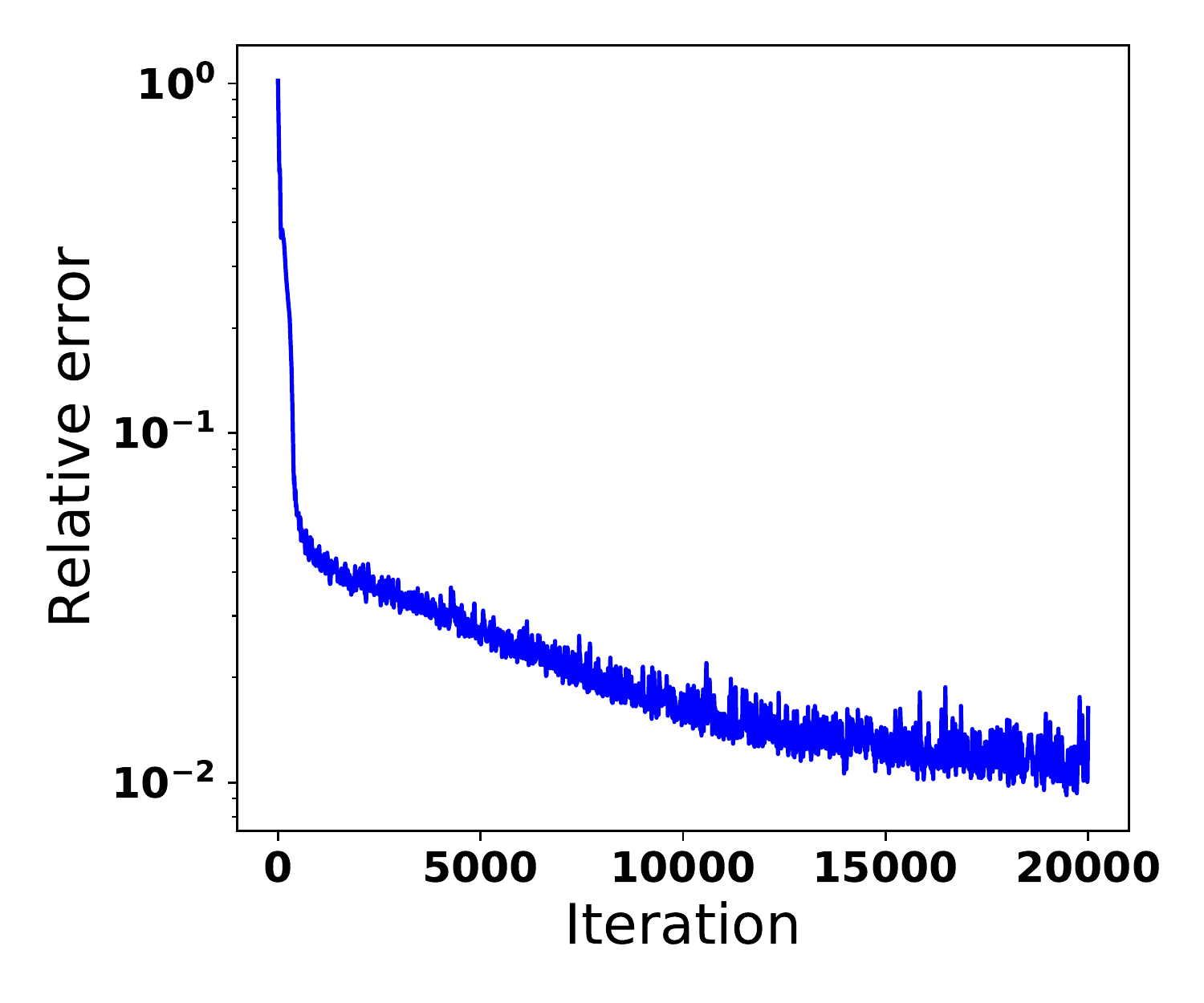}
\caption{Relative error of $\gamma_{\theta}$.}
\label{subfig:smooth_error}
\end{subfigure}
\caption{Test 1 result on \eqref{eq:eit} with smooth $\gamma^*$ and problem dimension $d=5$.}
\label{fig:continue}
\end{figure}

\textbf{Test 1: Inverse conductivity problem with smooth $\gamma$.}
We first test our method on the inverse conductivity problem \eqref{eq:eit} with a smooth conductivity distribution $\gamma$.
{In this test, we set $\Omega=(-1,1)^d\subset \mathbb{R}^d$ with problem dimension $d=5$. The setup for the ground truth conductivity distribution $\gamma^*$, the ground truth potential $u^*$ and the source term $f$ are provided in the table \ref{tab:prob_set} which is in the appendix \ref{app:prob_set}.}
%
%
%
We set $N_r= 100,000$ and $N_b= 100d$, $\beta'=10$, and run Algorithm \ref{alg:iwan} for $20,000$ iterations.
The true $\gamma^*$ and the point-wise error $|\gamma^*-\gamma_\theta|$ (with relative error $2.54\%$) are shown in Figure \ref{subfig:smooth_gamma} and \ref{subfig:smooth_abs_error} respectively.
The progress of the relative error of $\gamma_{\theta}$ versus iteration number is shown in Figure \ref{subfig:smooth_error} (all plots of relative error versus iteration number in this section are shown using moving average with a window size $7$).
%
For the demonstration purpose, only the $(x_1,x_2)$ cross sections that have main spatial variations are shown (same for the other test results below).

\begin{figure}[th]
\centering
\begin{subfigure}[b]{\textwidth}
\centering
\includegraphics[width=0.3\textwidth]{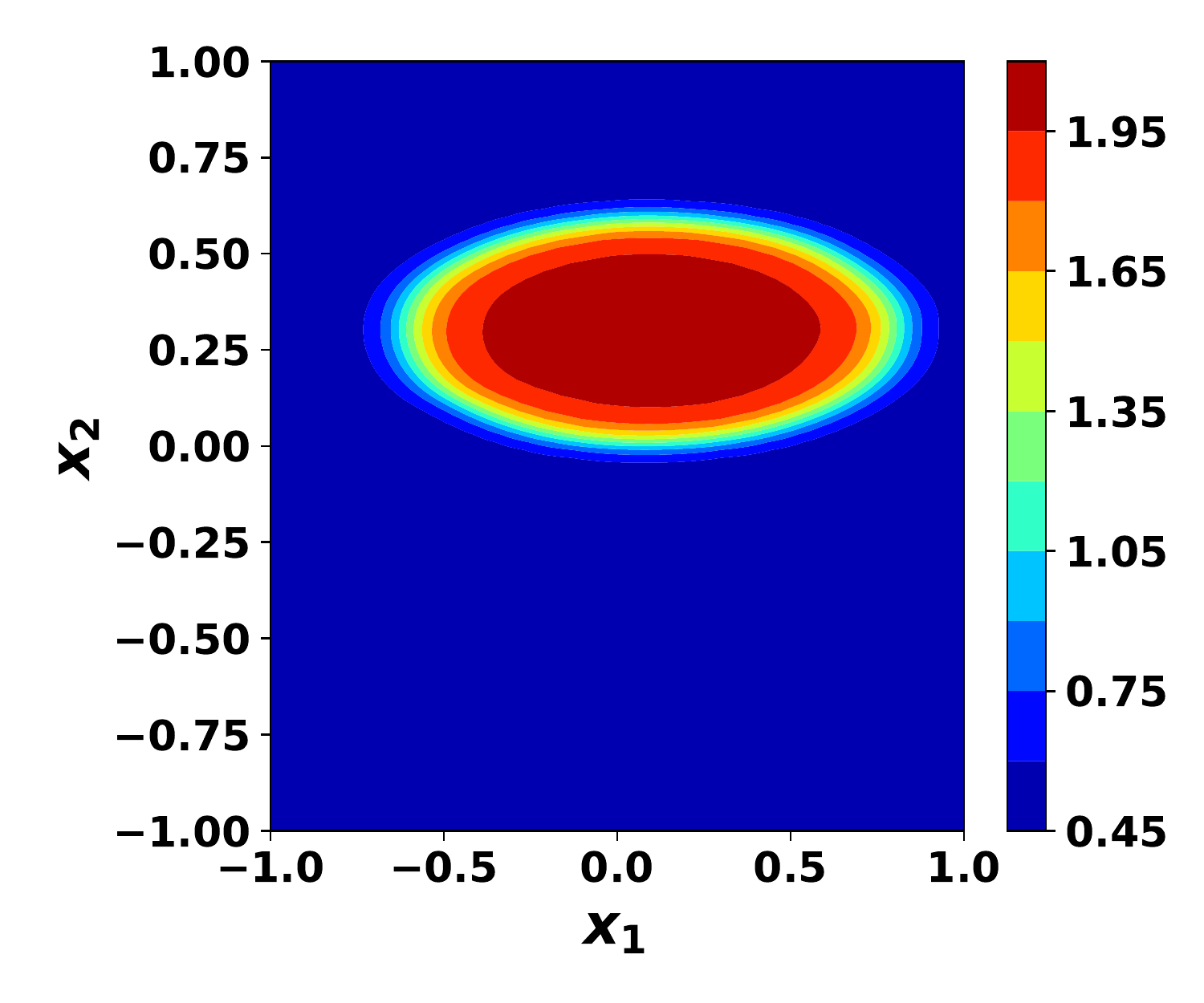}
\includegraphics[width=0.3\linewidth]{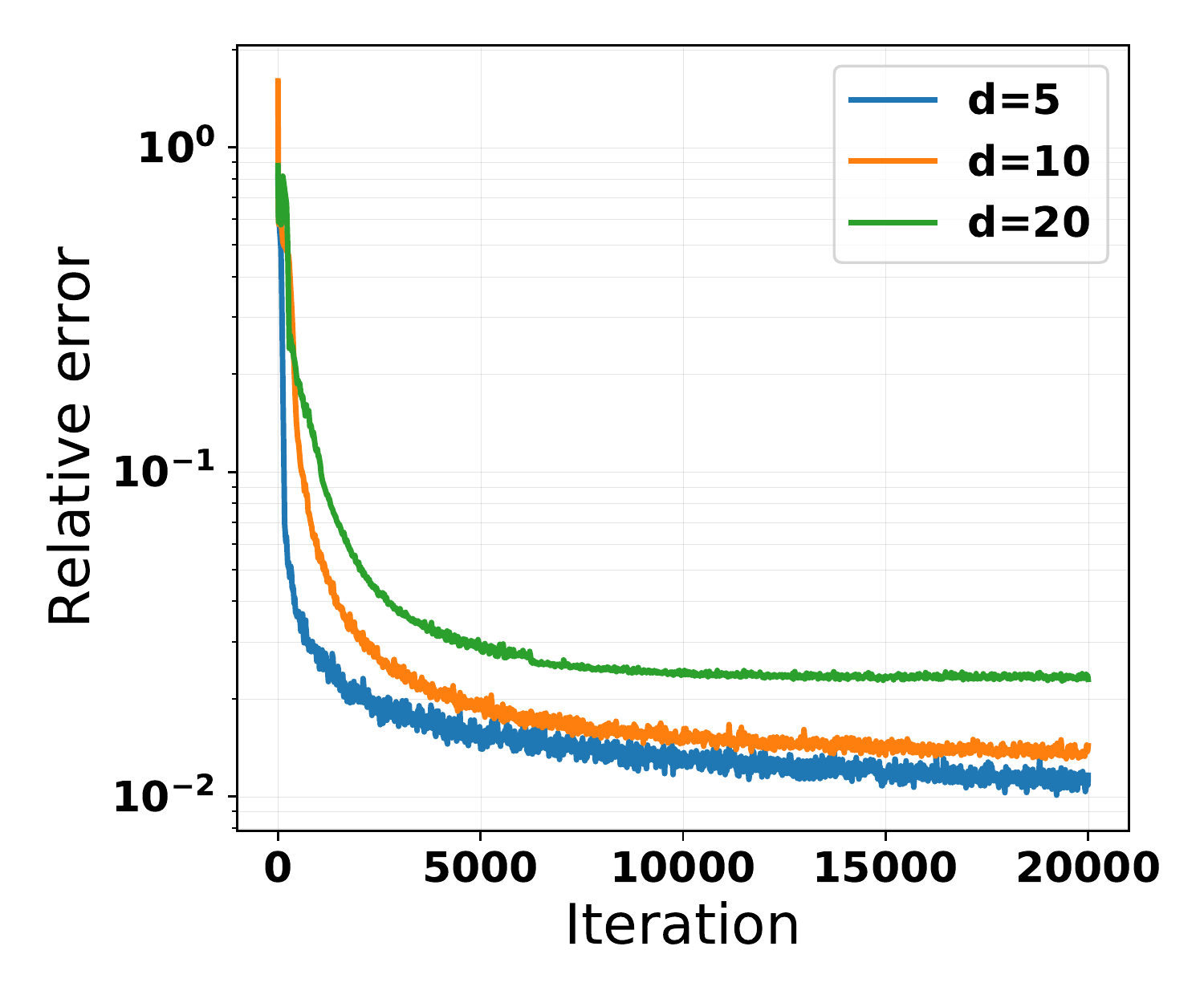}
\caption{Ground truth conductivity $\gamma^*$ (left) and relative error versus iteration (right).}
\label{subfig:pc_error}
\end{subfigure}
%
%
\begin{subfigure}[b]{\textwidth}
\includegraphics[width=0.9\linewidth]{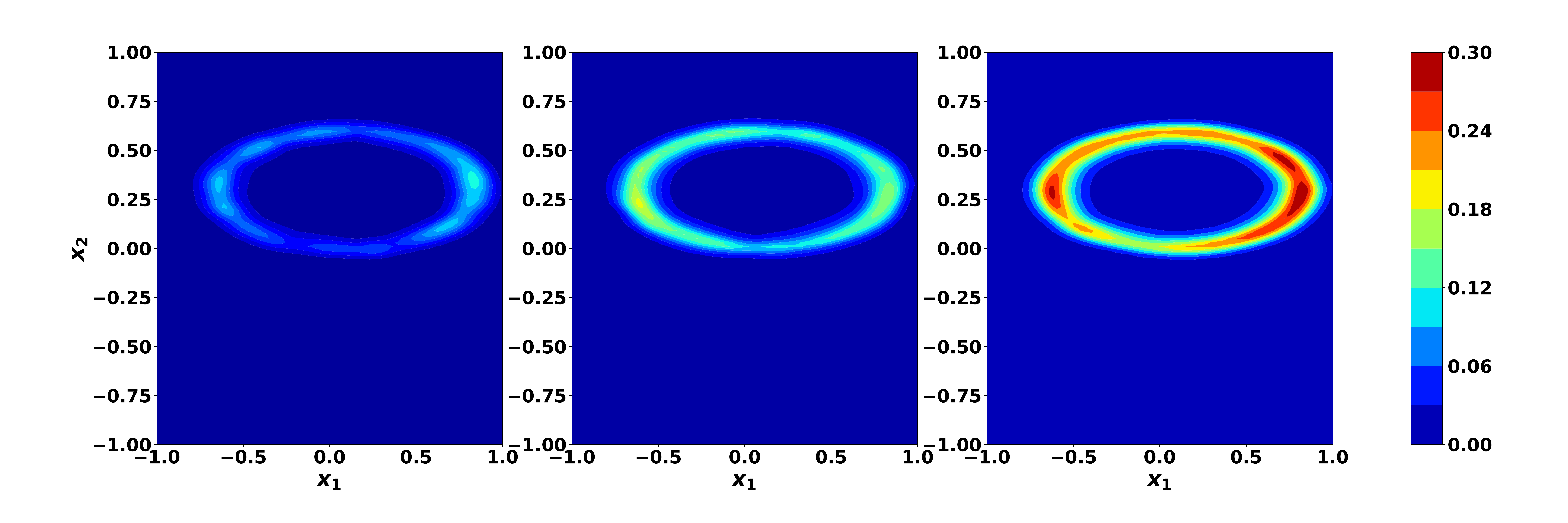}
\caption{From left to right: absolute error $|\gamma_\theta-\gamma^*|$ for $d=5, 10, 20$.}
\label{subfig:pc_abs_error}
\end{subfigure}
\caption{Test 2 result on \eqref{eq:eit} with nearly piecewise constant $\gamma^*$ and problem dimension $d=5,10,20$ without measurement noise.}
\label{fig:discrete_noise0}
\end{figure}

\medskip
\noindent
\textbf{Test 2: Inverse conductivity problem with nearly piecewise constant $\gamma$.}
We consider \eqref{eq:eit} with a less smooth, nearly piecewise constant conductivity $\gamma$.
In this test, we set $\Omega=(-1,1)^d$, define $\Omega_0=\{x\in\Omega:|x-c|_\Sigma^2 \le 0.6^2\}$ where $\Sigma=\mathrm{diag}(0.81,2,0.09,\dots,0.09)$ and $c=(0.1,0.3,0,\dots,0)$, and set $\gamma^*$ to $2$ in $\Omega_0$ and $0.5$ in $\Omega_0^c$.
{For ease of implementation, we slightly smooth the ground truth conductivity. One can find the setup of the smoothed conductivity $\gamma^*$, the ground truth potential $u^*$, and the source term $f$ in the table \ref{tab:prob_set}.}
%
%
%
%
%
%
Then we solve the inverse problem \eqref{eq:eit} with dimensionality $d=5, 10, 20$.
We set $N_r= 20,000d$ and $N_b= 100d$, and $\beta'=10, 1, 0.005$ for $d=5, 10, 20$ respectively.
In each case, we run Algorithm \ref{alg:iwan} for $20,000$ iterations, and obtain relative errors $1.16\%, 1.43\%, 2.29\%$ for $d=5, 10, 20$, respectively.
The recovery results are shown in Figure \ref{fig:discrete_noise0}.
Figure \ref{subfig:pc_error} shows the ground truth $\gamma^*$ (left) and the progress of relative errors versus iteration number for different $d$ (right).
The pointwise absolute errors $|\gamma_\theta-\gamma^*|$ for these dimensions are shown in Figure \ref{subfig:pc_error}.
%
%

\begin{figure}[th]
\centering
\begin{subfigure}[b]{\textwidth}
\centering
\includegraphics[width=0.3\linewidth]{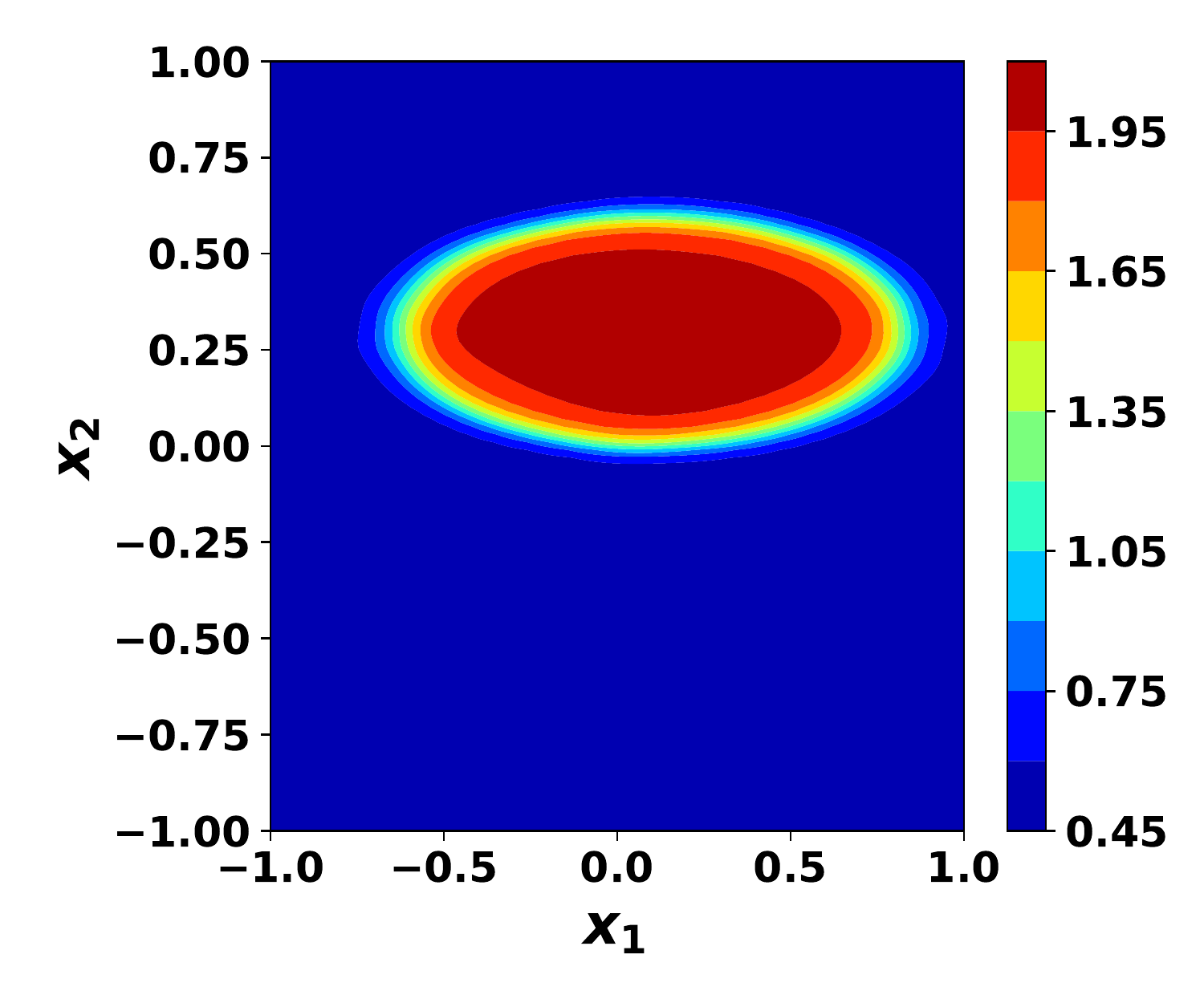}
\includegraphics[width=0.3\linewidth]{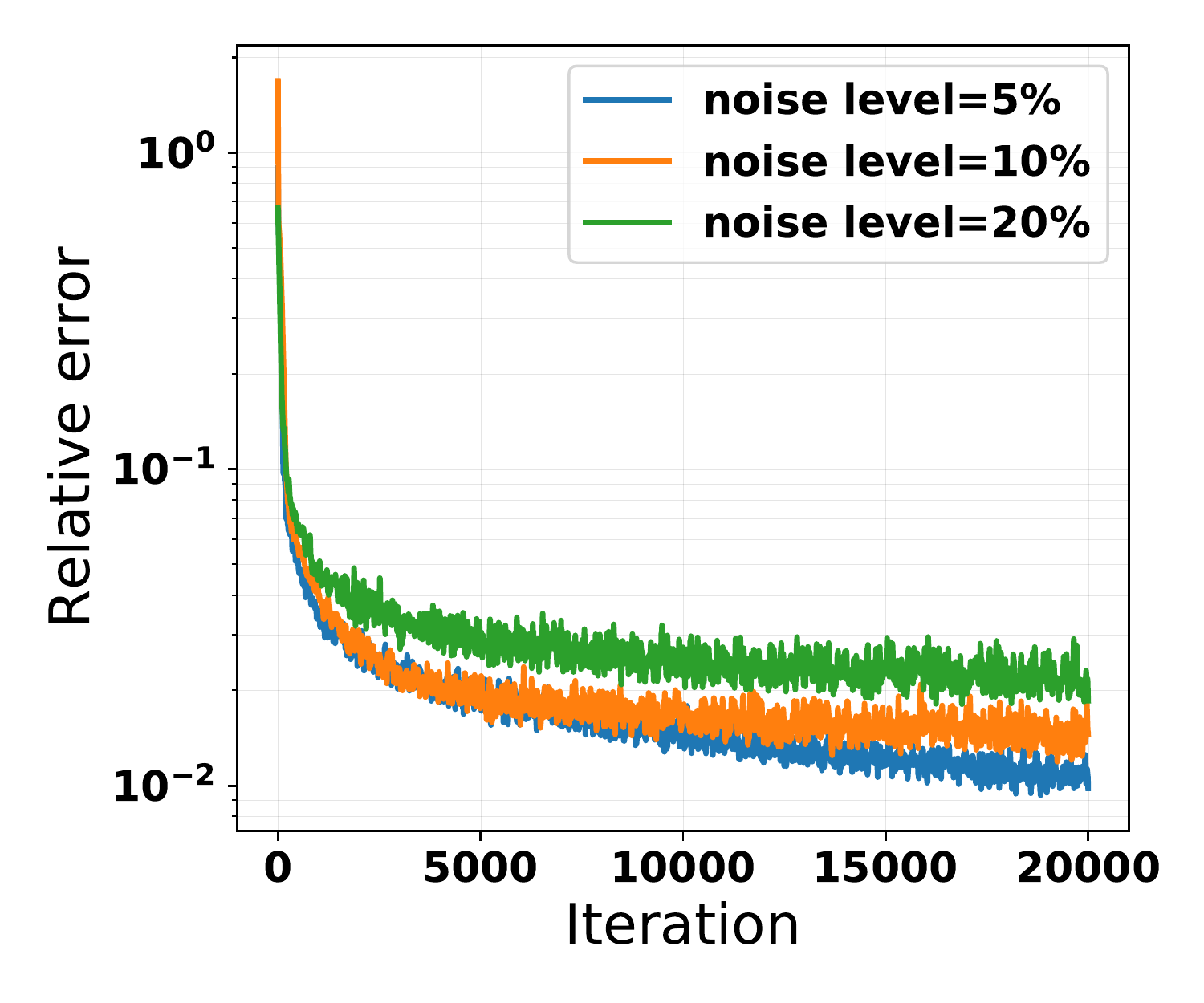}
\caption{Ground truth conductivity $\gamma^*$ (left) and relative error of $\gamma_\theta$ versus iteration (right).}
\label{subfig:pc_noisy_error}
\end{subfigure}
%
%
\begin{subfigure}[b]{\textwidth}
\includegraphics[width=0.9\linewidth]{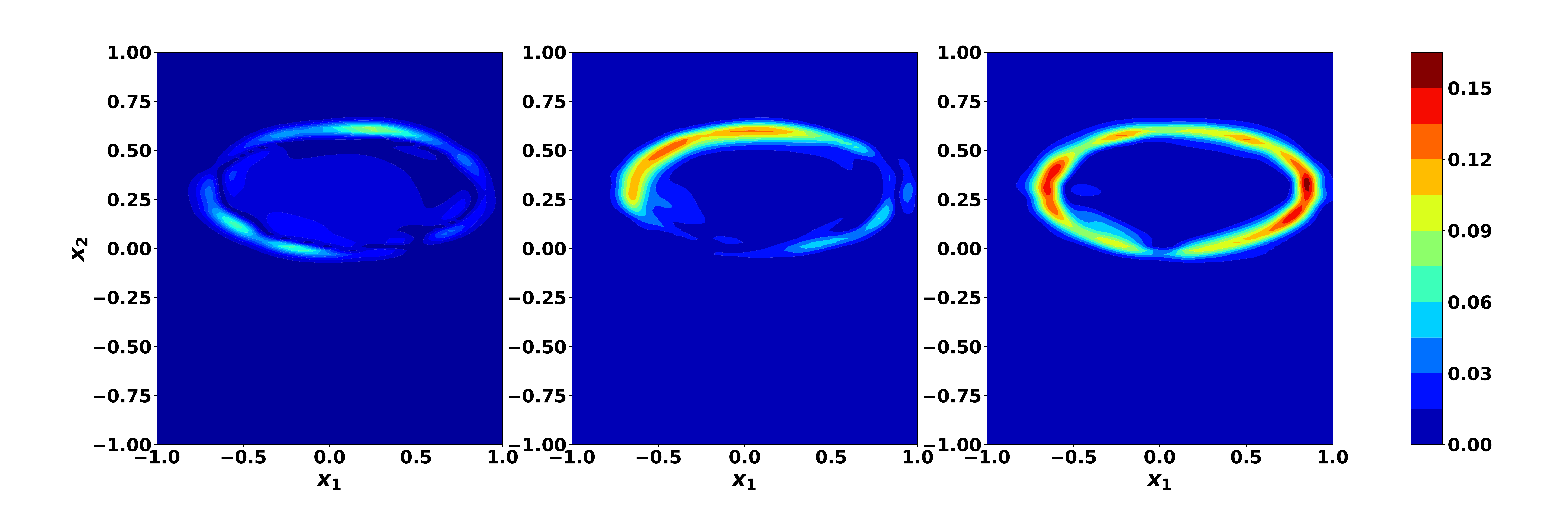}
\caption{From left to right: absolute error $|\gamma_\theta-\gamma^*|$ with measurement noise level $5\%, 10\%, 20\%$.}
\label{subfig:pc_noisy_abs_error}
\end{subfigure}
\caption{Test 3 result on \eqref{eq:eit} with nearly piecewise constant $\gamma^*$ and noisy data.}
\label{fig:discrete_noise_d5}
\end{figure}

\medskip
\noindent
\textbf{Test 3: Inverse conductivity problem with noisy measurements.}
Under the same experiment setting as Test 2, we solve the inverse problem \eqref{eq:eit} where the measurement data are perturbed by random noise for the $d=5$ case.
Specifically, we scale every measurement data value
by $1+5\% e$, $1+10\% e$, $1+20\% e$ where $e$ is drawn independently from the standard normal distribution every time, followed by a truncation into interval $[-100,100]$. We do not perturb $f$.
The results are given in Figure \ref{fig:discrete_noise_d5} in parallel to the noiseless case above, where Figure \ref{subfig:pc_noisy_error} shows the ground truth conductivity $\gamma^*$ (left) and the progress of relative error of $\gamma_{\theta}$ versus iteration number (right). The pointwise absolute error after 20,000 iterations $|\gamma_\theta-\gamma^*|$ with noise levels $5\%,10\%,20\%$ are shown in Figure \ref{subfig:pc_noisy_abs_error}.
We observe that the progress becomes more oscillatory due to the random measurement noise in Figure \ref{subfig:pc_noisy_error}, and the final reconstruction error is larger for higher noise level in Figure \ref{subfig:pc_noisy_abs_error} as expected.

\begin{figure}[t]
\centering
\begin{subfigure}[b]{.3\textwidth}
\includegraphics[width=\textwidth]{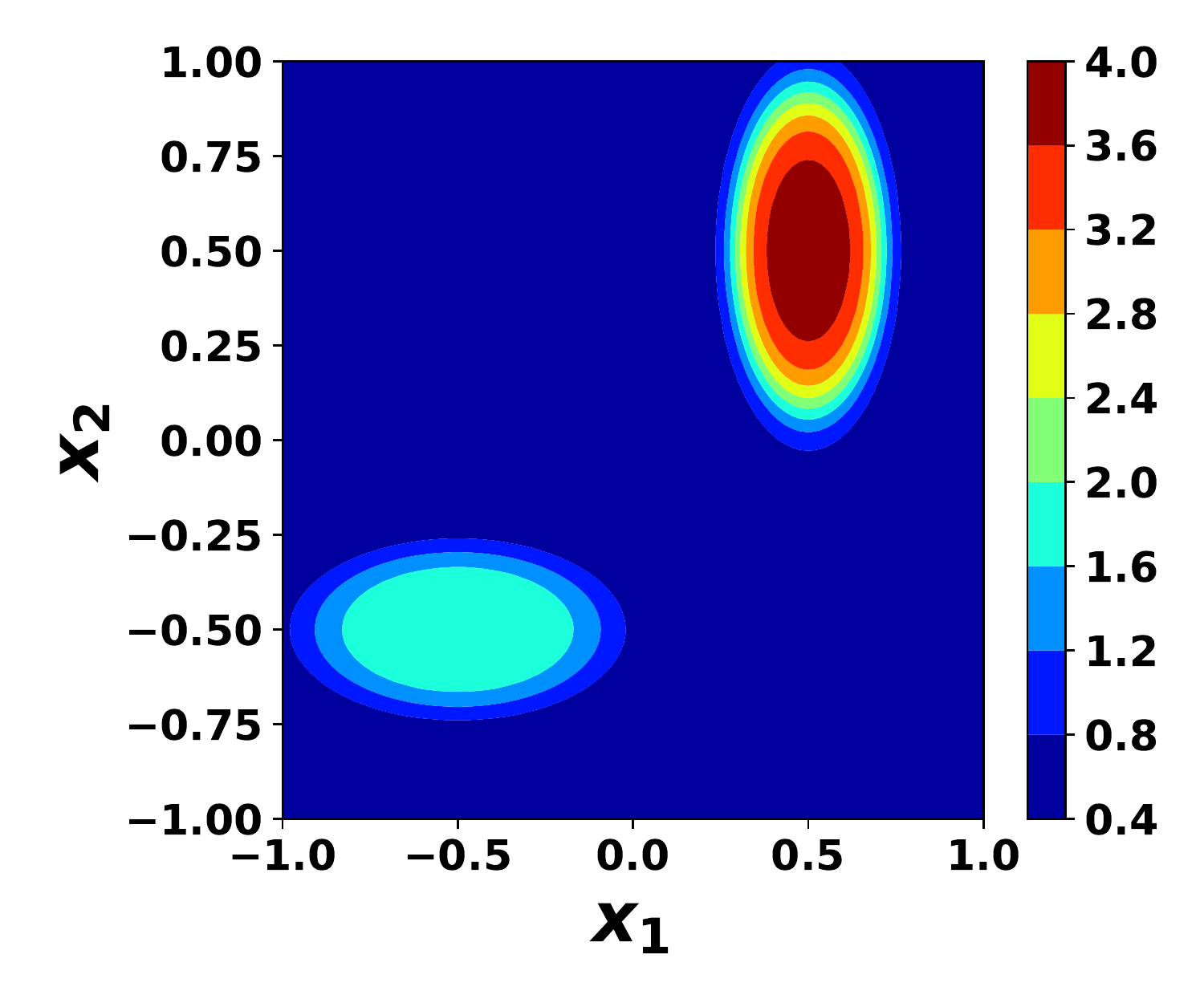}
\caption{Ground truth $\gamma^*$.}
\label{subfig:2e_gamma}
\end{subfigure}
\begin{subfigure}[b]{.3\textwidth}
\includegraphics[width=\textwidth]{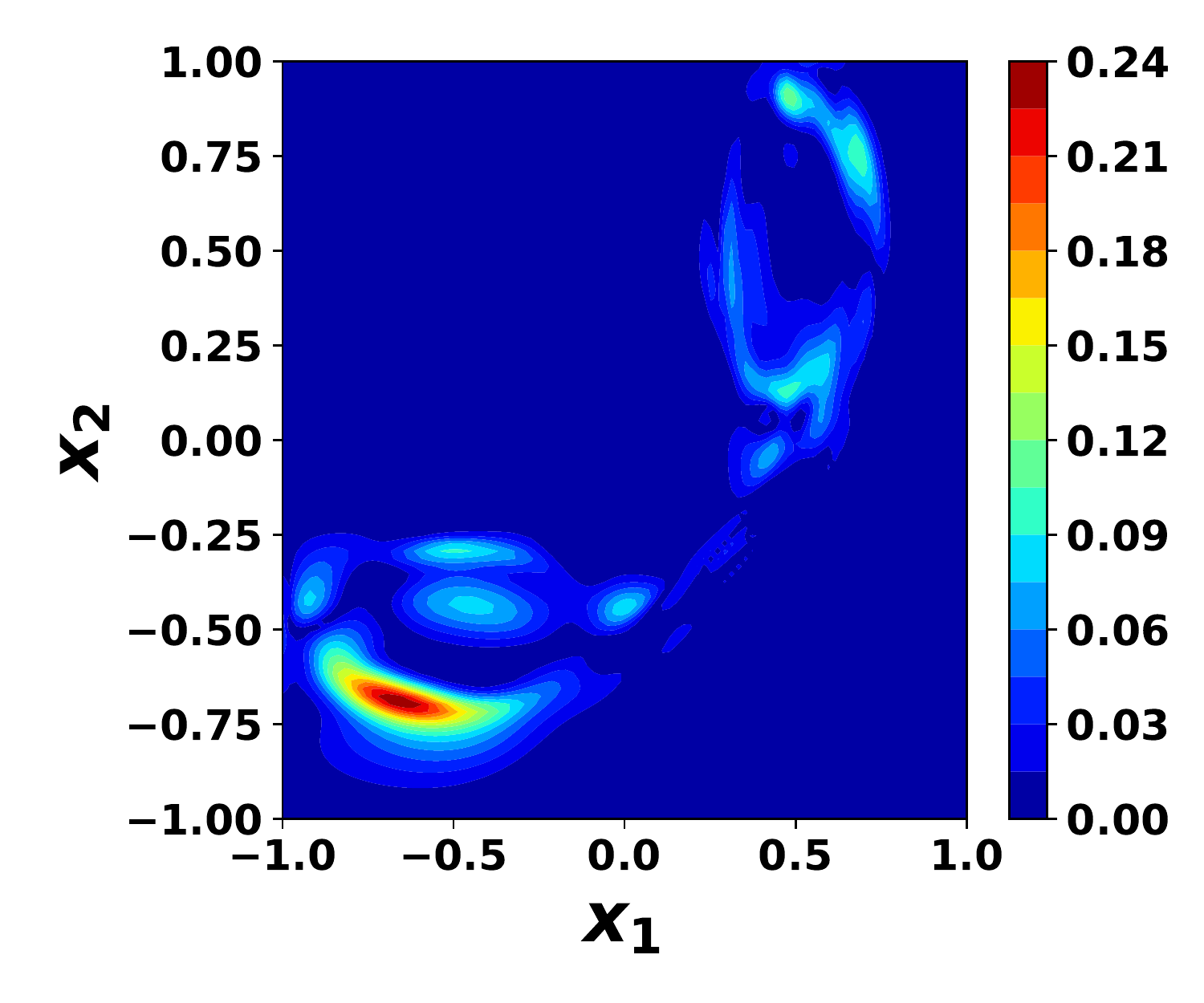}
\caption{Pointwise error $|\gamma^*-\gamma_\theta|$.}
\label{subfig:2e_abs_error}
\end{subfigure}
\begin{subfigure}[b]{.3\textwidth}
\includegraphics[width=\textwidth]{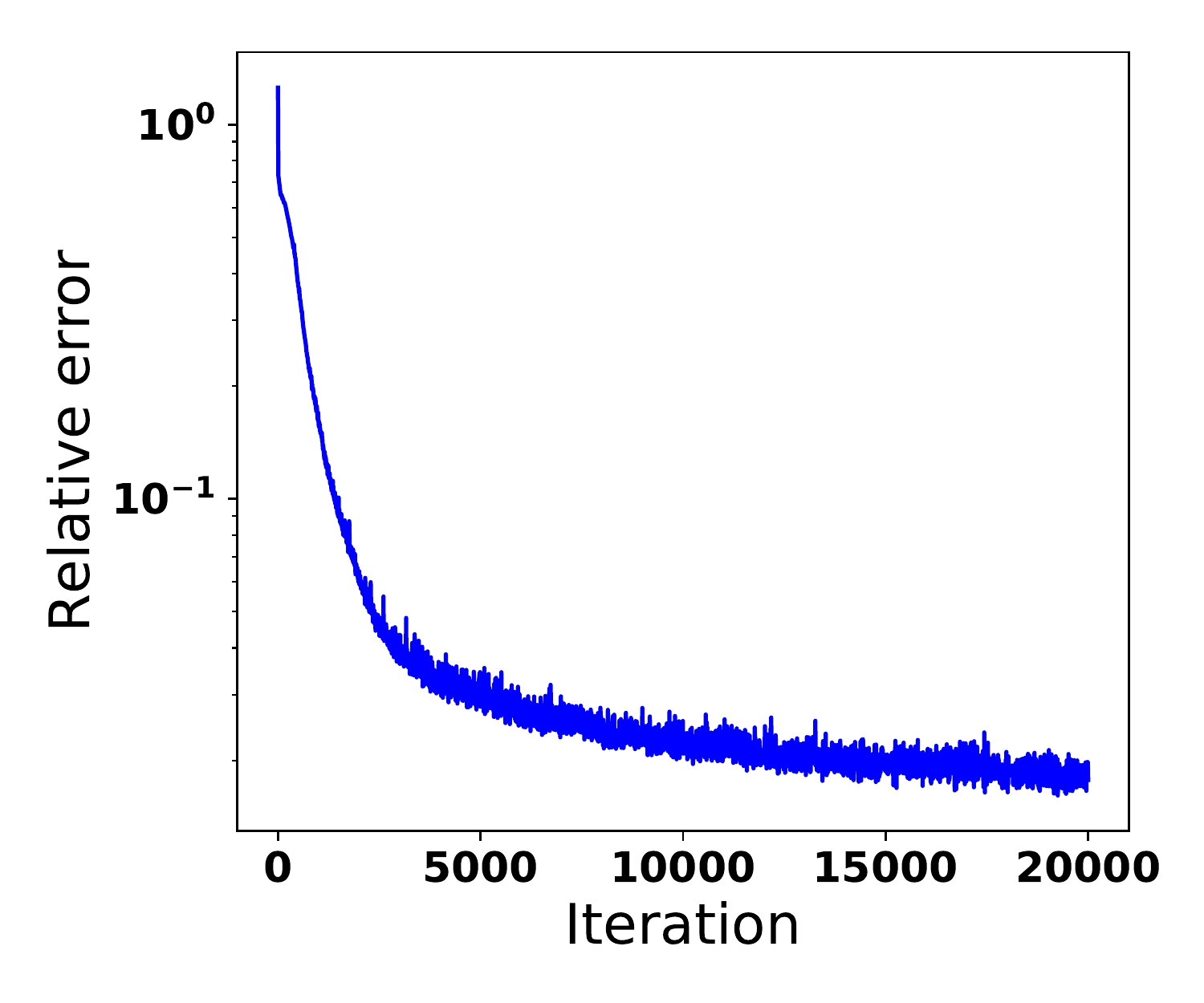}
\caption{$L_2$ relative error vs iteration.}
\label{subfig:2e_error}
\end{subfigure}
\begin{subfigure}[b]{.3\textwidth}
\includegraphics[width=\textwidth]{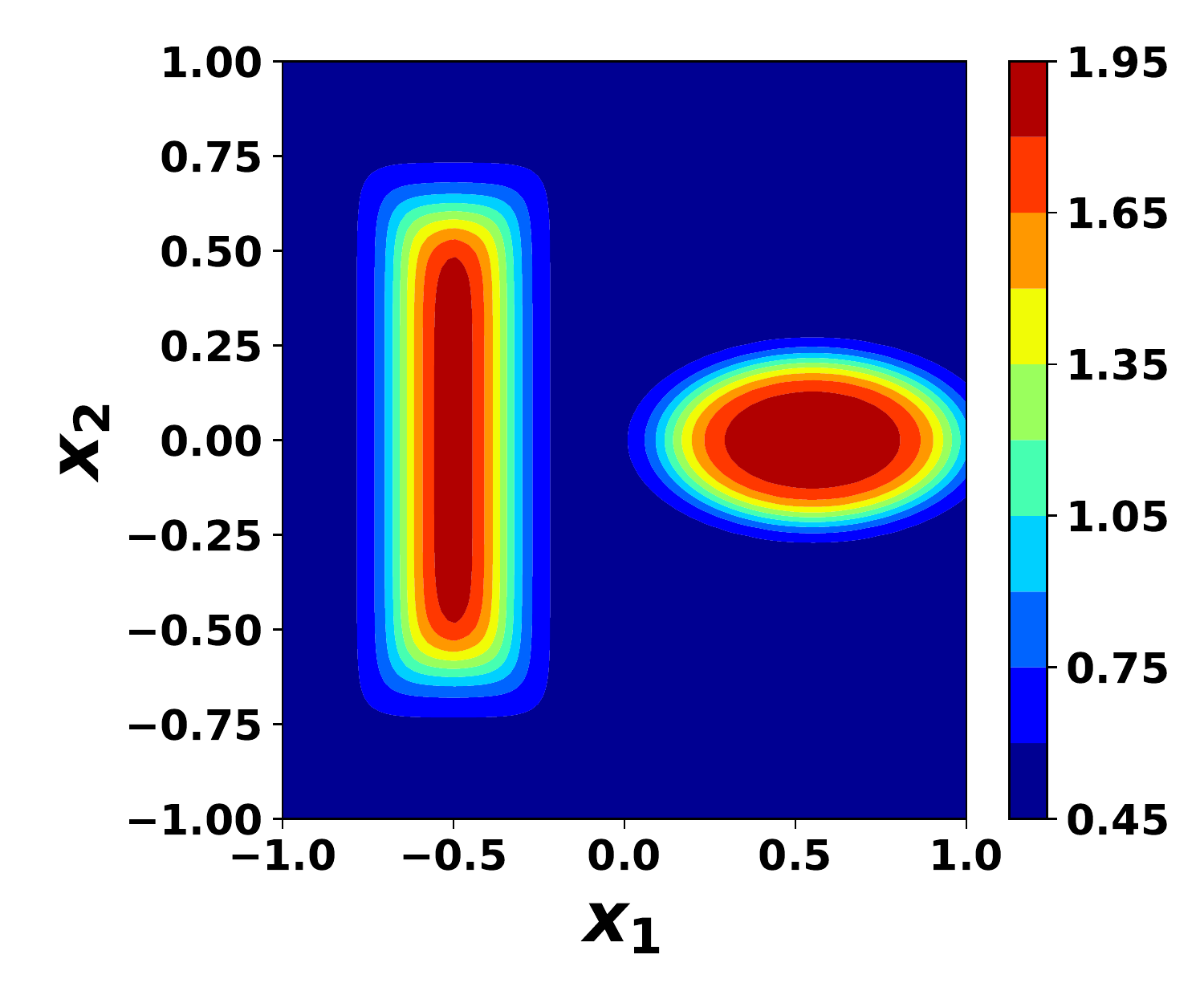}
\caption{Ground truth $\gamma^*$.}
\label{subfig:1e1s_gamma}
\end{subfigure}
\begin{subfigure}[b]{.3\textwidth}
\includegraphics[width=\textwidth]{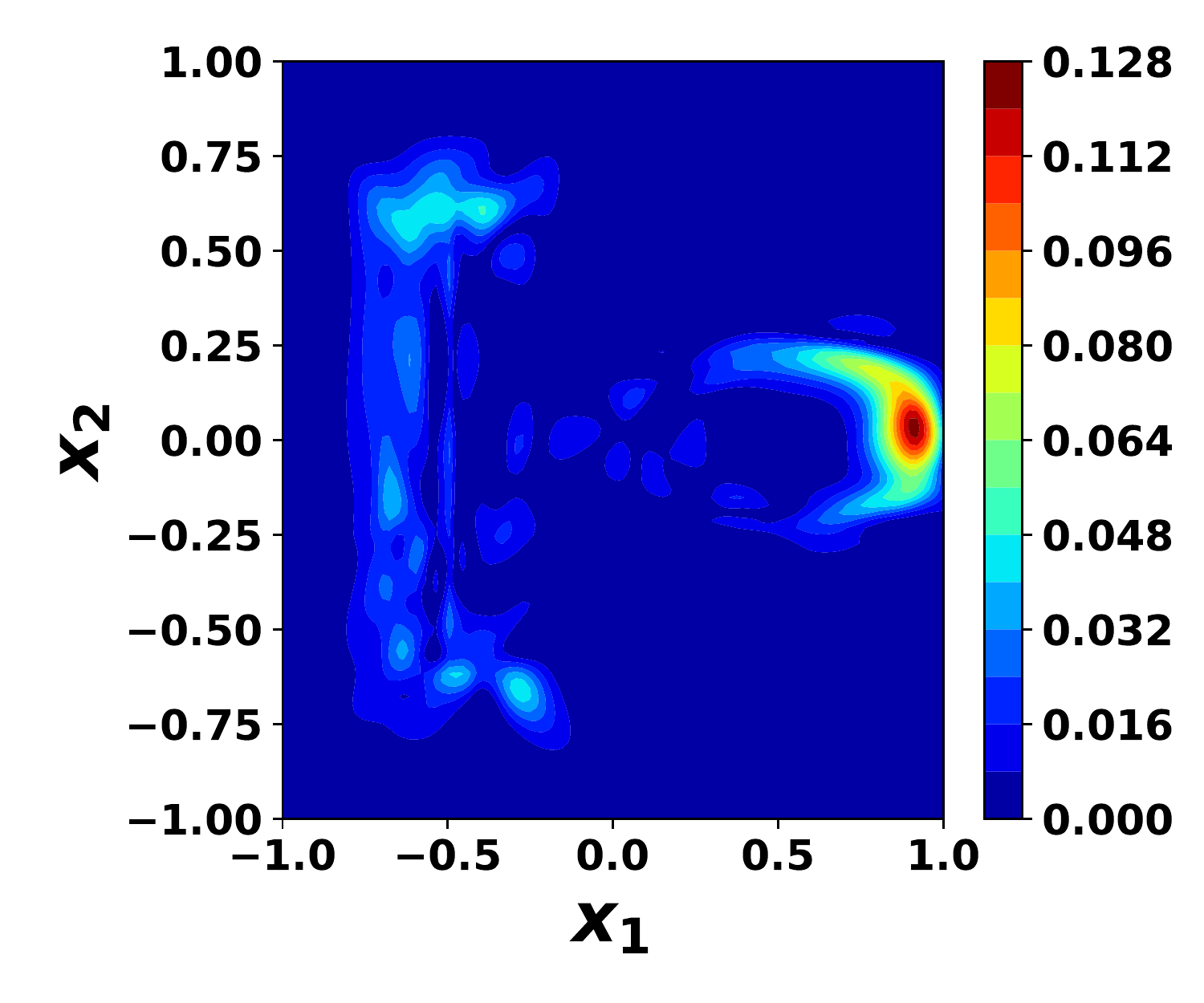}
\caption{Pointwise error $|\gamma^*-\gamma_\theta|$.}
\label{subfig:1e1s_abs_error}
\end{subfigure}
\begin{subfigure}[b]{.3\textwidth}
\includegraphics[width=\textwidth]{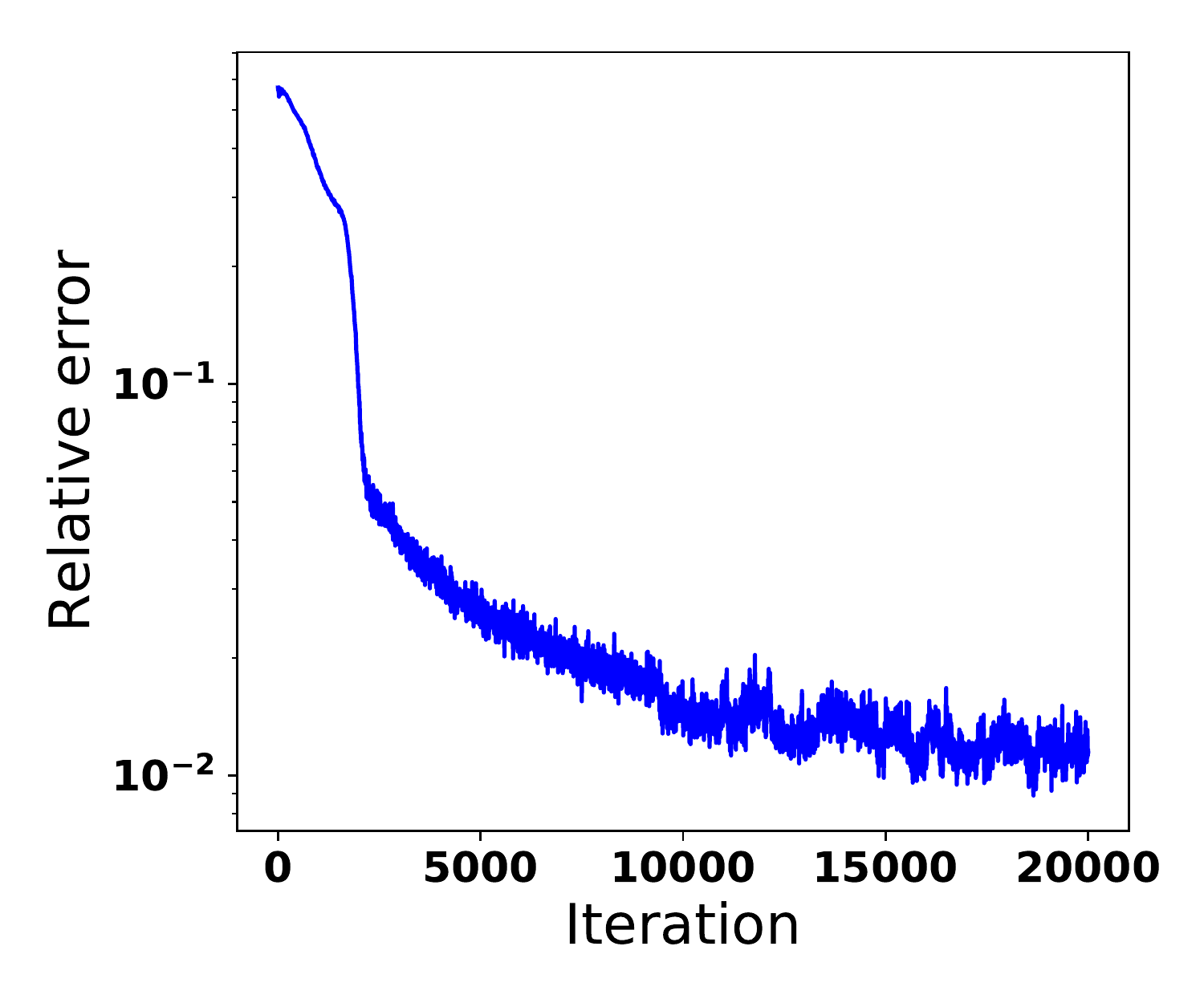}
\caption{$L_2$ relative error vs iteration.}
\label{subfig:1e1s_error}
\end{subfigure}
\begin{subfigure}[b]{.3\textwidth}
\includegraphics[width=\textwidth]{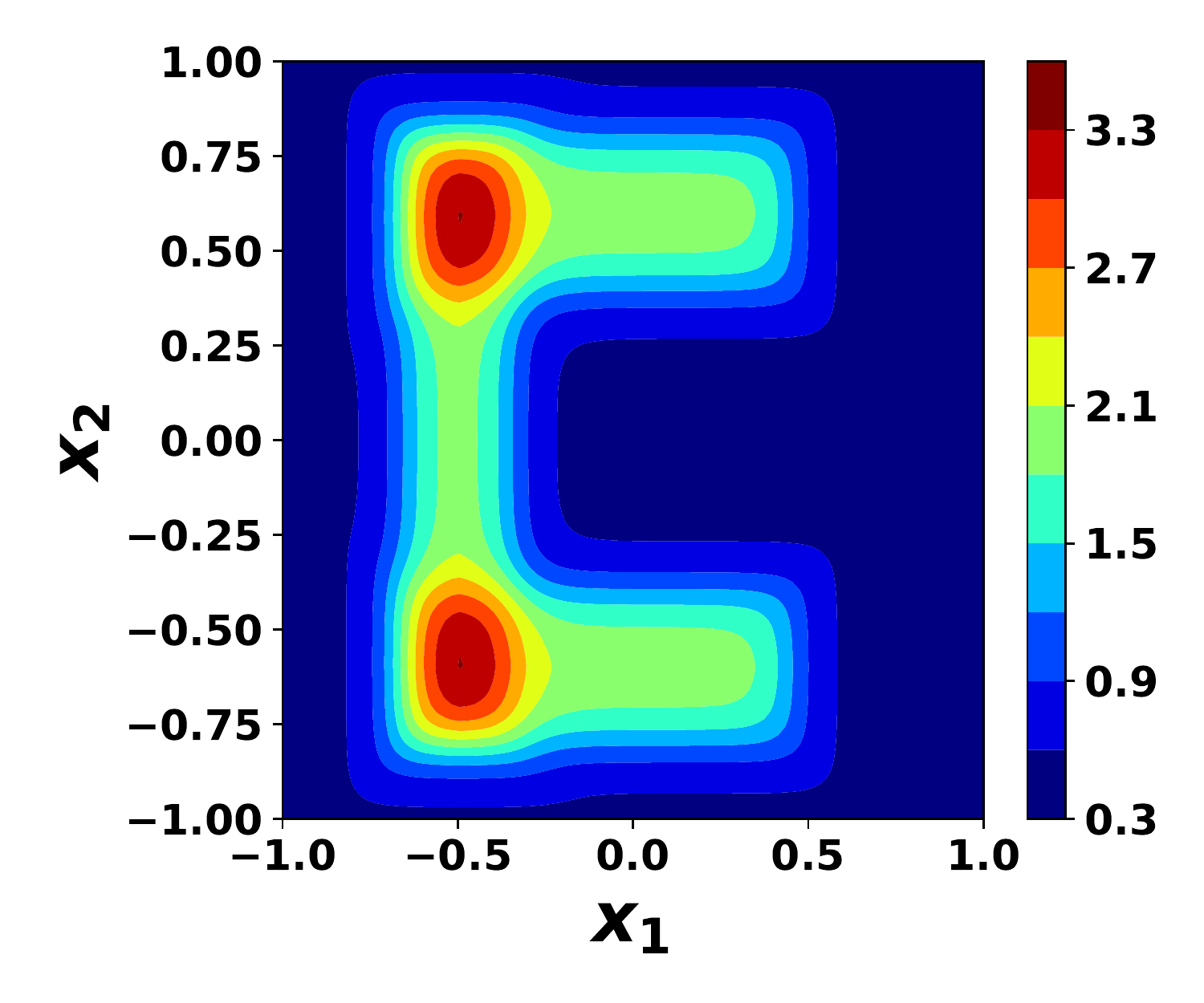}
\caption{Ground truth $\gamma^*$.}
\label{subfig:1c_gamma}
\end{subfigure}
\begin{subfigure}[b]{.3\textwidth}
\includegraphics[width=\textwidth]{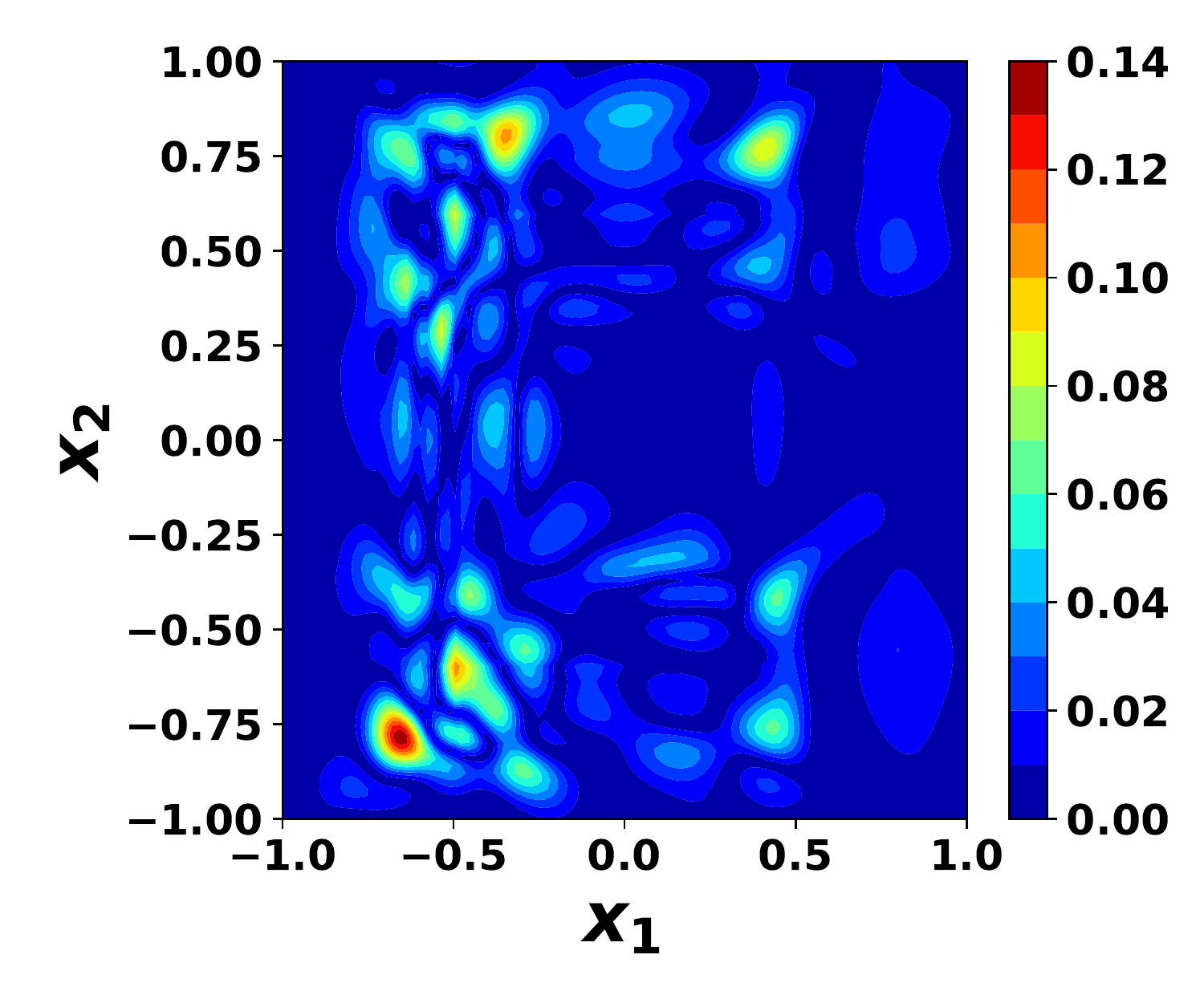}
\caption{Pointwise error $|\gamma^*-\gamma_\theta|$.}
\label{subfig:1c_abs_error}
\end{subfigure}
\begin{subfigure}[b]{.3\textwidth}
\includegraphics[width=\textwidth]{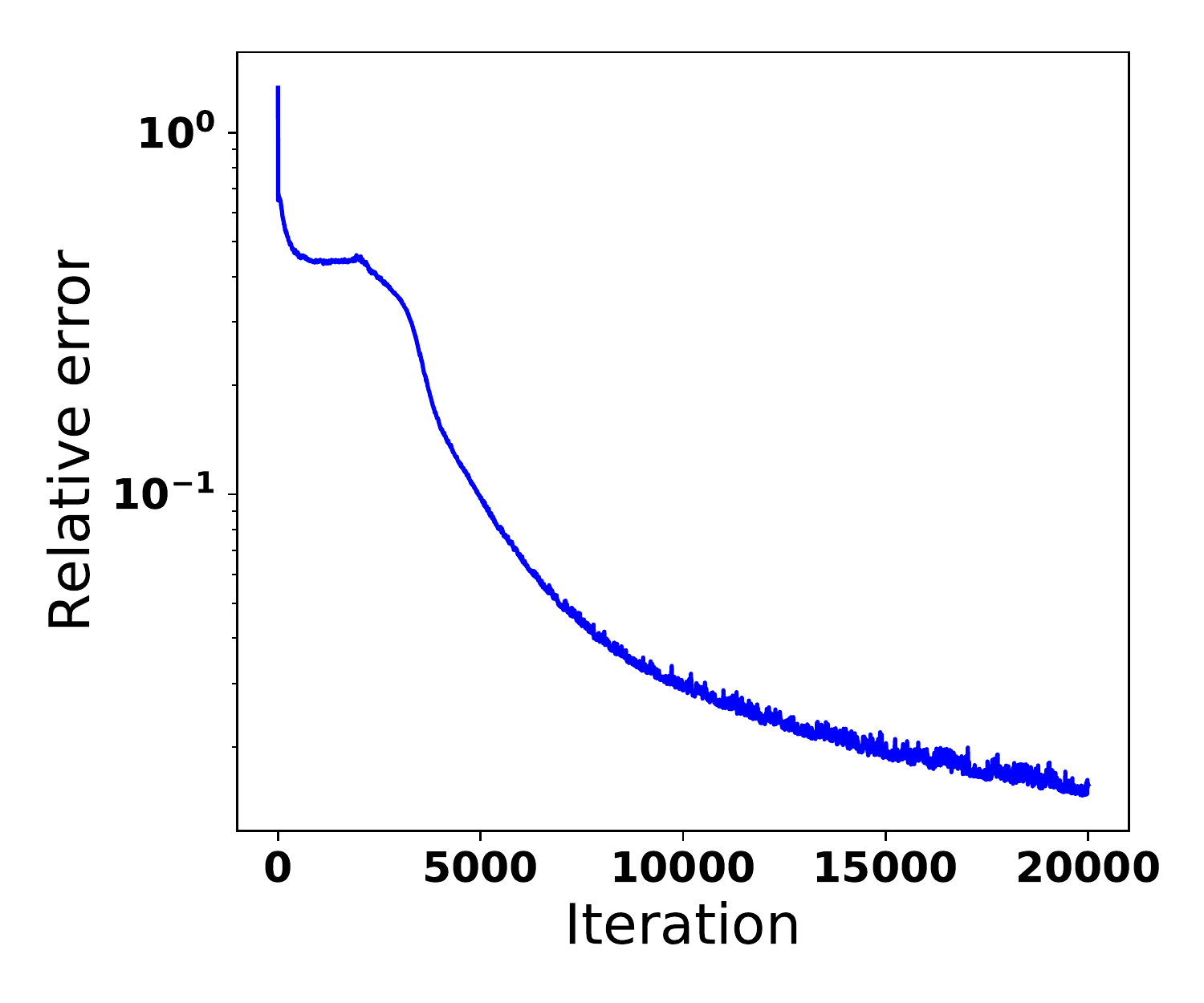}
\caption{$L_2$ relative error vs iteration.}
\label{subfig:1c_error}
\end{subfigure}
\caption{Test 4 results on \eqref{eq:eit} with problem dimension $d=5$. (a)(b)(c) two separate modes in $\gamma^*$; (d)(e)(f) two separate modes (one has sharp corner) in $\gamma^*$; (g)(h)(i) nonconvex shaped $\gamma^*$. }
\end{figure}

\medskip
\noindent
\textbf{Test 4: Inverse conductivity problem with different features in $\gamma$.}
We consider several cases with more challenging ground truth conductivity $\gamma^*$.
The first case has $\gamma^*$ with two disjoint modes.
We define $\Omega=(-1,1)^5$, and set $\Omega_1 = \{x:|x-c_1|_{\Sigma_1}^2 \le 0.4^2\}$ and $\Omega_2 = \{x:|x-c_2|_{\Sigma_2}^2 \le 0.4^2\}$, where $c_1=(-0.5,-0.5,0,0,0)$, $c_2=(0.5,0.5,0,0,0)$, $\Sigma_1=\mathrm{diag}(0.81,2,009,0.09,0.09)$, and $\Sigma_2=\mathrm{diag}(2,0.81,009,0.09,0.09)$.
We set the conductivity $\gamma^*$ to $4$ in $\Omega_1$, $2$ in $\Omega_2$, and $0.5$ in $(\Omega_1\cup \Omega_2)^c$.
We also smooth $\gamma^*$ using a Gaussian kernel.
{The setup of the ground truth conductivity $\gamma^*$, the ground truth potential $u^*$, and the source function $f(x)$ are provided in the table \ref{tab:prob_set} (see Test 4(1)).}
%
%
We set $N_r= 200,000$, $N_b= 100d$, $\beta'=10$, and run Algorithm \ref{alg:iwan} for $20,000$ iterations.
Figure \ref{subfig:2e_gamma} shows the ground truth $\gamma^*$ (left) and the recovered conductivity $\gamma_\theta$ with relative error $1.77\%$ (right).
Figure \ref{subfig:2e_error} plots the progress of relative error of $\gamma_{\theta}$ versus iteration number.
In the second case, we follow the same setting but define $\Omega_1=\{x:|x_1+0.5|\le 0.15,|x_2|\le0.6\}$ (which has sharp corner) and $\Omega_2=\{x:|x-c|_{\Sigma}^2 \le 0.4^2\}$ where $c=(0.55, 0, 0, 0, 0)$, $\Sigma=\mathrm{diag}(1, 4, 0, 0, 0)$, and set $\gamma^*=2$ in $\Omega_1\cup\Omega_2$ and $0.5$ in $(\Omega_1\cup\Omega_2)^c$.
{We again smooth $\gamma^*$ and provide the ground truth conductivity $\gamma^*$, the ground truth potential $u^*$, and the source function $f$ in the table \ref{tab:prob_set} (see Test 4(2)).}
%
%
We set $\beta'=1$ and again run Algorithm \ref{alg:iwan} for $20,000$ iterations.
The recovered $\gamma_\theta$ (with relative error $1.15\%$) and the progress of relative error are shown in Figure \ref{subfig:1e1s_gamma} and \ref{subfig:1e1s_error} respectively.
Lastly, we consider a non-convex shaped $\gamma$, and show the recovered $\gamma_\theta$ (with relative error $1.57\%$) and the progress of relative error in Figure \ref{subfig:1c_gamma} and \ref{subfig:1c_error} respectively.
We set the domain $\Omega=(-1,1)^5$ and define $\Omega_j=\{x\in \Omega:\sum^{2}_{i=1}|x_i-c_j(i)|\le r_j(i)\}, j= 1,2,3$, where $c_1=(-0.5, 0), c_2=(-0.1, 0.6), c_3= (-0.1, -0.6)$ and $r_1=(0.15, 0.8), r_2=r_3=(0.55, 0.2)$. We set $\gamma^*=4$ in $\Omega_0=(\Omega_1\cap\Omega_2)\cup(\Omega_1\cap\Omega_3)$ and $2$ in $\Omega_1\cup\Omega_2\cup\Omega_3/\Omega_0$ and $0.5$ in $(\Omega_1\cup\Omega_2\cup\Omega_3)^c$.
{Once again, we sightly smooth the conductivity $\gamma^*$ and provide the setup of the ground truth conductivity $\gamma^*$, the ground truth potential $u^*$ and the source function $f$ in the table \ref{tab:prob_set} (see Test 4(3)).}
%
%

\begin{figure}[t]
\centering
\begin{subfigure}[b]{.3\textwidth}
\includegraphics[width=\textwidth]{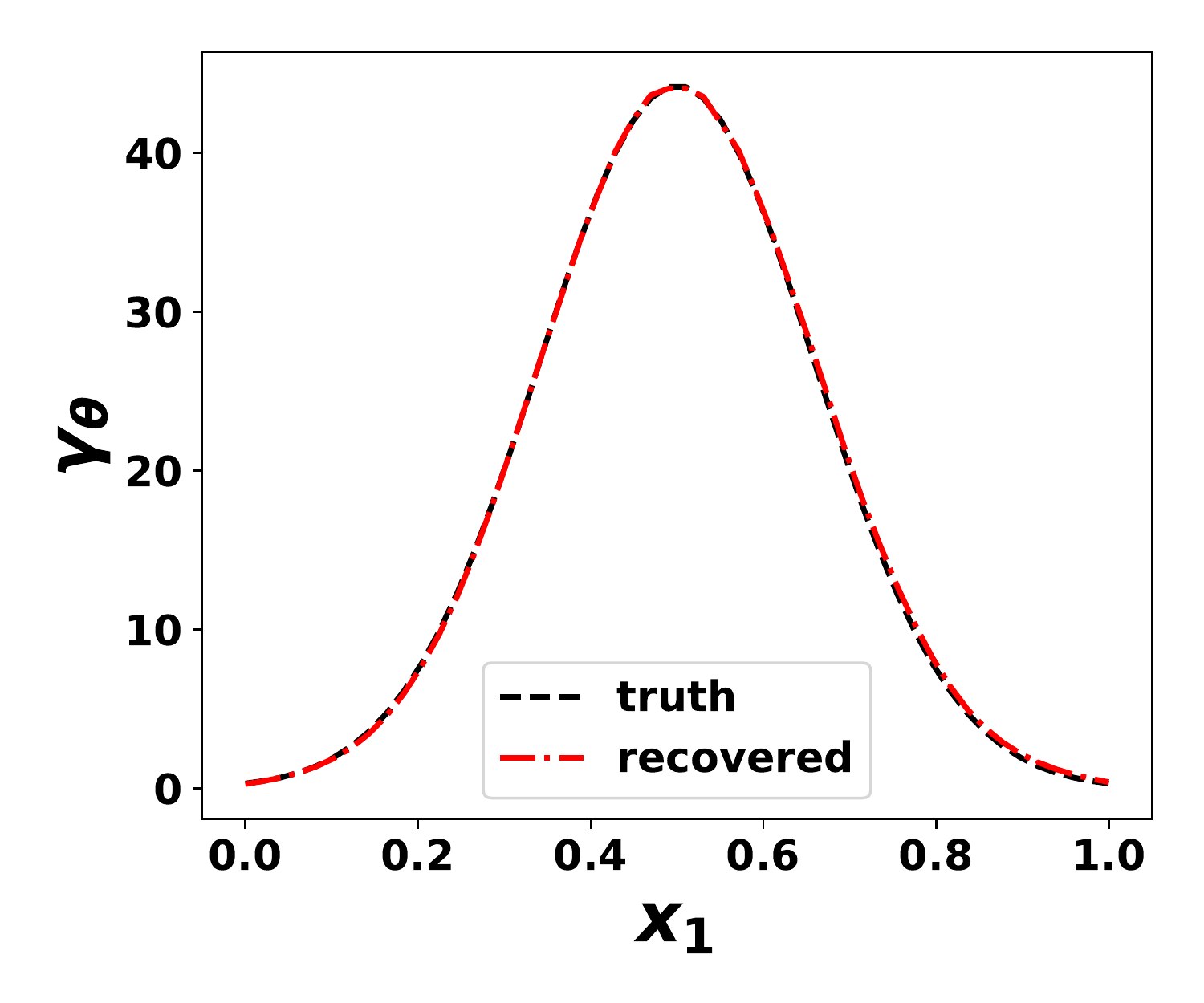}
\caption{True $\gamma^*$ vs recovered $\gamma_\theta$.}
\label{subfig:eit_gamma}
\end{subfigure}
\begin{subfigure}[b]{.3\textwidth}
\includegraphics[width=\textwidth]{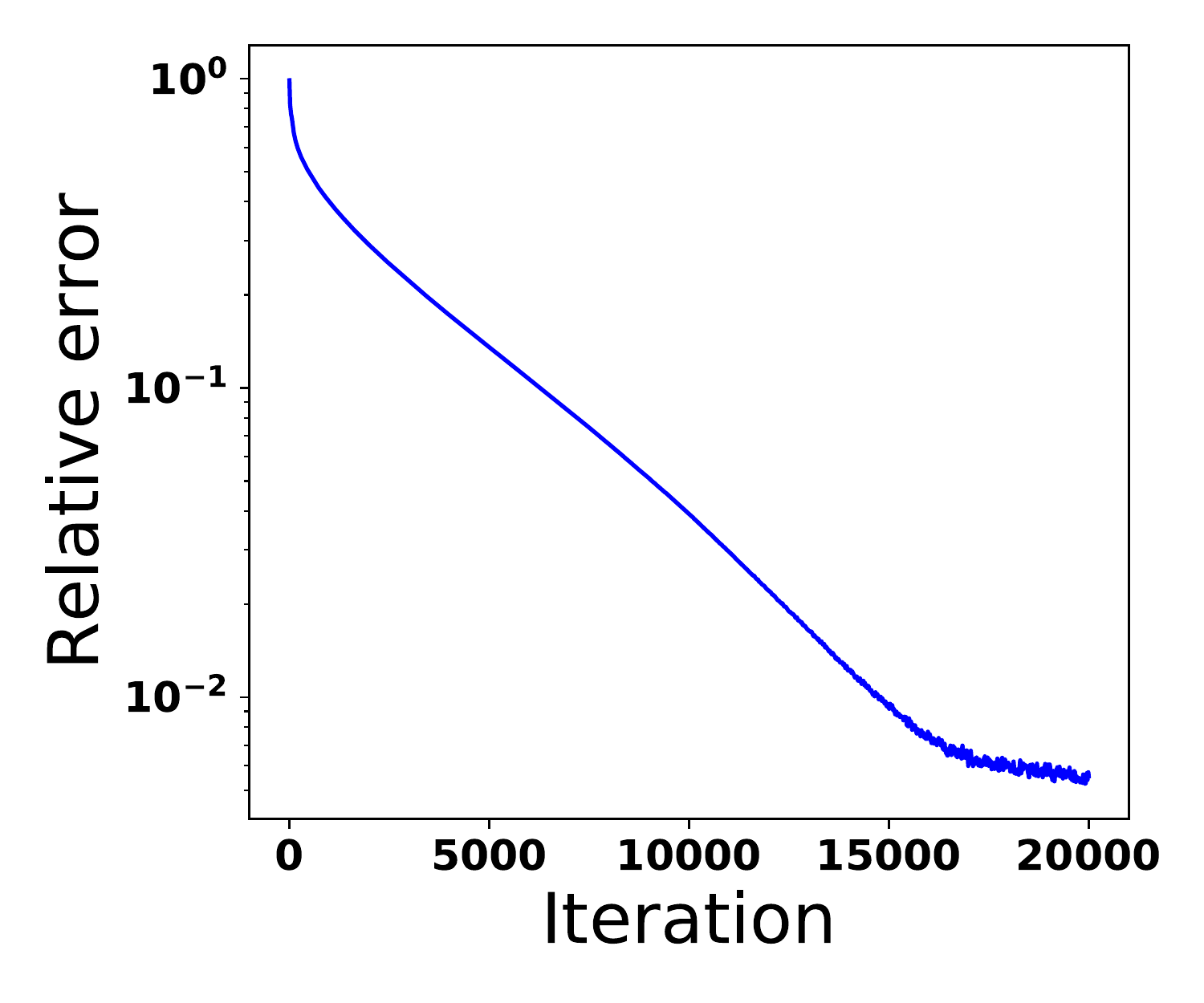}
\caption{$L_2$ relative error vs iteration.}
\label{subfig:eit_error}
\end{subfigure}
\caption{Test 5 result on artificial 5D EIT problem \eqref{eq:eit} with boundary condition on $\gamma \nabla u\cdot \vec{n}$.}
\label{fig:eit}
\end{figure}

\medskip
\noindent
\textbf{Test 5: EIT problem.}
We consider an artificial 5D EIT problem \eqref{eq:eit} on $\Omega=(0,1)^5$ but replace \eqref{eq:ip2b} with a different boundary condition given by $\gamma \nabla u \cdot \vec{n}$ on $\partial \Omega$, where $\vec{n}$ is the outer normal vector at the boundary point.
We define $\Gamma_1=\{x\in\partial\Omega:x_1=0,1\}$ and $\Gamma_2=\partial\Omega \setminus \Gamma_1$.
We set ground truth conductivity $\gamma^*(x)=\pi^{-1}\exp\{(d-1)\pi^2(x_1-x_1^2)/2\}$ and potential $u^*(x)=\exp\{(d-1)\pi^2(x_1^2-x_1)/2\}\cdot \Pi^{d}_{i=2}\sin(x_i)$, and compute the corresponding boundary value as our input data.
We set the source term $f=0$ in \eqref{eq:ip2a}, $N_r= 100,000$, $N_b= 100d$, $\beta'=10$, and run Algorithm \ref{alg:iwan} for 20,000 iterations. The $x_1$ cross section of the recovered $\gamma_\theta$ (with relative error $0.56\%$) and the  progress of relative error versus iteration number are shown in Figure \ref{subfig:eit_gamma} and \ref{subfig:eit_error}, respectively.

\begin{figure}[t]
\centering
\begin{subfigure}[b]{.3\textwidth}
\includegraphics[width=\textwidth]{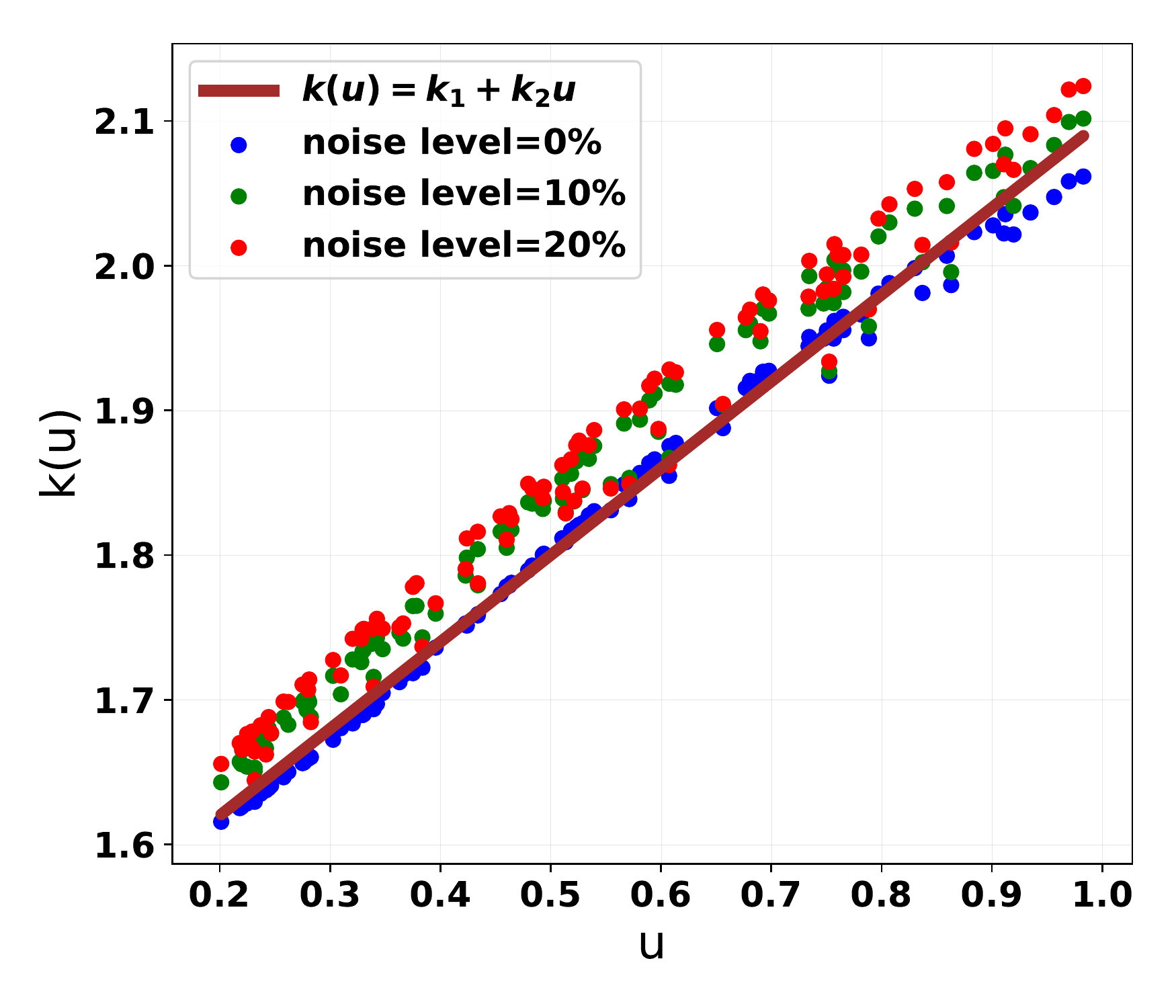}
\caption{Recovered $(u_\theta,\gamma_\theta)$.}
\label{subfig:itp_gamma}
\end{subfigure}
\begin{subfigure}[b]{.3\textwidth}
\includegraphics[width=\textwidth]{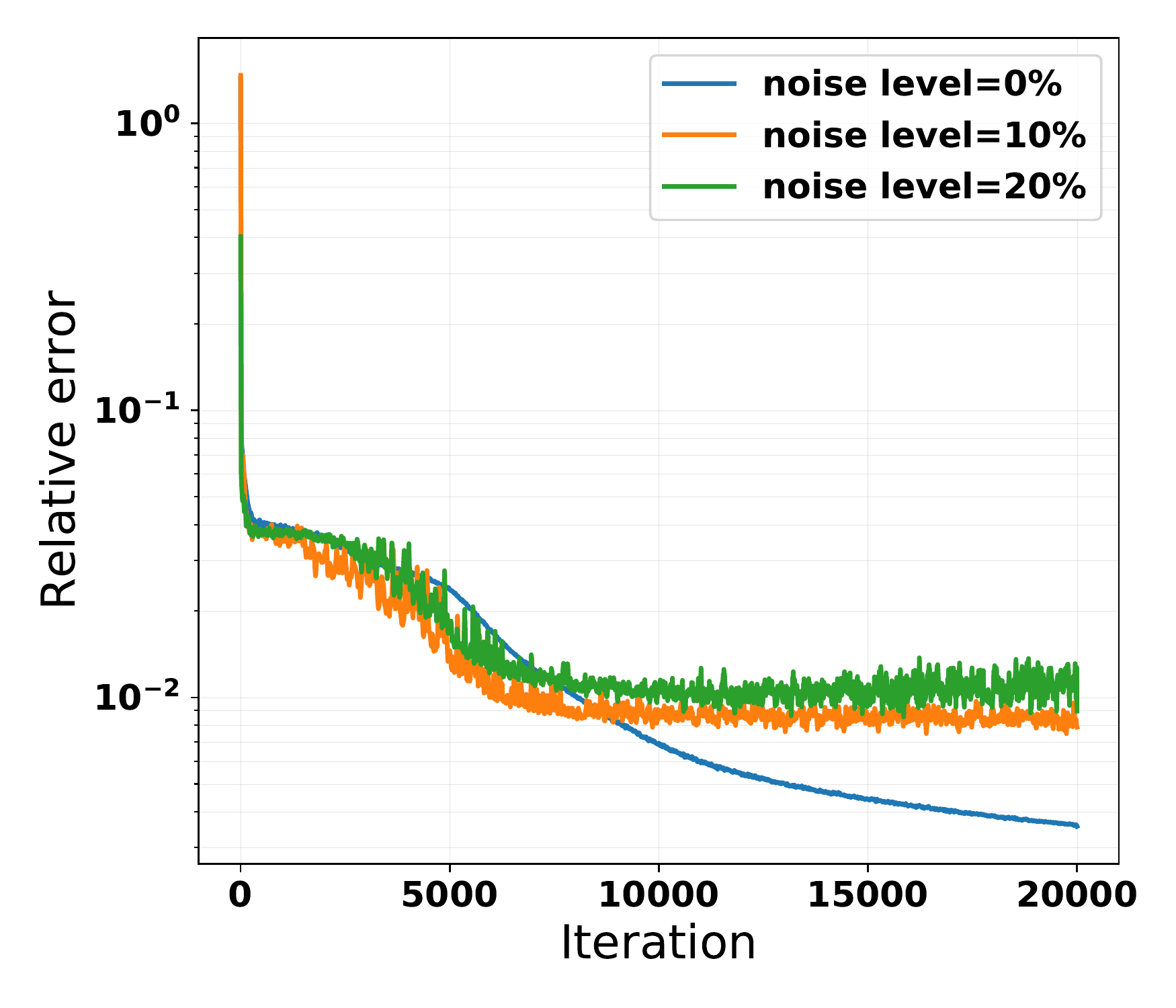}
\caption{relative error of $\gamma_\theta$.}
\label{subfig:itp_error}
\end{subfigure}
\caption{Test 6 result on inverse thermal conductivity problem with dimension $d=5$ and noise levels $0\%, 10\%$ and $20\%$. Ground truth relation between $u^*$ and $\gamma^*$ is $\gamma^*=k_1+k_2u^*$ where $k_1=1.5$ and $k_2=0.6$.}
\label{fig:itp}
\end{figure}

\medskip
\noindent
\textbf{Test 6: Inverse thermal conductivity problem involving time.}
We consider an inverse thermal conductivity problem of \eqref{eq:ip} with temporally varying $u(x,t)$ as follows,
\begin{subequations}
\label{eq:problem_heat}
\begin{align}
\partial_t u - \nabla \cdot (\gamma \nabla u) - f = 0, & \quad \text{in}\ \Omega_{T} = \Omega\times[0,T]\\
u-u_i = 0, & \quad \text{in}\ \Omega \times \{0\} \\
\nabla u \cdot \vec{n} - u_n = 0,\ u-u_b = 0, \ \gamma - \gamma_b = 0, & \quad \text{on}\ \partial\Omega\times[0,T]
\end{align}
\end{subequations}
where $\Omega = (0,1)^5\subset \mathbb{R}^5$ and {the final time $T=1$, $\gamma$ is the thermal conductivity, $u$ indicates the temperature and $f(x,t)$ is the source function which indicates the rate of heat generation per unit volume,} where $\vec{n}$ is the unit outer normal vector of $\partial\Omega$, $u_i(x)$ for $x\in\Omega$ is the given initial value, and $u_n(x,t),u_b(x,t),\gamma_b(x,t)$ for $(x,t)\in \partial\Omega \times[0,T]$ are given boundary values.
In this problem, we would like to recover the thermal conductivity $\gamma(u)$ which is a function of the temperature $u$ in the standard setting, but we simply treat $\gamma(x,t)\coloneqq\gamma(u(x,t))$ as a function of $(x,t)$ in our experiment here.
%
We set the source function $f(x,t)$ as follows,
\begin{equation*}
f(x,t)= \pi^2 \del[2]{k_1+k_2 s(t)\sum^{d}_{i=1}\sin(\pi x_i)-\lambda_2} s(t)\sum^{d}_{i=1}\sin(\pi x_i)-k_2\pi^2 s^2(t)\sum^{d}_{i=1}\cos^2(\pi x_i)
\end{equation*}
where $s(t)\coloneqq \exp(-3t/2)/5$, the initial value $u_i=\lambda_1\sum^{d}_{i=1}\sin(\pi x_i)$, the Neumann boundary value of $u$ as $u_n(x,t)= \pi s(t) (\cos(\pi x_1), \cdots, \cos(\pi x_d)) \cdot \vec{n}$.
We set the ground truth $\gamma^*(x,t)=k_1 + k_2 u^*(x,t)$ where $k_1=1.5$, $k_2=0.6$, and $u^*(x,t)=s(t)\sum^{d}_{i=1}\sin(\pi x_i)$, and use noisy boundary value measurements $u_b = (1+\sigma e)u^*$ and $\gamma_b = (1+\sigma e)\gamma^*$, where $e$ is independently drawn for $u_b$ and $\gamma_b$ and all $x\in\partial\Omega$ from the standard normal distribution followed by a truncation to $[-100, 100]$, and the noise levels are set to $\sigma=0\%, 10\%, 20\%$.
We consider the case with problem dimension $d=5$, and set $N_r=100,000$, $N_b= 100d$, $\beta'=1$ and $\beta=1,000$.
The results are shown in Figure \ref{fig:itp}, where Figure \ref{subfig:itp_gamma} plots the sampled values of recovered $(u_\theta,\gamma_\theta)$ in comparison with the ground truth relation $\gamma^*=k_1+k_2u^*$, and Figure \ref{subfig:itp_error} shows the progress of relative error of $\gamma_\theta$ versus iteration number for the three different noise levels.
The reconstructions of $\gamma_\theta$ are all faithful, while higher noise levels decreases the accuracy and make convergence to true solution more challenging.

\begin{figure}
\centering
\begin{subfigure}[b]{.6\textwidth}
\includegraphics[width=.875\textwidth]{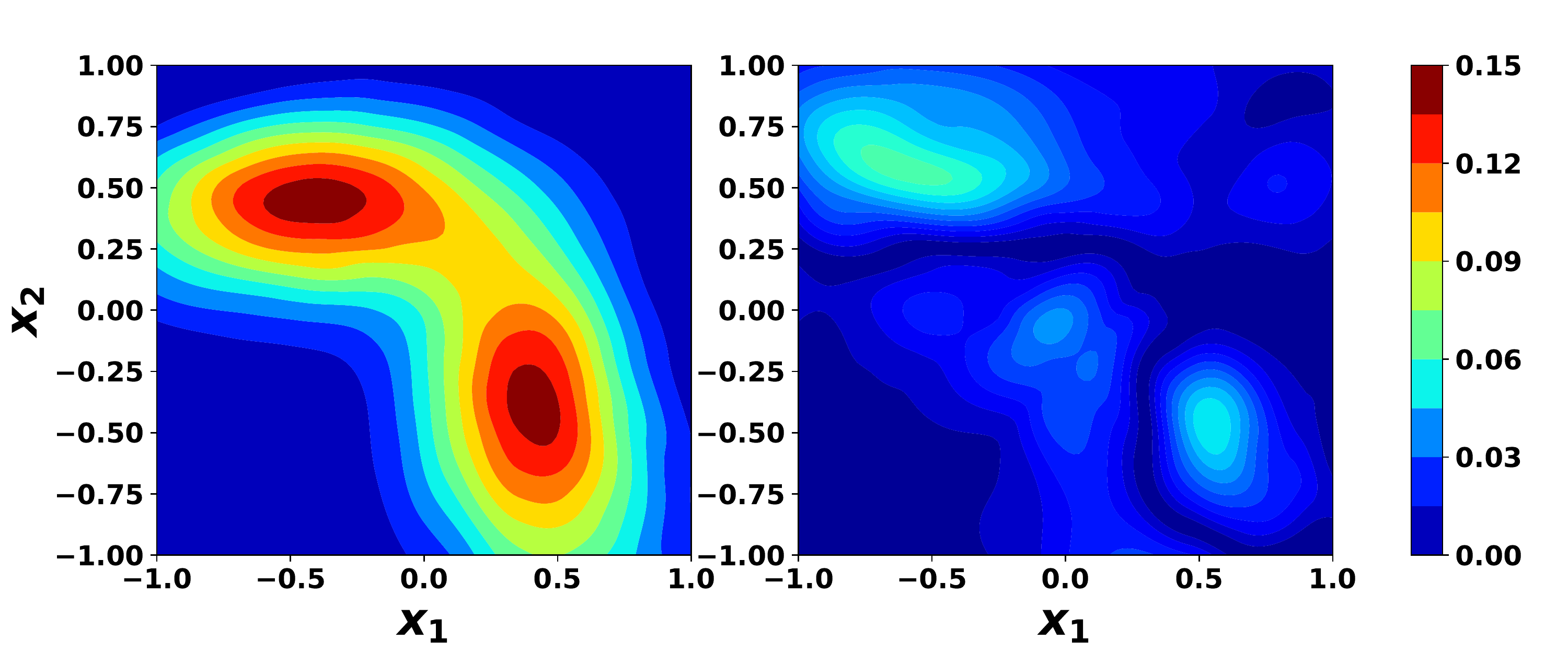}
\caption{$|\gamma-\gamma^*|$ of PINN (left) and IWAN (right) on smooth $\gamma^*$}
\label{subfig:pinn_smooth_gamma}
\end{subfigure}
\begin{subfigure}[b]{.32\textwidth}
\includegraphics[width=.875\textwidth]{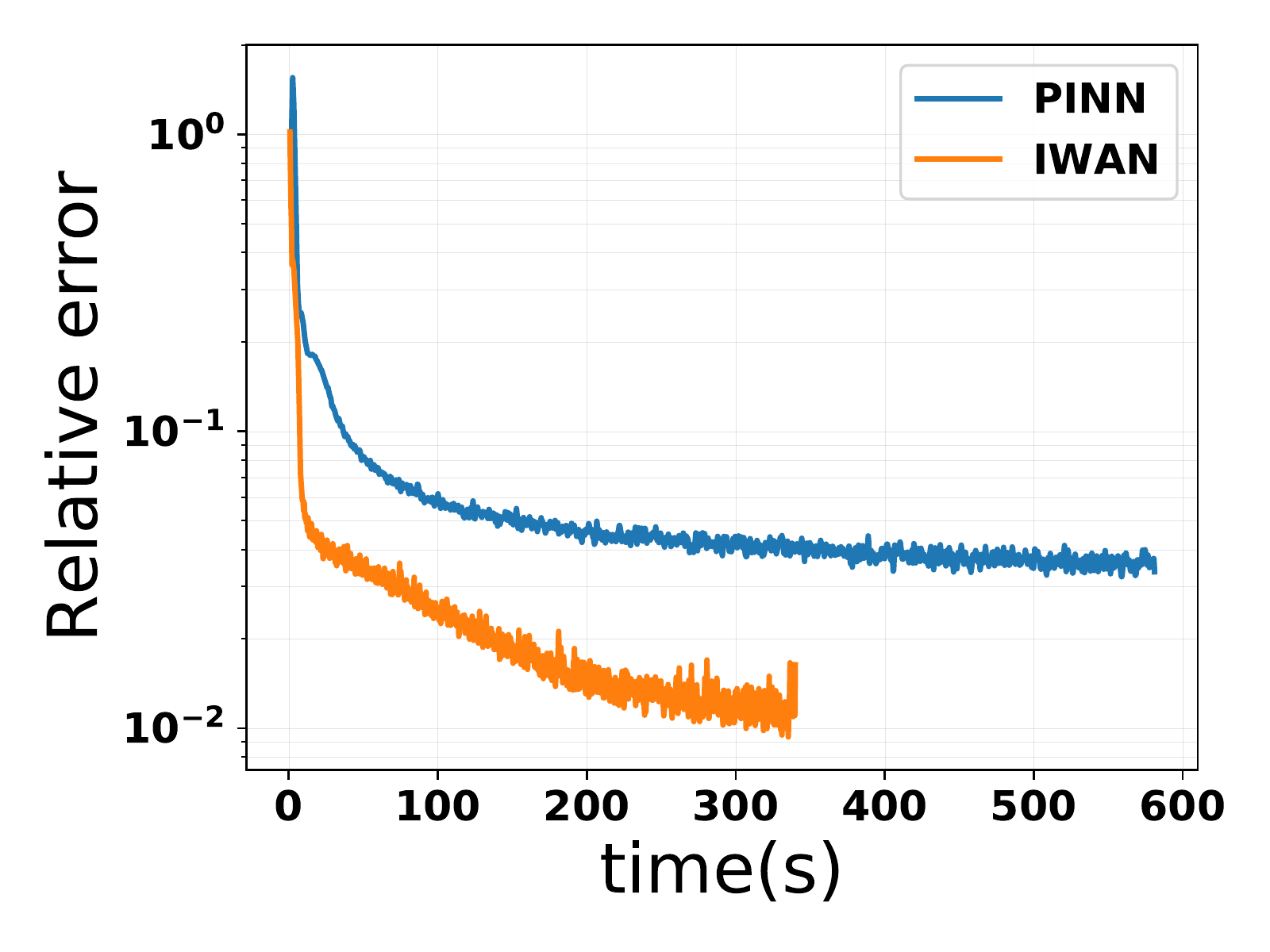}
\caption{Relative error vs time (s).}
\label{subfig:pinn_smooth_time}
\end{subfigure}
\begin{subfigure}[b]{.6\textwidth}
\includegraphics[width=.875\textwidth]{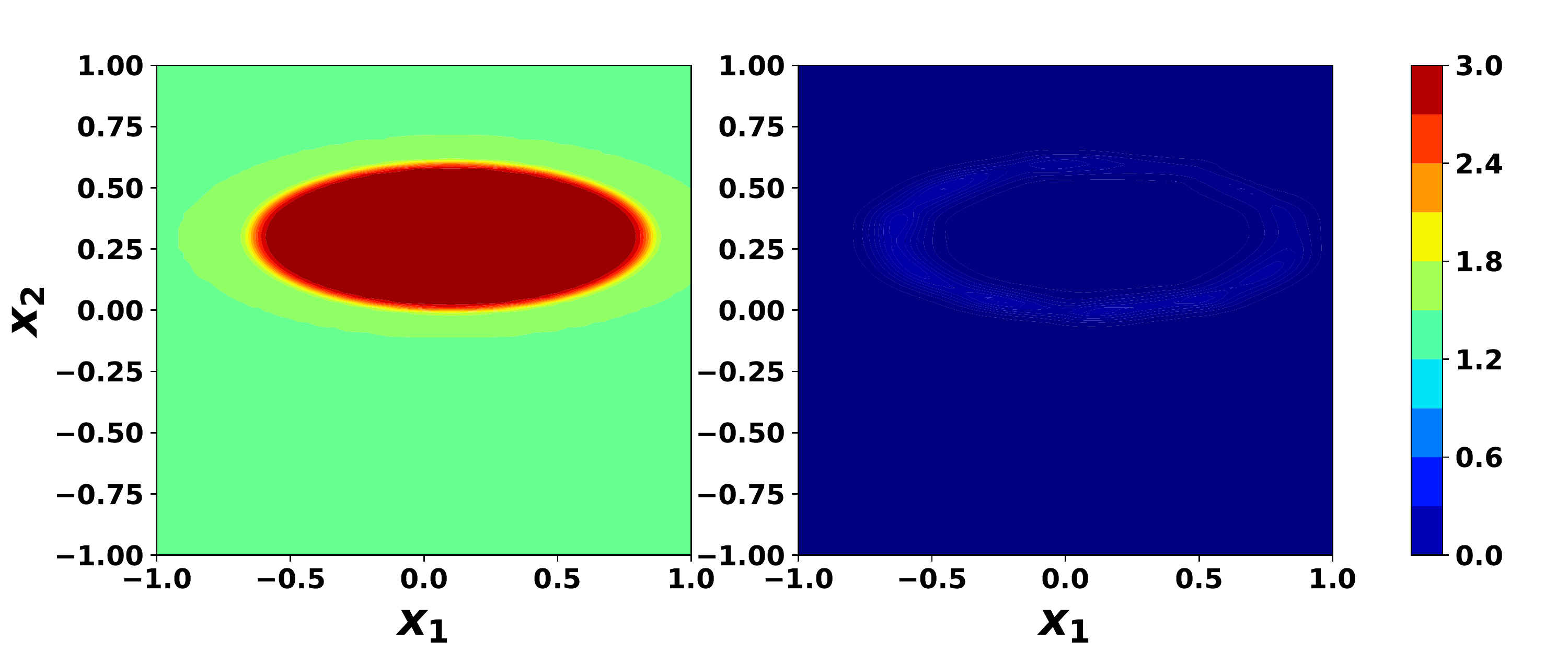}
\caption{$|\gamma-\gamma^*|$ of PINN (left) and IWAN (right) on less smooth $\gamma^*$}
\label{subfig:pinn_nonsmooth_gamma}
\end{subfigure}
\begin{subfigure}[b]{.32\textwidth}
\includegraphics[width=.875\textwidth]{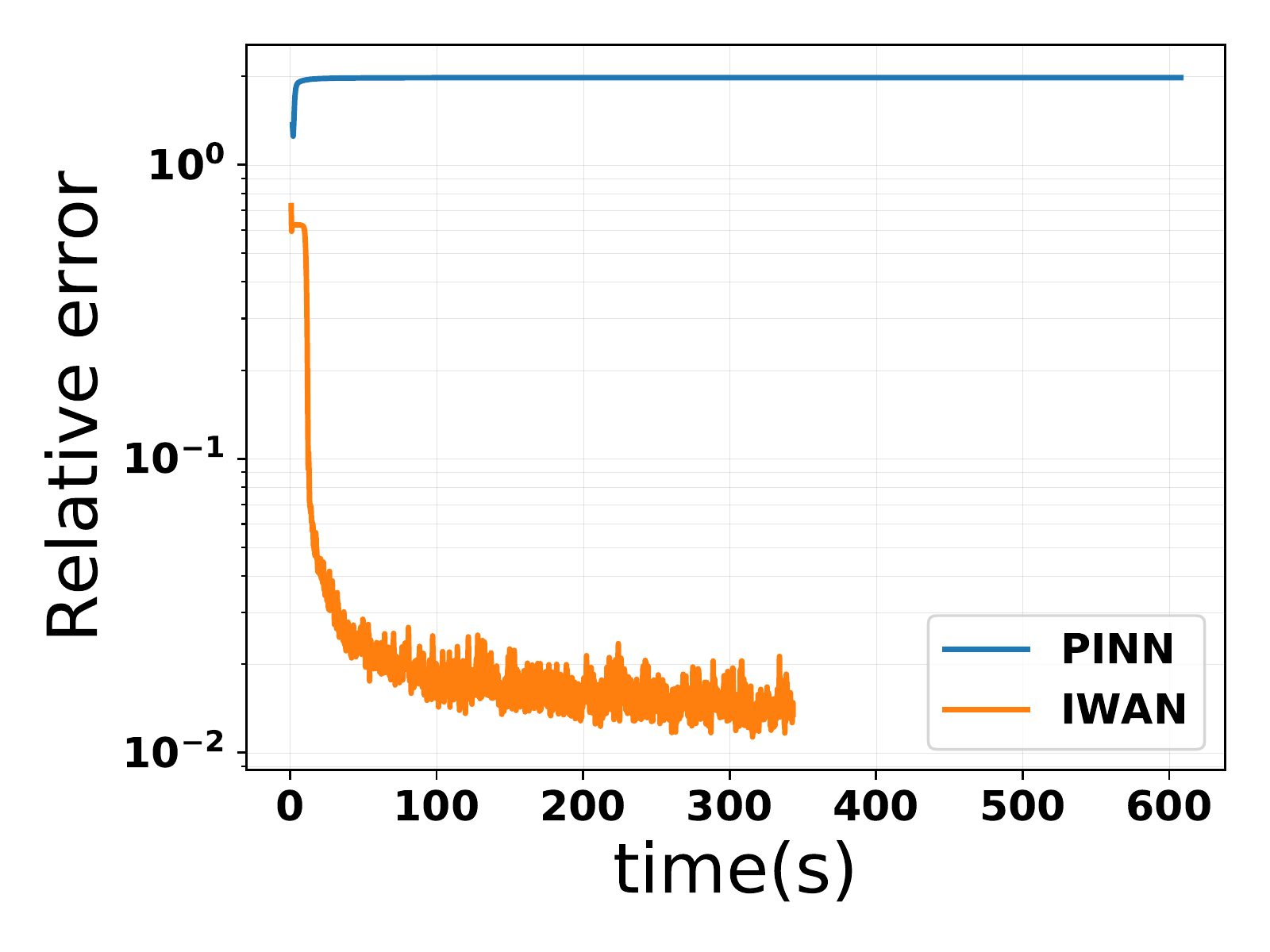}
\caption{Relative error vs time (s).}
\label{subfig:pinn_nonsmooth_time_20k_iter}
\end{subfigure}
\begin{subfigure}[b]{.6\textwidth}
\includegraphics[width=.875\textwidth]{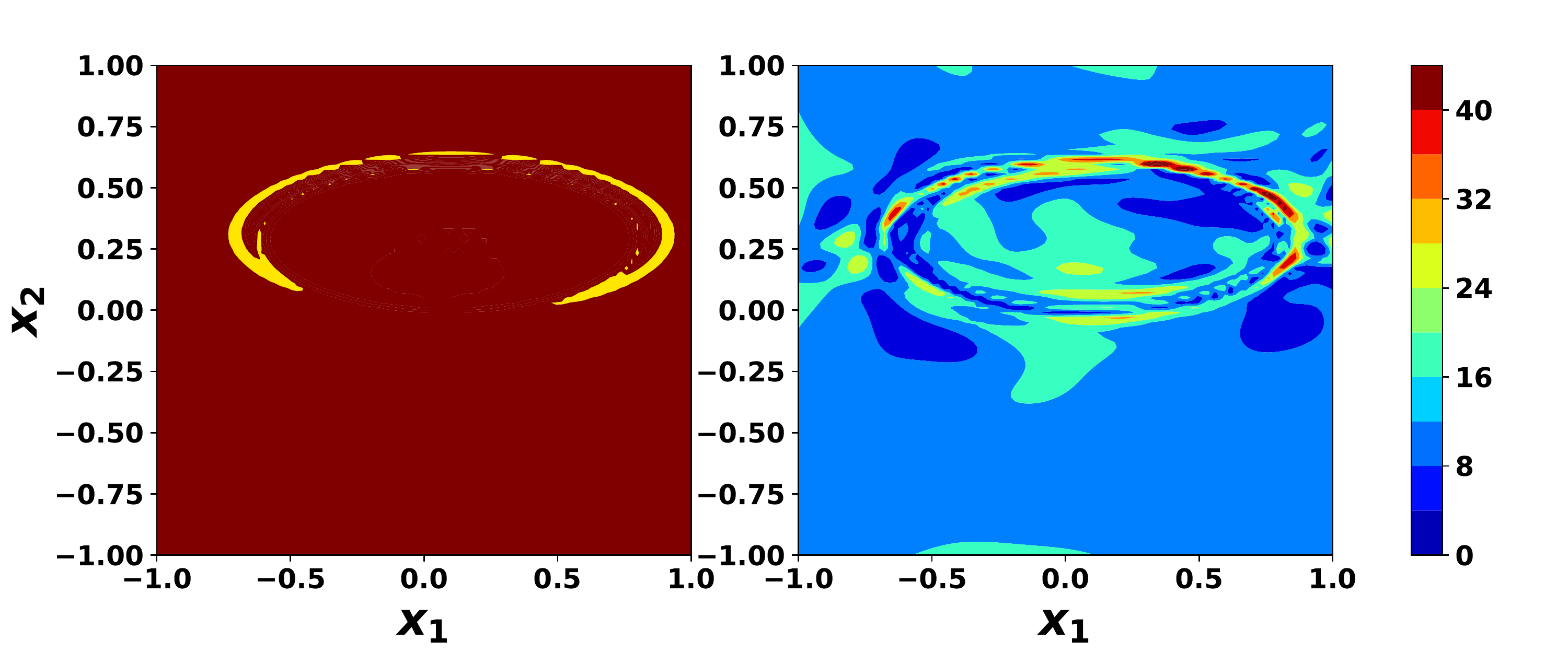}
\caption{$|f|$ (left) and $|-\nabla(\gamma\nabla u)-f|$ by the PINN (right)}
\label{subfig:pinn_nonsmooth_PDE_error}
\end{subfigure}
\begin{subfigure}[b]{.32\textwidth}
\includegraphics[width=.875\textwidth]{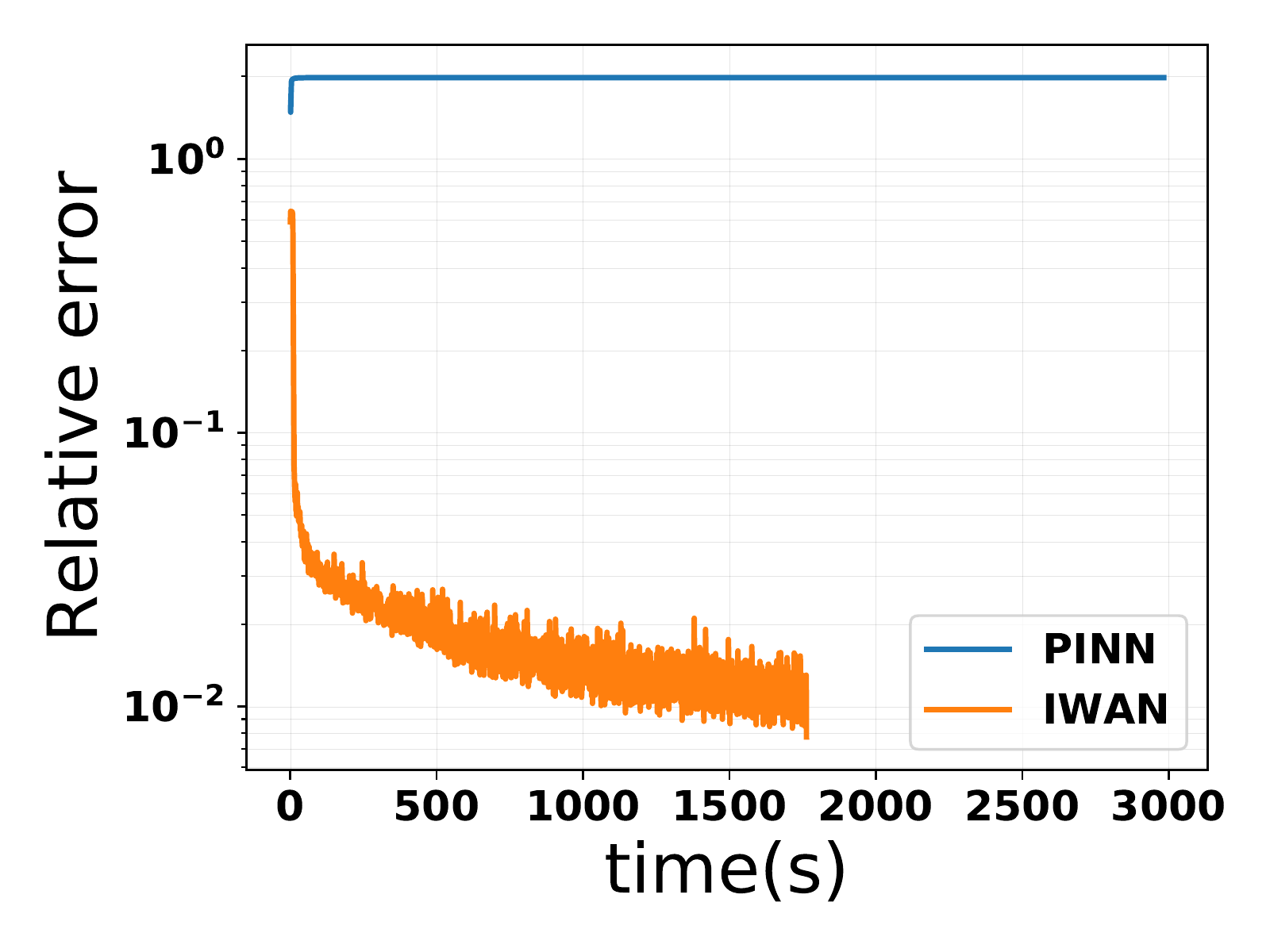}
\caption{Relative error vs time (s).}
\label{subfig:pinn_nonsmooth_time_100k_iter}
\end{subfigure}
\begin{subfigure}[b]{.3\textwidth}
\includegraphics[width=.875\textwidth]{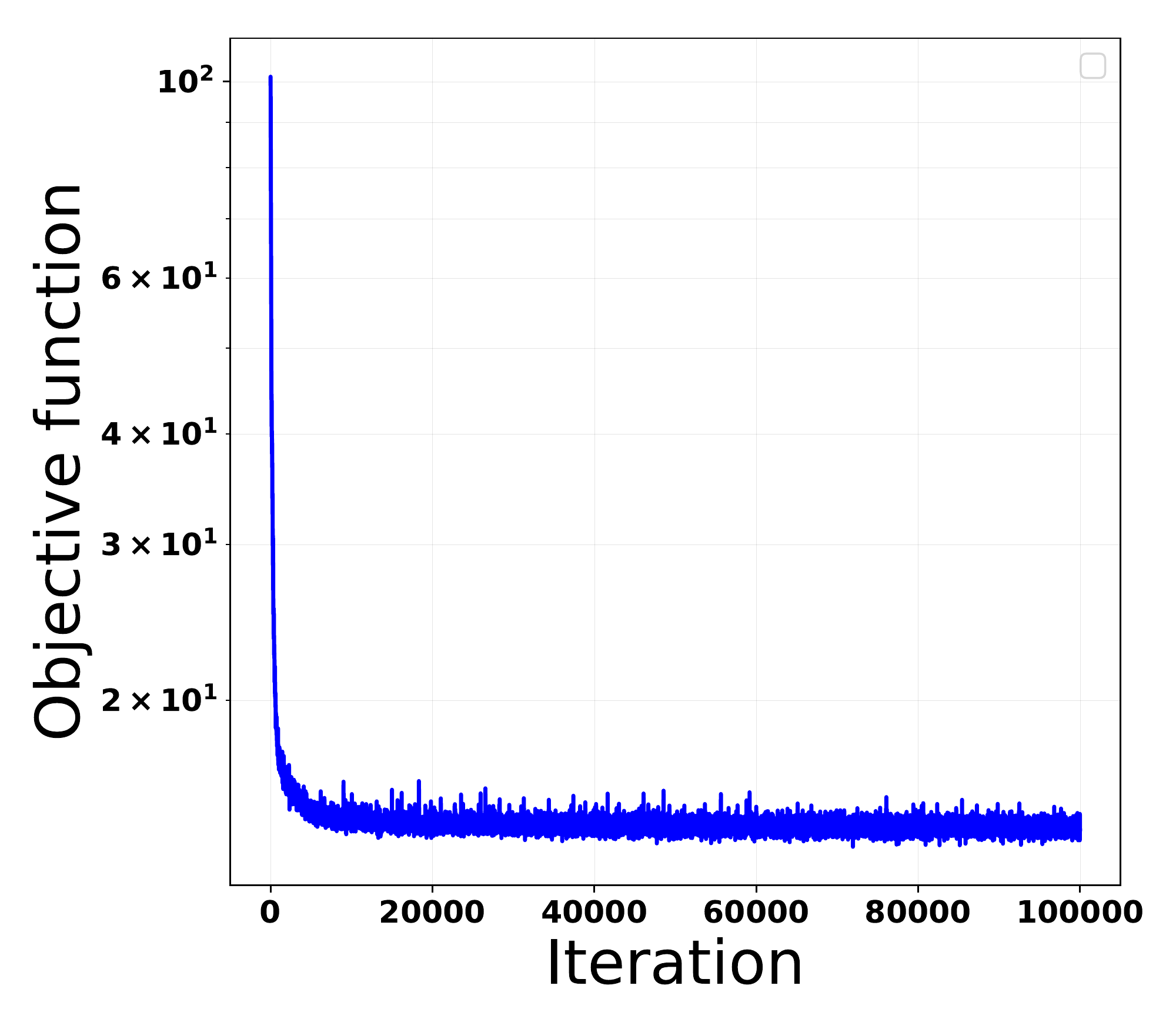}
\caption{Objective vs.~iter\\ by PINN}
\label{subfig:pinn_loss_nonsmooth}
\end{subfigure}
\begin{subfigure}[b]{.3\textwidth}
\includegraphics[width=.875\textwidth]{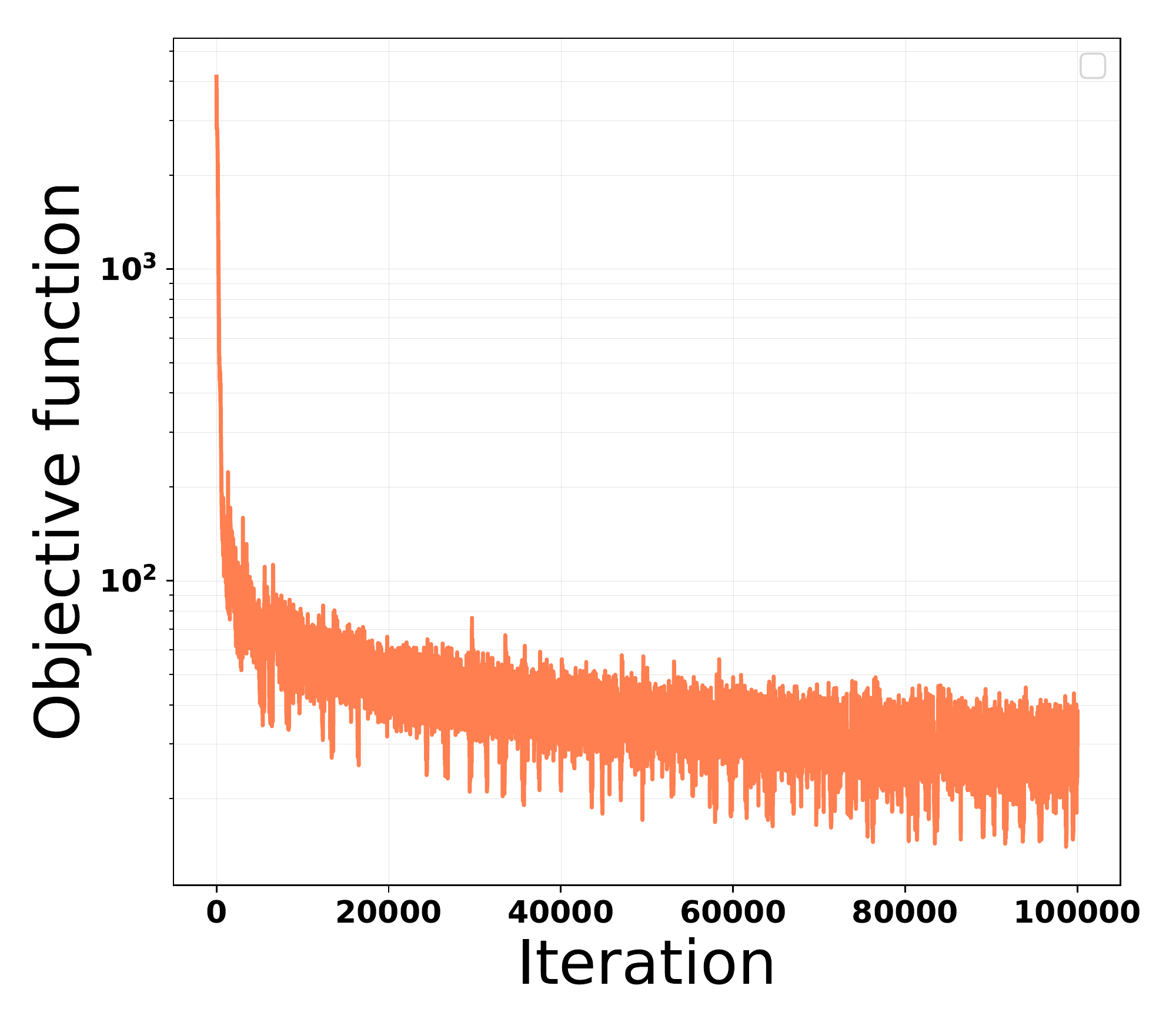}
\caption{Objective vs.~iter \\ by IWAN}
\label{subfig:iwan_loss_nonsmooth}
\end{subfigure}
\caption{Test 7 on the comparison of IWAN and PINN on the recovery smooth ((a) and (b)) and less smooth ((c)--(h)) conductivity $\gamma^*$. (a) Pointwise absolute error $|\gamma-\gamma^*|$ with $\gamma$ obtained by PINN (left) and the proposed method IWAN (right) for smooth $\gamma^*$. (b) Relative error versus time in seconds for smooth $\gamma^*$. (c) Pointwise absolute error $|\gamma-\gamma^*|$ with $\gamma$ obtained by PINN (left) and the proposed method IWAN (right) for less smooth, nearly piecewise constant $\gamma^*$. (d) Relative error versus time in seconds for 20,000 iterations for less smooth, nearly piecewise constant $\gamma^*$. (e) $|f|$ (left) and The map of $|-\nabla(\gamma\nabla u)-f|$ by the PINN (right) for the less smooth $\gamma^*$. (f) Relative error versus time in seconds for 100,000 iterations for less smooth, nearly piecewise constant $\gamma^*$. (g) Objective function value versus iteration number by PINN for nearly piecewise constant $\gamma^*$. (h) Objective function value versus iteration number by the proposed IWAN for nearly piecewise constant $\gamma^*$.}
\label{fig:con_dnn}
\end{figure}

\medskip
\noindent
\textbf{Test 7: Comparison with PINN.}
We compare the proposed method with a state-of-the-art method called physics-informed neural networks (PINNs) \cite{raissi2019physics}. PINN is also a deep-learning based method designed for solving forward problem as well as inverse problem for PDEs. PINN is based on the strong form of the PDEs, where the loss function in the minimization problem of PINN consists of the sum of squared errors in the violation of the PDE and the boundary condition at points sampled inside $\Omega$ and on $\partial \Omega$, respectively. In contrast, our method is based on the weak form of PDE which employs a test function and yields a min-max problem to better tackle singularities of the problem.
We first compare the proposed method and the PINN method in the problem in {Test 1}. Note that PINN was only applied to inverse conductivity problem with constant conductivity in \cite{raissi2019physics}, it is straightforward to extend this method to non-constant conductivity by also parameterizing the conductivity $\gamma$ as an additional deep neural network. In this problem, we take $N_r= 10,000, N_b=100*d$ with $d=5$. For the PINN method, we also parameterize $u(x)$ as a 9-layer fully-connected network with $20$ neurons in each hidden layer and tanh as activation functions in all hidden layers. For $\gamma$, we parameterized it by using the same network structure as that used for the proposed method in {Test 1}. We let the weight of the boundary term in the loss function is $1.0$ and use the builtin Adam optimizer of TensorFlow with learning rate $0.001$ to update the network parameters in PINN. For fair comparison, we also use the Adam optimizer with learning rate $0.001$ for updating $\theta$ and $\eta$ in the proposed method. The results after 20,000 iterations of both methods are given in Figure \ref{fig:con_dnn}. In Figure \ref{subfig:pinn_smooth_gamma}, we can see the error of $\gamma$ obtained by IWAN is much lower than that by PINN. This can also be seen from Figure \ref{subfig:pinn_smooth_time}, where the error decays very fast for the proposed IWAN. We tried a variety of network structures and parameter settings of PINN and obtain similar results.

We also compared the proposed method and PINN on the inverse conductivity problem in {Test 2}, where the ground truth conductivity $\gamma^*$ is less smooth and nearly piecewise constant. We use the same parameter settings for both methods as above except for the network structure of $\gamma$, which follows the one in {Test 2}. The conductivity $\gamma$ recovered by PINN and IWAN and the progresses of their relative error versus computation time (in seconds) are given in Figures \ref{subfig:pinn_nonsmooth_gamma} and \ref{subfig:pinn_nonsmooth_time_20k_iter}, respectively.
From Figure \ref{subfig:pinn_nonsmooth_time_20k_iter}, it appears that PINN cannot get close to the ground truth $\gamma^*$ within 20,000 iterations. Therefore, we rerun PINN for 100,000 iterations, and plot the relative error versus computation time in Figure \ref{subfig:pinn_nonsmooth_time_100k_iter}, from which it seems that PINN still cannot converge to the desired solution. However, the result obtained by PINN does satisfy the PDE closely, as shown in Figure \ref{subfig:pinn_nonsmooth_PDE_error}: the difference between the two sides of PDE (left), i.e., $|-\nabla(\gamma\nabla u)-f|$, is much smaller than $|f|$ (right), but PINN cannot capture the irregularities and singularities of the solution since it is based on the strong form of the PDE.
In contrast, IWAN can overcome this issue and recover the weak solution properly.
We also show the objective function value of PINN and IWAN in Figures \ref{subfig:pinn_loss_nonsmooth} and \ref{subfig:iwan_loss_nonsmooth}, respectively (Note that the objective function $L(\theta,\eta)$ in IWAN is defined as $E(\theta,\eta)+\beta \Lbdry(\theta)$ for min-max optimization, and the objective function of PINN is for minimization only and hence different from IWAN). 

\begin{figure}[t]
    \centering
    \includegraphics[width=0.3\textwidth]{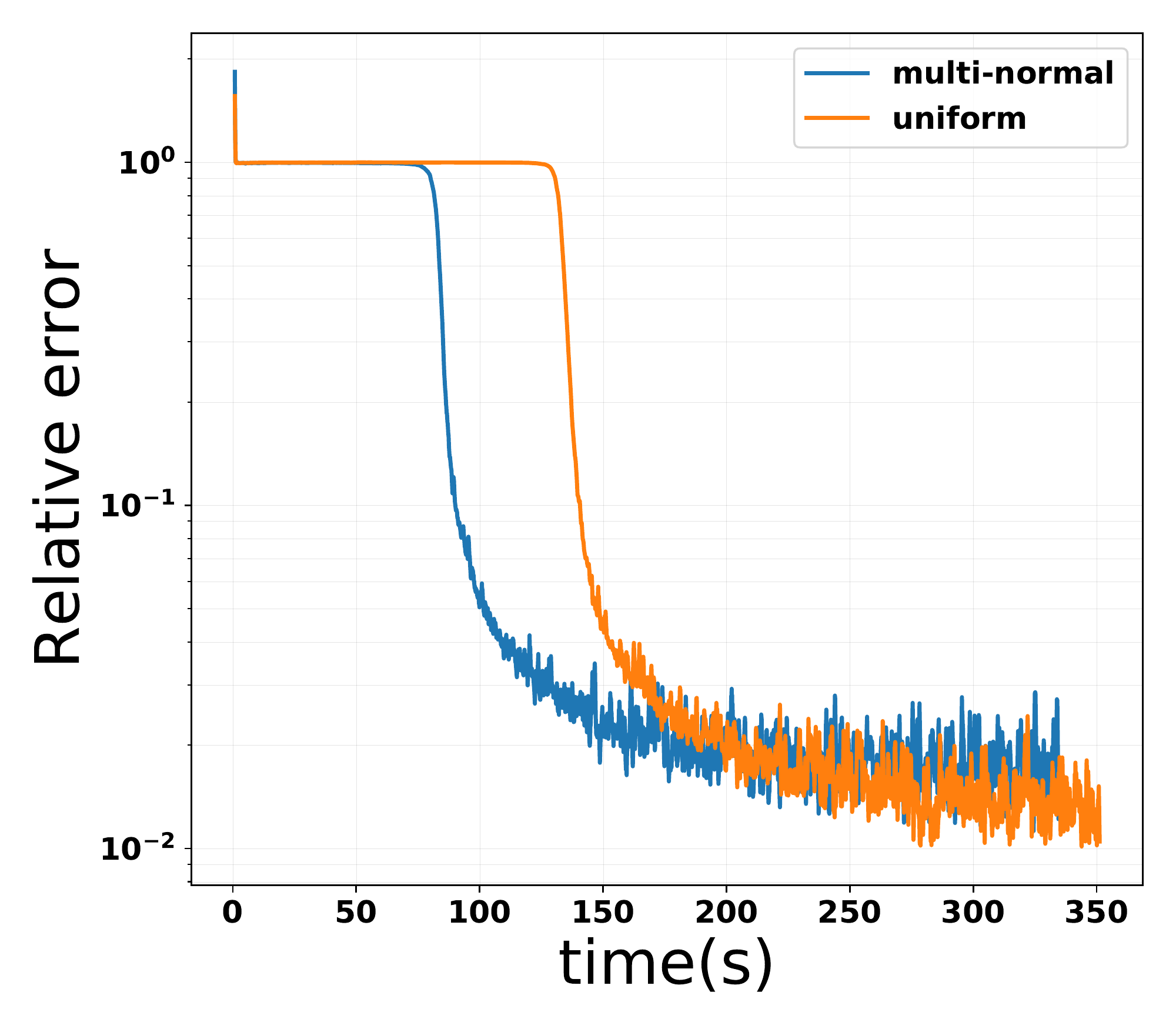}
    \caption{Test 8 result on the difference of relative errors versus computation time (s) using collocations points $\{x_r^{(i)}\in \Omega: i\in[N_r]\}$ sampled from uniform distribution (orange) and adaptive multivariate normal distribution (blue).}
    \label{fig:important_sample}
\end{figure}

\medskip
\noindent
\textbf{Test 8: Efficiency improvement using important sampling.}
As shown in Lemma \ref{lem:integral}, adaptive sampling may reduce the variance of the sample-based approximation of integrals, which in turn can improve the convergence of stochastic gradient descent. To demonstrate this, we consider the inverse conductivity problem on $\Omega=(-1,1)^d\subset \mathbb{R}^d$ with problem dimension $d=5$, ground truth conductivity distribution $\gamma^*(x)=2\exp(-|x-c|_{\Sigma}^2/2) $, where $\Sigma=\mathrm{diag}(4.0,100.0,0,0,0)$, $c_1=(-0.2,0.2,0,0,0)$.
We set $u_n(x)=-2\sin(|x|^2)(x_1,x_2,\cdots,x_d)\cdot\vec{n}$, $u_b=\cos(|x|^2)$ on $\partial\Omega$ and $f(x)= 8\sum^{2}_{i=1}\Sigma_{ii}(x_i-c_i)x_i\sin(|x|^2)\exp(-|x-c|^2_{\Sigma}/2)+2\exp(-|x-c|^2_{\Sigma}/2)*(2d\sin(|x|^2)+4|x|^2\cos(|x|^2))$.
For the parameter setup, we use the same setup as that in {Test 7}. Then we solve this inverse problem using the proposed method IWAN with points in the domain $\Omega$ sampled from the uniform distribution as above and also a multivariate normal distribution respectively. Specifically, to obtain multivariate normal samples, we first sample $N_r$ points of $(x_1,x_2)$ from the multivariate normal distribution with mean value $\mu=(-0.2, 0.2)$ and inverse covariance matrix $\Sigma=\text{diag}(1.0, 25.0)$ (points outside of $\Omega$ is discarded), and then draw each of the remaining coordinates randomly from interval $(-1,1)$ independently.
The result was shown in the figure \ref{fig:important_sample}. The progress of relative error versus computation time (in second) for $20,000$ iterations is shown in Figure \ref{fig:important_sample}, which shows that the convergence using adaptive multivariate normal distribution is faster than that with uniform distribution.

\subsection{Empirical Robustness Analysis}
We conduct a series of experiments to evaluate the robustness of Algorithm \ref{alg:iwan} in terms of network structure (number of layers and neurons) and the number of sampled collocation points.

\begin{figure}[t]
\centering
\begin{subfigure}[b]{.23\textwidth}
\includegraphics[width=\textwidth]{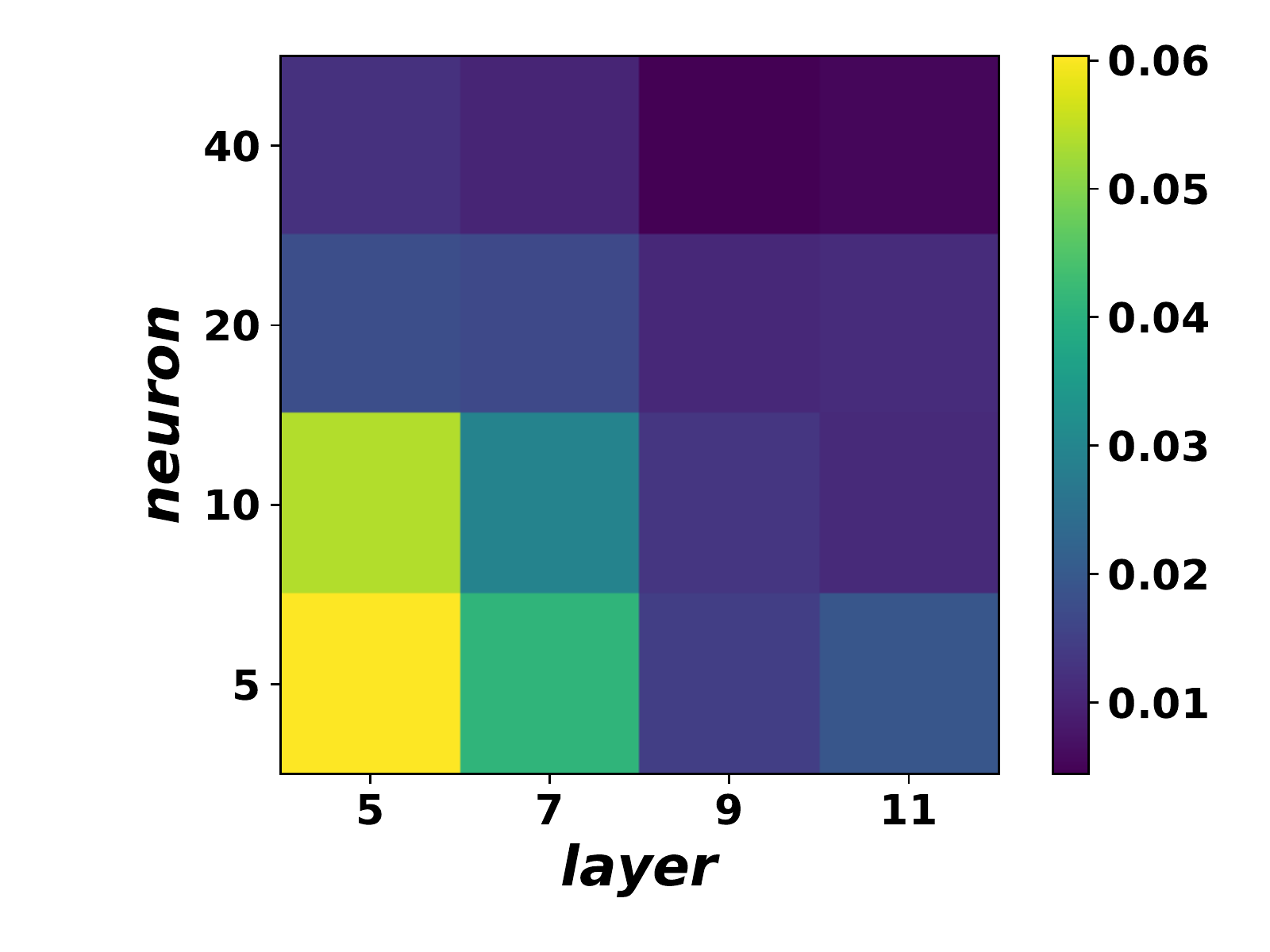}
\caption{Relative error}
\label{subfig:error_layer_neuron}
\end{subfigure}
%
\begin{subfigure}[b]{.23\textwidth}
\includegraphics[width=\textwidth]{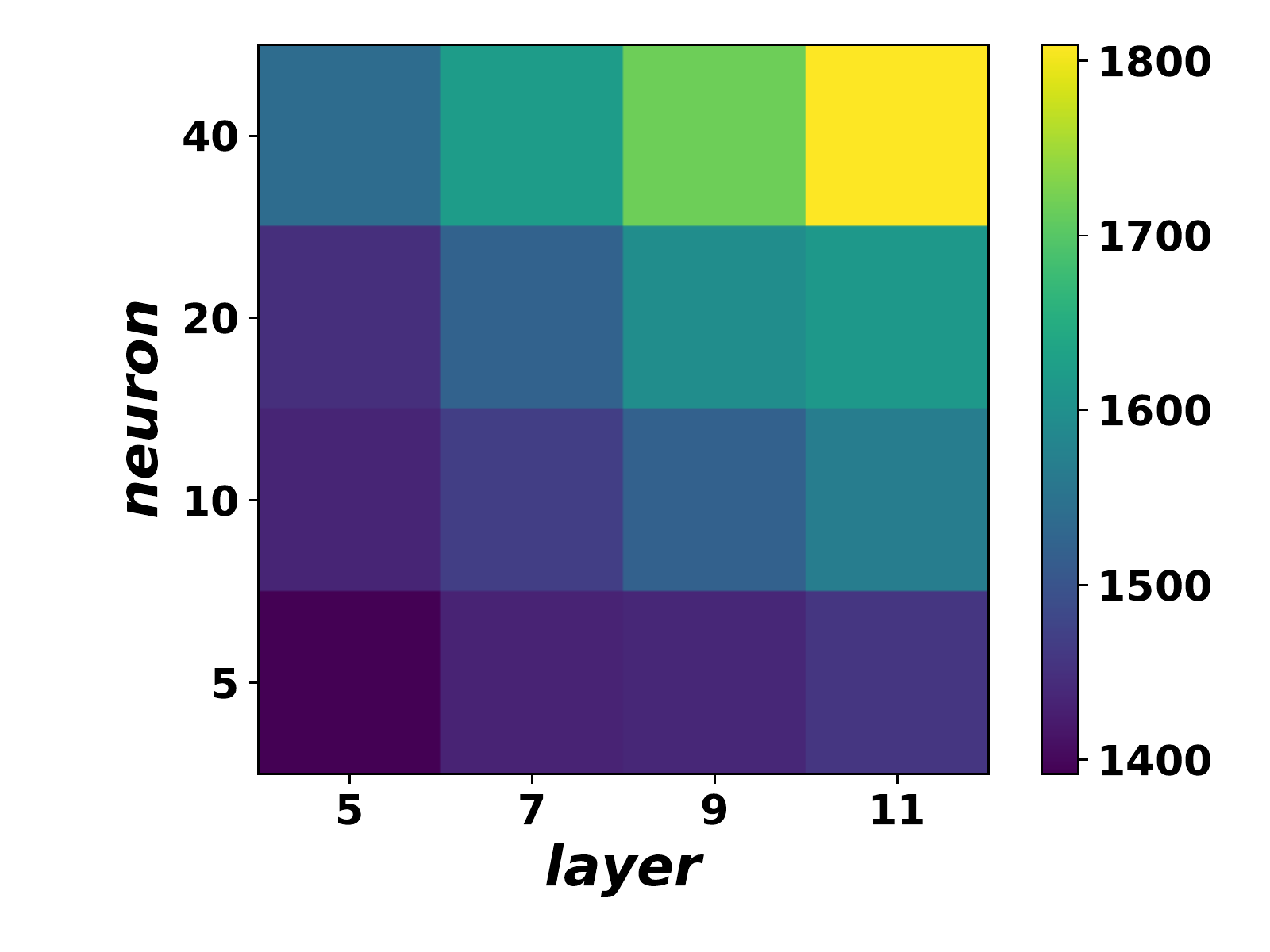}
\caption{Time(s)}
\label{subfig:time_layer_neuron}
\end{subfigure}
\begin{subfigure}[b]{.23\textwidth}
\includegraphics[width=\textwidth]{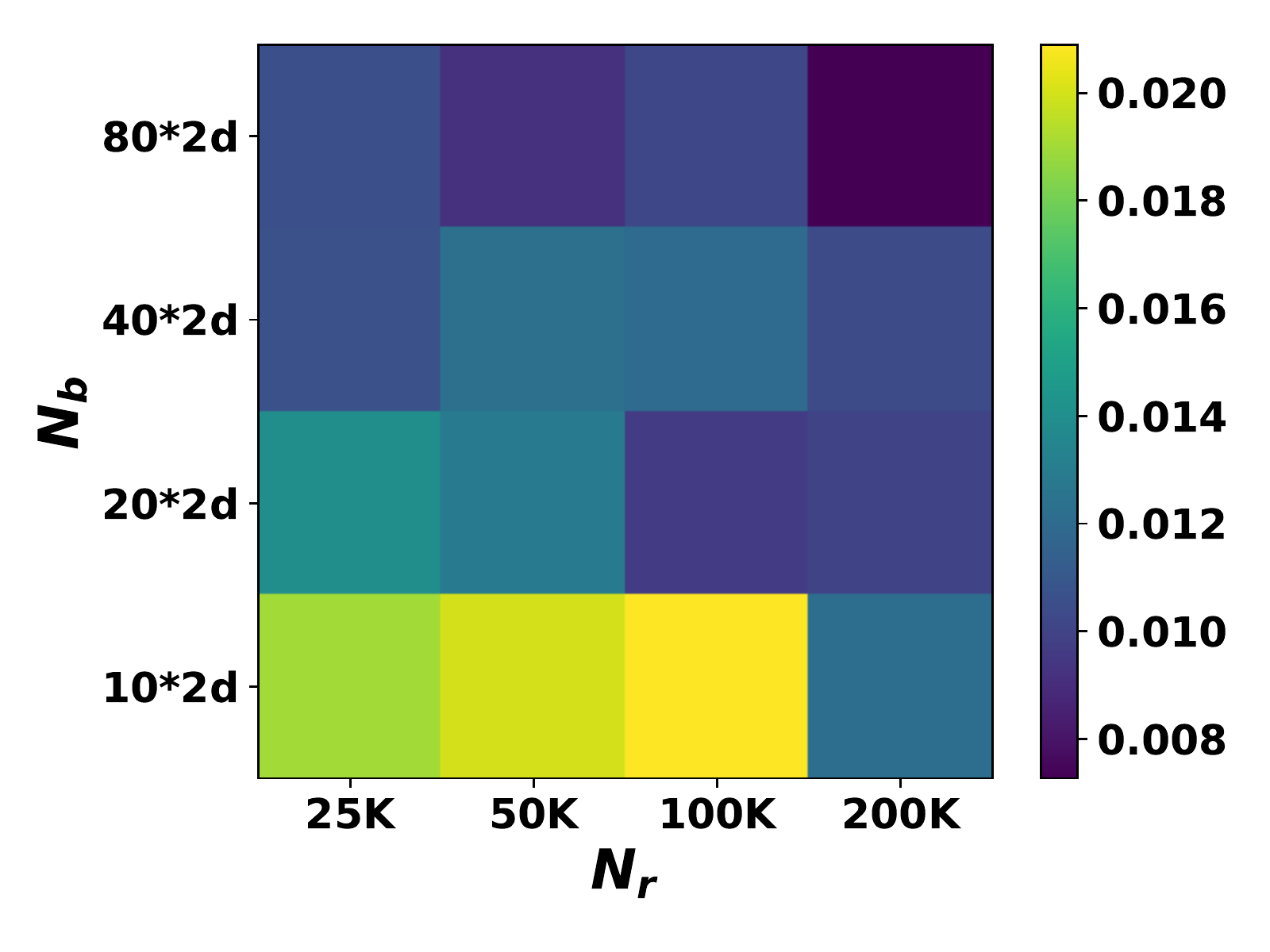}
\caption{Relative error}
\label{subfig:error_Nr_Nb}
\end{subfigure}
%
%
\begin{subfigure}[b]{.23\textwidth}
\includegraphics[width=\textwidth]{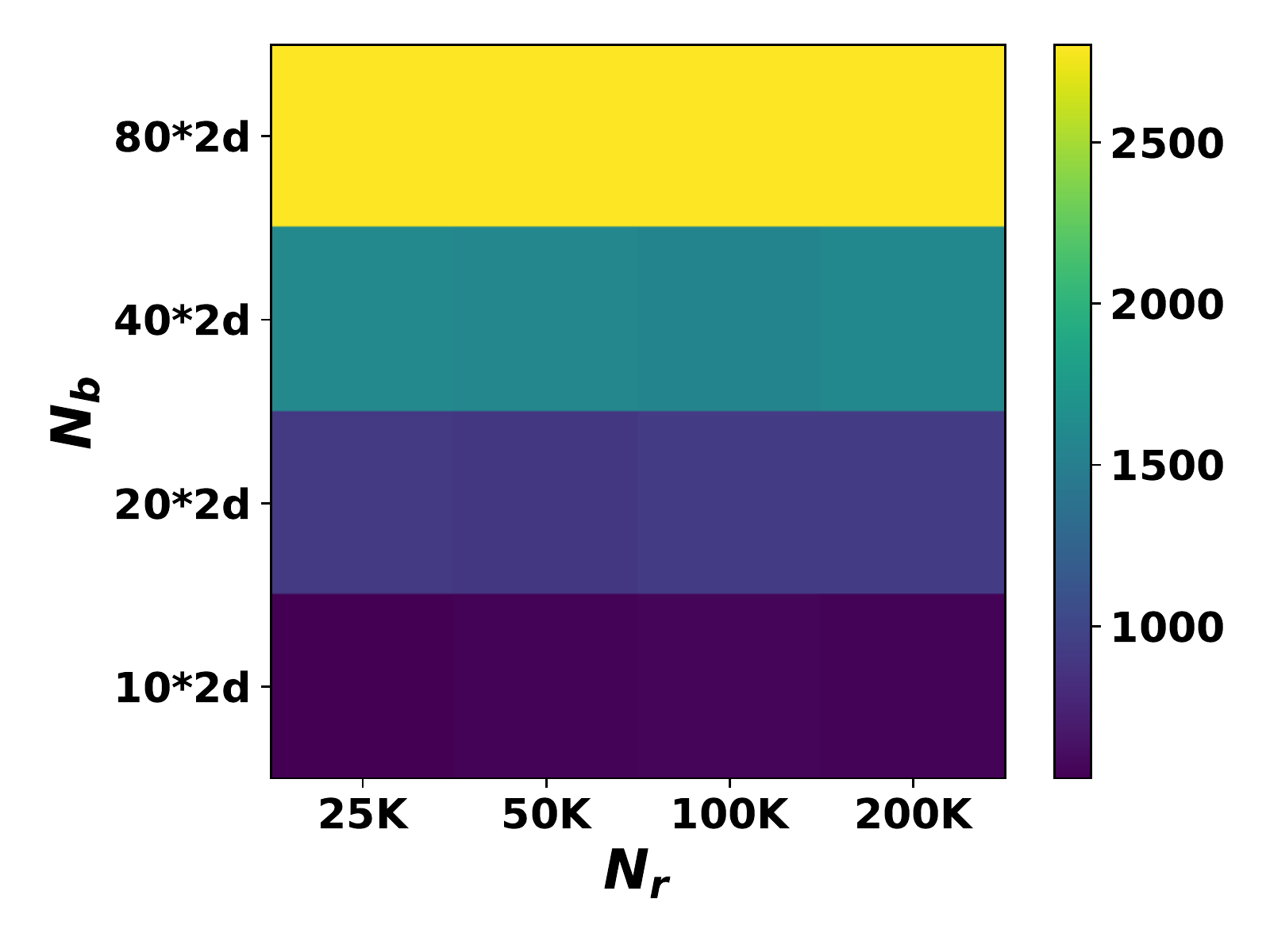}
\caption{Time(s)}
\label{subfig:time_Nr_Nb}
\end{subfigure}
\caption{(a),(b) Test 9 result on relative error of recovered conductivity $\gamma_\theta$ and the corresponding running time  using various combinations of $(K,d')$, where $K$ is the layer number and $d'$ is the per-layer neuron number. (c),(d) Test 10 result on relative error of recovered conductivity $\gamma_\theta$ and the corresponding running time using various combination of $(N_r,N_b)$, where $N_r$ is the number of sampled collocation points inside the region $\Omega$ and $N_b$ is the number of those on the boundary $\partial \Omega$.}
\label{fig:error_blcok}
\end{figure}

\begin{figure}[bht]
\centering
\begin{subfigure}[b]{.3\textwidth}
\includegraphics[width=\textwidth]{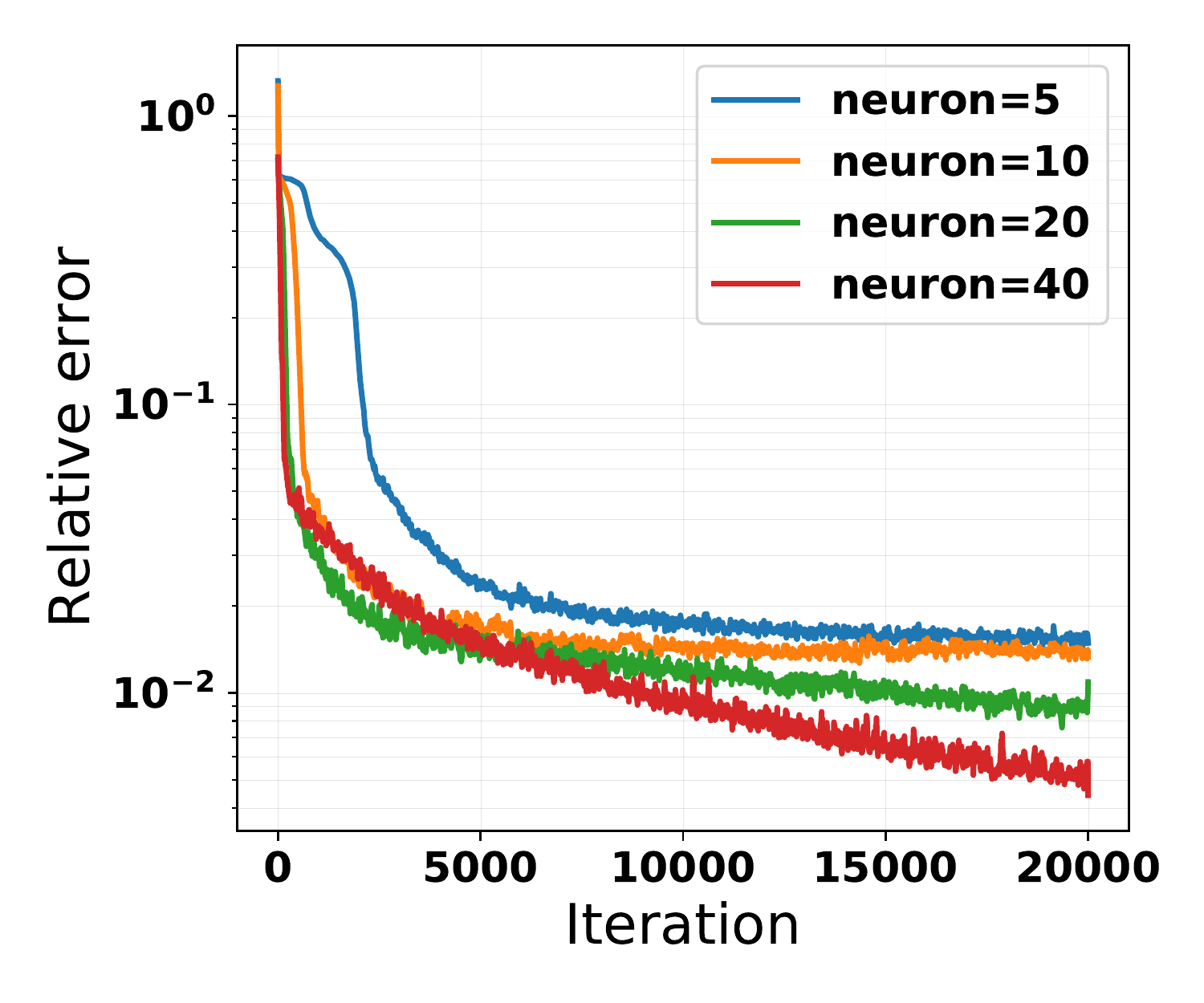}
\caption{Relative error for $K=9$.}
\label{subfig:fix_layer}
\end{subfigure}
\begin{subfigure}[b]{.3\textwidth}
\includegraphics[width=\textwidth]{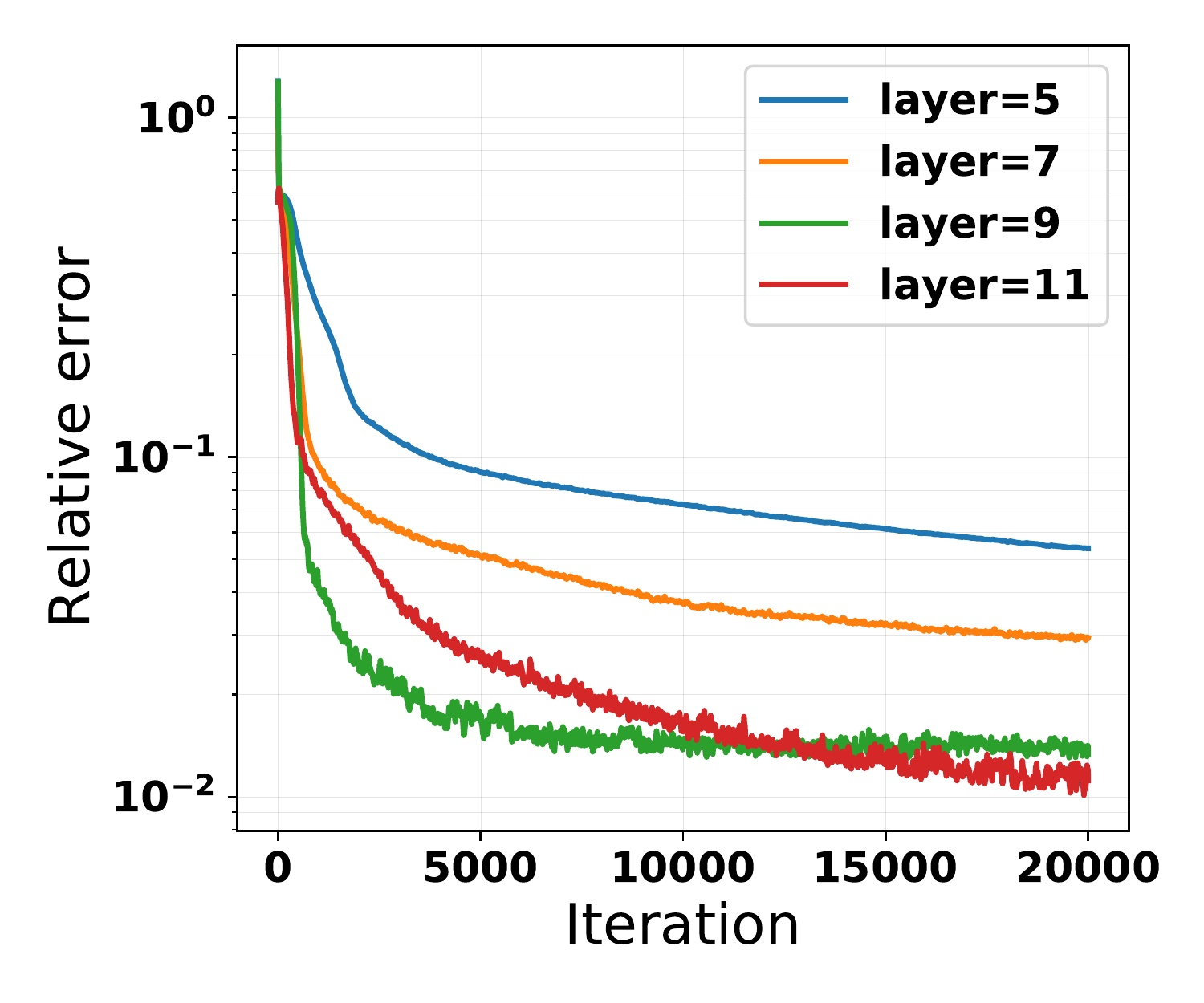}
\caption{Relative error for $d'=10$.}
\label{subfig:fix_neuron}
\end{subfigure}
\caption{Test 9 result on the effect of network structure. (a) relative error of $\gamma_\theta$ versus iteration number with fixed layer number $K=9$ and varying per-layer neuron number $d'=5,10,20,40$; (b) relative error of $\gamma_\theta$ versus iteration number with fixed per-layer neuron number $d'=10$ and varying layer number $K=5,7,9,11$.}
\label{fig:layer_neuron}
\end{figure}

\medskip
\noindent
\textbf{Test 9: Network structure.}
In this experiment, we test the performance of Algorithm \ref{alg:iwan} with different network structures, i.e., the layer number (network depth) $K$ and the per-layer neuron number (network width) $d_k$.
We test different combinations of $K$ and $d'$ (in each combination we set $d_k=d'$ for all $k=1,\dots,K-1$).
More specifically, we apply Algorithm \ref{alg:iwan} to the inverse conductivity problem \eqref{eq:eit} in Test 2 above with problem dimension $d=5$ and a total of 16 combinations $(K,d')$ with $K=5,7,9,11$ and $d'=5,10,20,40$.
{For each combination $(K,d')$, we run Algorithm \ref{alg:iwan} for $20,000$ iterations and plot the relative error of $\gamma_\theta$ in Figure \ref{subfig:error_layer_neuron}} and the corresponding running time (in seconds) in Figure \ref{subfig:time_layer_neuron}. The exact values of errors and running times are present in Table \ref{tab:error_time_layer_neuron} in Appendix \ref{app:error_time_table}.
Based on Figure \ref{subfig:error_layer_neuron}, it seems that deeper (larger $K$) and/or wider (larger $d'$) neural networks yield lower reconstruction error (but at the expense of higher per-iteration computational cost).
For fixed layer number $K=9$, we show the progress of relative error versus iteration number with varying per-layer neuron number $d'$ in Figure \ref{subfig:fix_layer}.
Similarly, for fixed per-layer neuron number $d'=10$, we also show the progress of relative error versus iteration number with varying layer number $K$ in Figure \ref{subfig:fix_neuron}.
These figures also suggest that larger $K$ and $d'$ yield better accuracy, although the per-iteration computational cost also increases and it may take more iterations to converge.
%
%

\begin{figure}[ht]
\centering
\begin{subfigure}[b]{.3\textwidth}
\includegraphics[width=\textwidth]{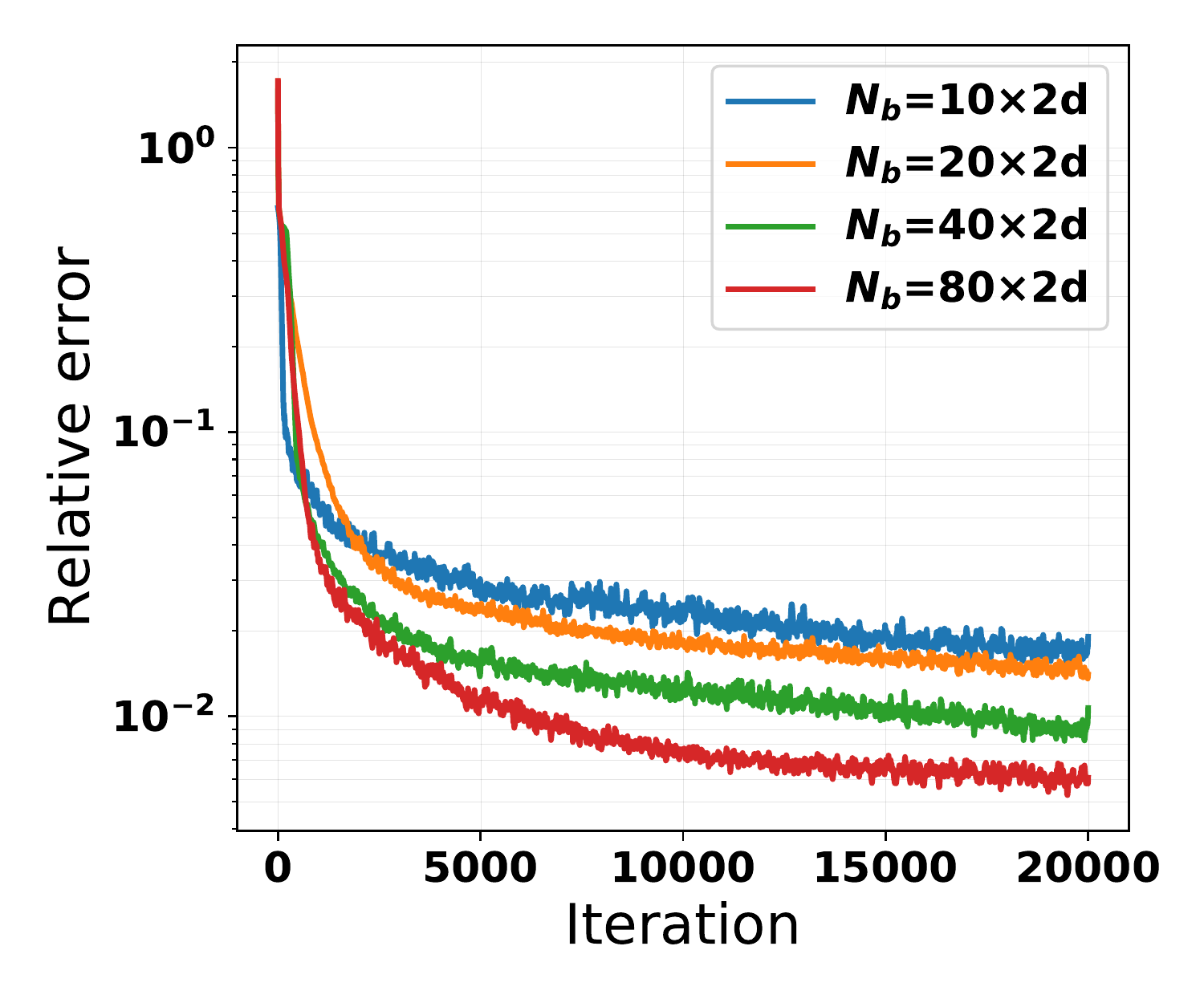}
\caption{Relative error for $N_r=25$K.}
\label{subfig:fix_Nr}
\end{subfigure}
\begin{subfigure}[b]{.3\textwidth}
\includegraphics[width=\textwidth]{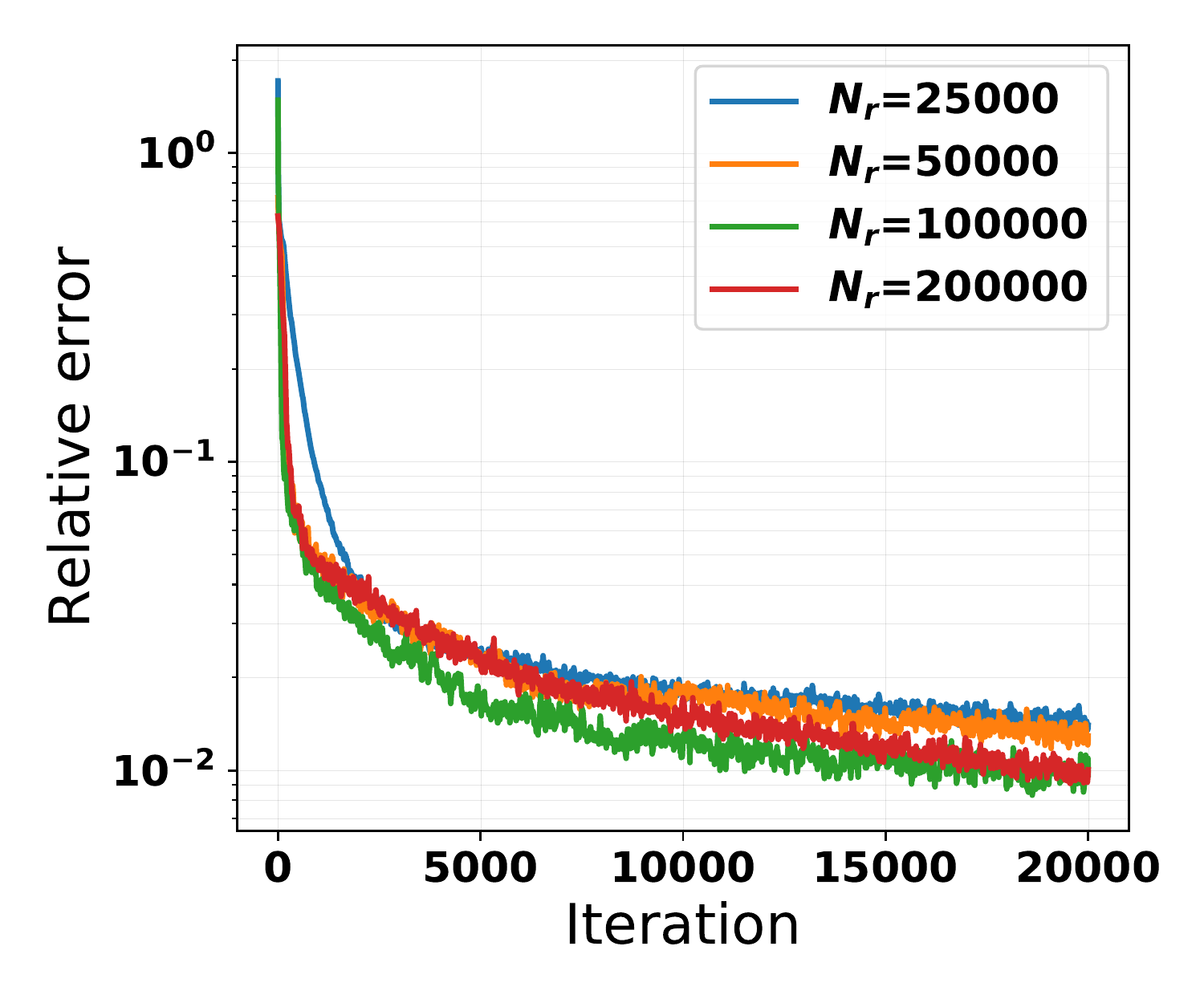}
\caption{Relative error for $N_b=200$.}
\label{subfig:fix_Nb}
\end{subfigure}
\caption{Test 10 result on the different numbers of collocation points $N_r$ and $N_b$. (a) relative error of $\gamma_\theta$ versus iteration number with fixed $N_r=25$K and varying $N_b$; (b) relative error of $\gamma_\theta$ versus iteration number with fixed $N_b=200$ and varying $N_r$.}
\label{fig:Ndm_Nbd}
\end{figure}

\medskip
\noindent
\textbf{Test 10: Number of sampled collocation points.}
In Section \ref{sec:iwan}, we showed that the number $N$ of sampled collocation points affects the variance of the integral estimator, so that the variance reduces at the order of $O(1/N)$.
In this experiment, we test the empirical effect of the collocation point numbers $N_r$ in the region $\Omega$ and $N_b$ on the boundary $\partial\Omega$ in Algorithm \ref{alg:iwan} for the same inverse conductivity problem \eqref{eq:eit} of dimension $d=5$ in Test 2 above.
We choose different combinations of $(N_r,N_b)$ for $N_r=25\text{K},50\text{K},100\text{K},200\text{K}$ (K=1,000) and $N_b=(10\times 2d, 20\times 2d, 40\times 2d, 80\times 2d)$, and keep all other parameters in Test 2 unchanged.
{We run Algorithm \ref{alg:iwan} for 20,000 iterations, and plot the final relative error of $\gamma_\theta$ in Figure \ref{subfig:error_Nr_Nb} and the corresponding running time (in seconds) in Figure \ref{subfig:time_Nr_Nb}. The exact values of errors and running times are present in Table \ref{tab:error_time_Nr_Nb} in Appendix \ref{app:error_time_table}.}
%
We also plot the progress of relative error of $\gamma_\theta$ versus iteration number for fixed $N_r=25$K and varying $N_b$ in Figure \ref{subfig:fix_Nr}, and that for fixed $N_b=20\times 2d=200$ and varying $N_r$ in Figure \ref{subfig:fix_Nb}.
The results in Figure \ref{subfig:error_Nr_Nb} and Figure \ref{fig:Ndm_Nbd} show that larger amounts of collocation points can generally improve accuracy of the reconstruction. 
%
%
%

%

\begin{figure}[t!]
\centering
\begin{subfigure}[b]{\textwidth}
\centering
\includegraphics[width=0.3\linewidth]{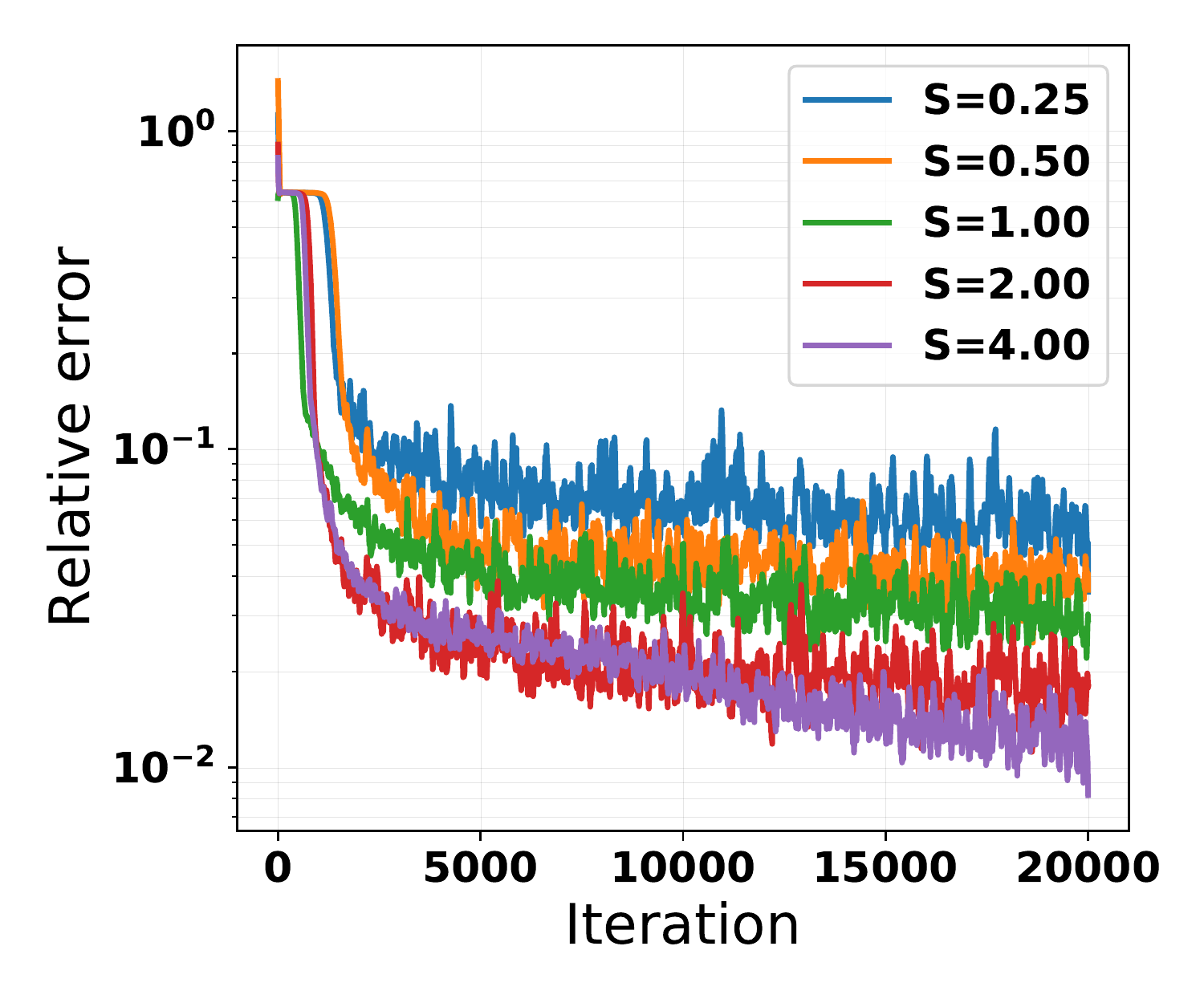}
\includegraphics[width=0.3\linewidth]{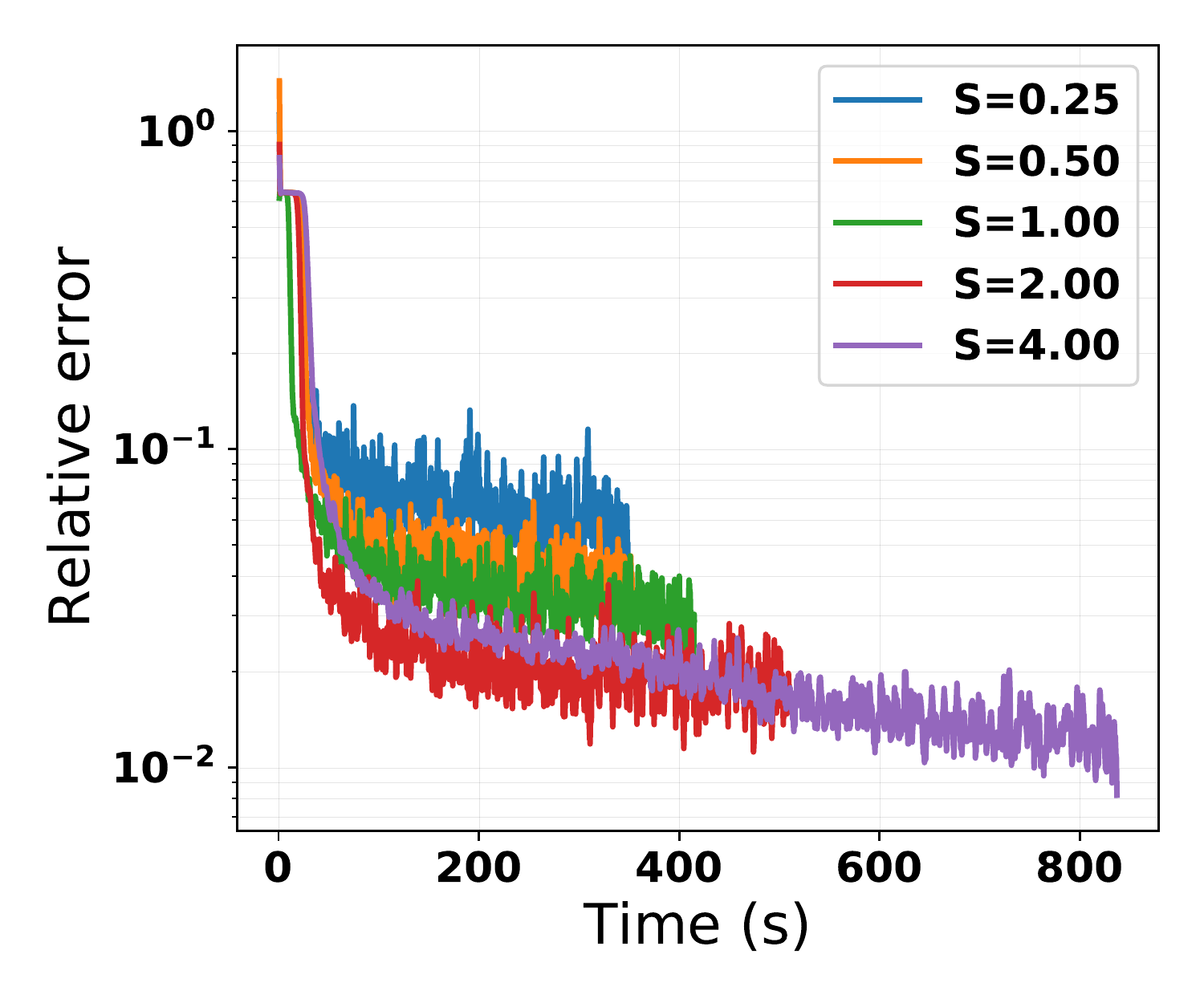}
\caption{Relatives vs. iteration (left) and time (right)}
\label{subfig:scale_err_time}
\end{subfigure}
\begin{subfigure}[b]{.3\textwidth}
\includegraphics[width=\textwidth]{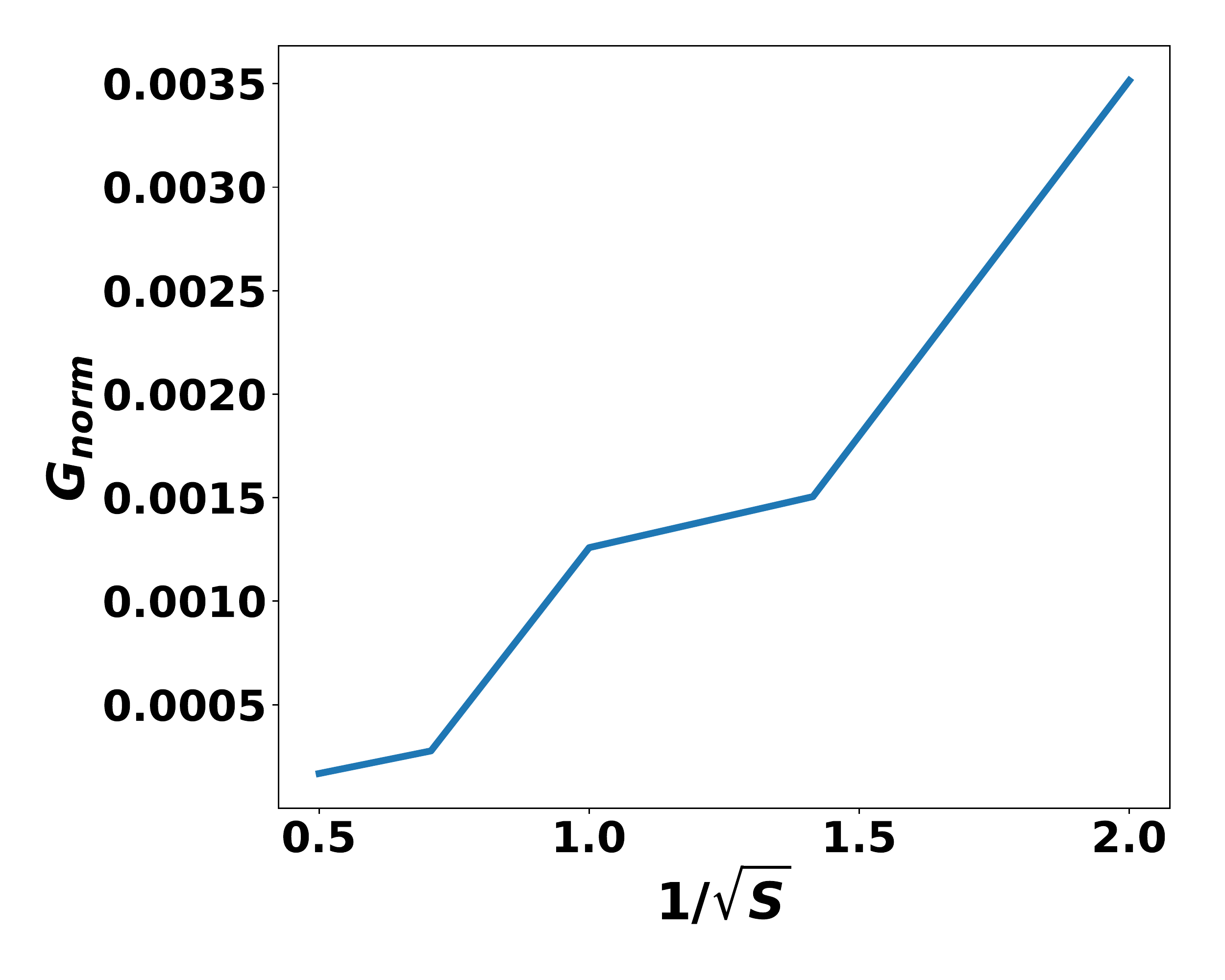}
\caption{$G_{norm}$ vs. $1/\sqrt{S}$}
\label{subfig:scale_G}
\end{subfigure}
\begin{subfigure}[b]{.3\textwidth}
\includegraphics[width=\textwidth]{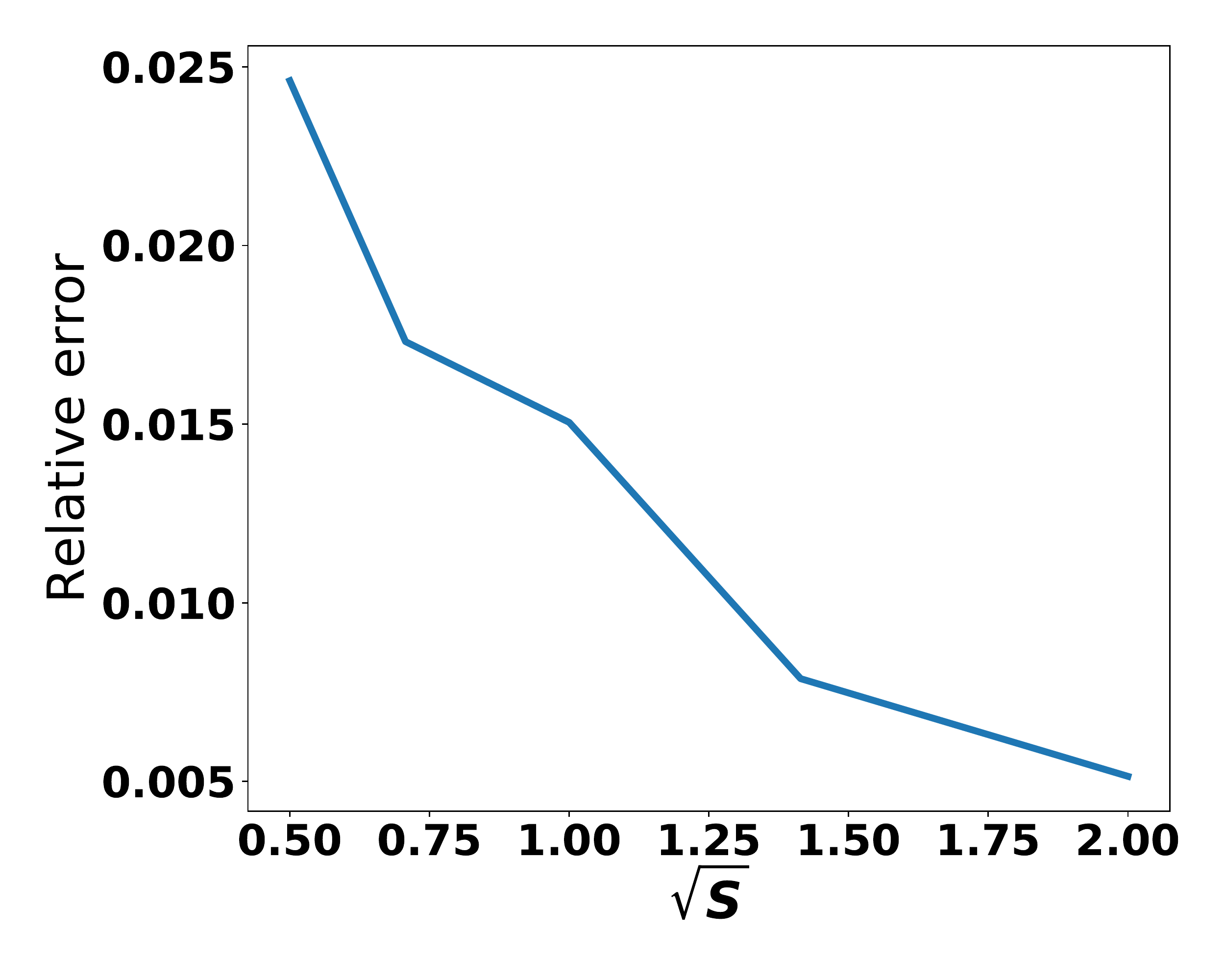}
\caption{Relative error vs. $\sqrt{S}$}
\label{subfig:scale_error}
\end{subfigure}
\begin{subfigure}[b]{.3\textwidth}
\includegraphics[width=\textwidth]{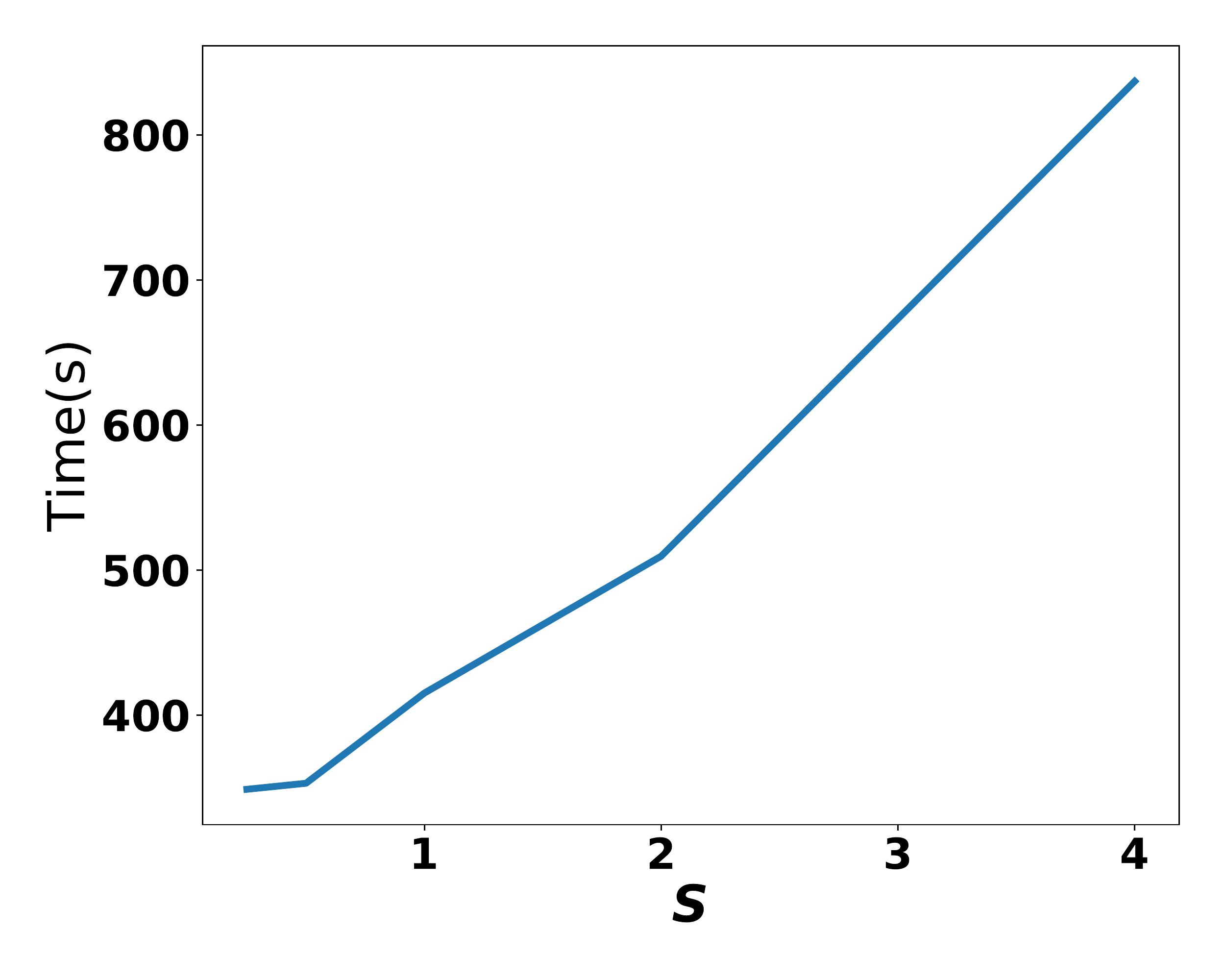}
\caption{Time vs. $S$}
\label{subfig:scale_time}
\end{subfigure}
\caption{Test 11 result on the effect of collocation points. (a) relative error of $\gamma_\theta$ versus iteration number (left) and running time (right) with collocation point numbers $N_r=S\times 25$K and $N_b=S\times 200$ with varying $S$; (b) The value of $G_{norm}$ versus $1/\sqrt{S}$; (c) Relative error versus $\sqrt{S}$; (d) Running time (in seconds) versus $S$.}
\label{fig:scale}
\end{figure}

\medskip
\noindent
\textbf{Test 11: Relation between relative error, gradient value and sample points.}
We conduct an experiment to show the relation between the relative error, the expected value of gradient mapping, and the numbers of collocation points $N_r$ and $N_b$. In this experiment, we use the Adam optimizer for $u_{\theta}, \gamma_{\theta}$ with learning rate $\tau_{\theta}=0.001$ and use the Adagrad optimizer for $\varphi_{\eta}$ with learning rate $\tau_{\eta}=0.008$ and $J_{\eta}=1$. We keep other settings unchanged as that used for case $N_r =25$K and $N_b=20\times 2d$ in {Test 10}. Then, we solve problem \eqref{eq:eit} of dimension $d=5$ using the proposed algorithm with $N_r= S\times 25$K and $N_b= S\times 200$, where $S$ takes value in $\{0.25, 0.5, 1.0, 2.0, 4.0\}$. After $J=20,000$ iterations, we evaluate $G_{norm}\coloneqq \sqrt{\min_{1\le j\le J}\mathbb{E}[|\Gcal(\theta_j)|^2]}$ (expectation is approximated by empirical average of 5 runs), 
the relative error, and the running time for each value of $S$, and plot their values versus $S$ in Figures \ref{subfig:scale_G}, \ref{subfig:scale_error}, \ref{subfig:scale_time}, respectively.
From Figure \ref{subfig:scale_G}, we can
see that $G_{norm}$ is approximately proportional to $1/\sqrt{S}$. Figures \ref{subfig:scale_error} and \ref{subfig:scale_time} show that the solution error generally decreases in $S$ at the expense of longer computational time. For reference, the relative error versus iteration and running time with varying $S$ are shown in Figure \ref{subfig:scale_err_time}.

\section{Concluding Remarks}
\label{sec:conclusion}
We have presented a weak adversarial network approach to solve a class of inverse problems numerically.
We leverage the weak formulation of PDEs in the inverse problems, and parameterize the unknown solution as primal neural network and the test function as adversarial network.
The weak formulation and the boundary conditions yield a saddle function in the parameters of the primal network and adversarial network, which only rely on the inverse problem itself but not any other training data.
These parameters are alternately updated until convergence.
We provide a series of theoretical justifications on the convergence of our proposed algorithm.
Our method does not require any spatial discretization, and can be applied to a large class of inverse problems, especially those with high dimensionality and less regularity on solutions.
Numerical experiments have been conducted by applying the proposed method to a variety of challenging inverse problems.
The results suggest promising accuracy and efficiency of our approach.

\section*{Acknowledgements}
GB and YZ are supported in part by NSFC Innovative Group Fund (No.11621101). XY is supported in part by NSF under grants DMS-1620342, CMMI-1745382, DMS-1818886 and DMS-1925263. HZ is supported in part by NSF under grants DMS-1620345 and DMS-1830225, and ONR N00014-18-1-2852.

\appendix
\section{Appendix: Proofs}
\label{app:proof}
\subsection{Proof of Theorem \ref{thm:min_op}}
\label{pf:thm:min_op}
For ease of presentation, our proof of Theorem \ref{thm:min_op} here is based on the problem formulation \eqref{eq:eit}. However, it can be easily modified for the PDEs in many other inverse problems.
\begin{proof}
For any fixed $u\in H^1(\Omega)\cap C(\bar{\Omega})$ and $\gamma\in C(\bar{\Omega})$, the maximum of $\langle \Acal[u,\gamma], \varphi \rangle$ is achievable over $Y \coloneqq \{\varphi\in H^{1}_{0}(\Omega): \|\varphi\|_{H^1}=1\}$ since $\langle \Acal[u,\gamma], \cdot \rangle$ is continuous and $Y$ is closed in $H^1_0(\Omega)$.
Define $h(u,\gamma)=\max_{\varphi \in Y}\langle \Acal[u,\gamma], \varphi \rangle$, then $h(u,\gamma) = \|\Acal[u,\gamma]\|_{op}$ due to the definition of operator norm in \eqref{eq:norm_op}.
On the other hand, let $X=\{(u,\gamma)\in H^1(\Omega) \times C(\Omega): \Bcal[u,\gamma]=0\}$, then it is clear that the minimum value $0$ of $h(u,\gamma)$ over $X$ can be attained at any of the weak solutions.
Hence the minimax problem \eqref{eq:min_op} is well-defined.

Now we show that $(u^*,\gamma^*)$ satisfying $\Bcal[u^*,\gamma^*]=0$ is the solution of the minimax problem \eqref{eq:min_op} if and only if it is a weak solution of the problem \eqref{eq:ip}.
Suppose $(u^{*},\gamma^*)$ is a weak solution of the problem \eqref{eq:ip}, namely $(u^{*},\gamma^*)$ satisfies \eqref{eq:weak_u} for all $\varphi \in Y$, then $\langle \mathcal{A}[u^*,\gamma^*], \varphi \rangle \equiv 0$ for all $\varphi \in Y$.
Therefore, $\|\mathcal{A}[u^{*},\gamma^*]\|_{op}=0$, and $(u^{*},\gamma^*)$ is the solution of the minimax problem \eqref{eq:min_op}.
On the other hand, suppose a weak solution $(\hat{u},\hat{\gamma})$ of \eqref{eq:ip} exists.
Assume that $(u^*,\gamma^*)$ is a minimizer of the problem \eqref{eq:min_op}, i.e., $(u^*,\gamma^*)=\argmin_{(u,\gamma)\in H^1 \times C}h(u,\gamma)$, but not a weak solution of the problem \eqref{eq:ip}, then there exists $\varphi^* \in Y$ such that $\langle \mathcal{A}[u^*,\gamma^*], \varphi^*\rangle > 0$.
Therefore $h(u^*,\gamma^*)=\max_{\varphi \in Y}|\langle \mathcal{A}[u^*,\gamma^*],\varphi \rangle|>0$.
However, as we showed above, $h(\hat{u},\hat{\gamma})=0$ since $(\hat{u},\hat{\gamma})$ is a weak solution of \eqref{eq:ip}, which contradicts to the assumption that $(u^*,\gamma^*)$ is the minimizer of \eqref{eq:min_op}.
Hence $(u^*,\gamma^*)$ must also be a weak solution of \eqref{eq:ip}, i.e., $(u^*,\gamma^*)$ satisfies \eqref{eq:weak_u}.
\end{proof}
\subsection{Proof of Lemma \ref{lem:dLint}}
\label{pf:lem:dLint}
\begin{proof}
Due to the definition of $\Lint(\theta)$ in \eqref{eq:Lint} and the optimality of $\eta(\theta)$, we know that $\Lint(\theta) = E(\theta,\eta(\theta))$.
Therefore, we have
\begin{equation}\label{eq:dLint}
\nabla_\theta \Lint(\theta) = \partial_{\theta} E(\theta,\eta(\theta)) + \partial_{\eta} E(\theta,\eta(\theta)) \nabla_{\theta} \eta(\theta).
\end{equation}
Now we form the Lagrange function $\Lcal(\theta,\eta,\mu) = E(\theta,\eta) + \mu(\frac{1}{2}|\eta|^2 - B)$ for the  maximization problem $\max_{|\eta|^2 \le 2B} E(\theta,\eta)$.
Then the Karush-Kuhn-Tucker (KKT) condition of $\eta(\theta)$ is given by
\begin{subequations}
\label{eq:kkt}
\begin{align}
\partial_\eta\Lcal(\theta,\eta(\theta),\mu) = \partial_\eta E(\theta,\eta(\theta)) + \mu(\theta) \eta(\theta) & = 0, \label{eq:kktd} \\
\mu(\theta) \del[1]{(1/2) \cdot |\eta(\theta)|^2 - B} & = 0, \label{eq:kkts}\\
\mu(\theta)\ge 0,\quad |\eta(\theta)|^2 & \le  2 B. \label{eq:kktc}
\end{align}
\end{subequations}
The complementary slackness condition \eqref{eq:kkts} implies that
\begin{equation}
\label{eq:dkkts}
\nabla_\theta \mu(\theta) \del[1]{(1/2)\cdot|\eta(\theta)|^2 - B} + \mu(\theta)\eta(\theta)\nabla_{\theta}\eta(\theta) = 0
\end{equation}
If $\mu(\theta)=0$, then we know $\partial_\eta E(\theta,\eta(\theta)) = 0$ due to \eqref{eq:kktd} and hence \eqref{eq:dLint} reduces to $\nabla_{\theta} \Lint(\theta) = \partial_{\theta} E(\theta, \eta(\theta))$.
If $\mu(\theta)>0$, then $|\eta(\theta)|^2 = 2B$ due to \eqref{eq:kkts}, and hence \eqref{eq:dkkts} implies $\mu(\theta)\eta(\theta)\nabla_{\theta}\eta(\theta)=0$.
Thus multiplying \eqref{eq:kktd} by $\nabla_\theta \eta(\theta)$ yields $\partial_{\eta} E(\theta,\eta(\theta)) \nabla_{\theta} \eta(\theta)=0$, from which we can see \eqref{eq:dLint} also reduces to $\nabla_{\theta} \Lint(\theta) = \partial_{\theta} E(\theta, \eta(\theta))$.
\end{proof}
\subsection{Proof of Lemma \ref{lem:integral}}
\label{pf:lem:integral}
\begin{proof}
The first moment, i.e., expectation of $\hat{\Psi}$, can be computed as follows:
\begin{equation}
\mathbb{E}[\hat{\Psi}] = \mathbb{E}[\psi/\rho] = \int_{\Omega} \rho \frac{\psi}{\rho} \dif x = \int_{\Omega}\psi \dif x = \Psi.
\end{equation}
To compute the second moment of $\hat{\Psi}$, we first observe that the variance of $\hat{\Psi}$ is
\begin{equation*}
\mathrm{V}(\hat{\Psi}) = \mathrm{V}\del[2]{\frac{1}{N}\sum_{i=1}^N\frac{\psi(x^{(i)})}{\rho(x^{(i)})}} = \frac{1}{N}\mathrm{V}\del[2]{\frac{\psi}{\rho}}.
\end{equation*}
Note that the variance of $\psi/\rho$ is
\begin{equation*}
\mathrm{V}\del[2]{\frac{\psi}{\rho}} = \mathbb{E}\sbr[2]{\del[2]{\frac{\psi}{\rho}}^2} - \del[2]{\mathbb{E}\sbr[2]{\frac{\psi}{\rho}}}^2 = \int_{\Omega} \frac{\psi^2}{\rho} \dif x - \del[2]{\int_{\Omega}\psi \dif x}^2 = \int_{\Omega} \frac{\psi^2}{\rho} \dif x - \Psi^2
\end{equation*}
Hence the second moment of $\hat{\Psi}$ is
\begin{equation*}
\mathbb{E}[\hat{\Psi}^2] = \mathrm{V}(\hat{\Psi}) + \mathbb{E}[\hat{\Psi}]^2 = \frac{1}{N}\del[2]{\int_{\Omega} \frac{\psi^2}{\rho} \dif x - \Psi^2} + \Psi^2 = \frac{N-1}{N}\Psi^2 + \frac{1}{N} \int_{\Omega} \frac{\psi(x)^2}{\rho(x)} \dif x,
\end{equation*}
which completes the proof.
\end{proof}

\subsection{Proof of Theorem \ref{thm:converge}}
\label{pf:thm:converge}
\begin{proof}
Due to the parameterization of $(u_{\theta},\gamma_{\theta})$ using finite-depth neural network \eqref{eq:u_param} and the compactness of $\Theta \coloneqq \{\theta : |\theta| \le \sqrt{2B} \}$, we know $\partial^{\alpha}u_{\theta}$ and $\gamma_{\theta}$ have Lipschitz continuous gradient with respect to $\theta$ for all $|\alpha|\le1$.
As $\Omega$ is bounded and $\partial^\alpha u_{\theta},\gamma_{\theta}\in C(\bar{\Omega})$, there exists $M>0$ such that $L(\theta)$ has $M$-Lipschitz continuous gradient $\nabla_{\theta} L(\theta)$, since $L(\theta)$ is composed of integrals of $\partial^\alpha u_{\theta}$ and $\gamma_{\theta}$ over $\Omega$.

Recall that the projected stochastic gradient descent step \eqref{eq:gd}, started from initial $\theta_1$, generates the sequence $\{\theta_j\}$ as follows:
\begin{equation}
\label{eq:theta_iter}
\tjp = \Pi(\tj - \tau G_j) = \argmin_{\theta \in \Theta} \del[2]{G_j^{\top}\theta + \frac{1}{2\tau} |\theta - \tj|^2}
\end{equation}
where $G_j$ denotes the stochastic gradient of $L(\theta)$ at $\tj$ using $N_r$ ($N_b$ resp.) sample collocation points in $\Omega$ (on $\partial\Omega$ resp.) with $N_r,N_b=O(N)$.
We let $g_j \coloneqq \nabla_{\theta}L(\tj)$ denote the true (but unknown) gradient of $L$ at $\tj$, and define a companion sequence $\{\btj\}$ using $g_j$ as
\begin{equation}
\label{eq:btheta_iter}
\btjp = \Pi(\tj - \tau g_j) = \argmin_{\theta \in \Theta} \del[2]{g_j^{\top}\theta + \frac{1}{2\tau} |\theta - \tj|^2}.
\end{equation}
Note that $\{\btj\}$ is not computed in practice (computation of $\btj$ is not possible as $g_j$ is unknown), but only defined for convergence analysis here.
Also note that $\Theta$ and $\Omega$ are bounded, and hence all integrals of $\partial^\alpha u_{\theta}$ and $\gamma_{\theta}$ are bounded, we know $G_j$ is an unbiased estimate of $g_j$ with bounded variance, denoted by $\sigma^2>0$, according to Lemma \ref{lem:integral}.
Moreover, Lemma \ref{lem:integral} implies that there exists $E>0$ dependent on $B$, $\Omega$, $\Acal$, and $\Bcal$ only ($E$ is the bound of $\|\Acal[u_\theta,\gamma_\theta]\|_{op}^2$ and $\|\Bcal[u_\theta,\gamma_\theta]\|_{L^2(\partial\Omega)}^2$ due to the boundedness of $\Theta$ and $\Omega$) the integral such that $\mathbb{E}[|G_j-g_j|^2] \le \sigma^2 \le E/N$ as $N_r,N_b=O(N)$.

Now we are ready to verify the convergence of the projected SGD iterations \eqref{eq:theta_iter}. First, the $M$-Lipschitz continuity of $\nabla_{\theta}L$ implies that
\begin{equation}
\label{eq:lip1}
L(\tjp) \le L(\tj) + g_j^{\top}e_j + \frac{M}{2}|e_j|^2,
\end{equation}
where we denote $e_j \coloneqq \tjp - \tj$ for all $j$. Also due to the $M$-Lipschitz continuity of $\nabla_{\theta} L$, we have
\begin{equation}
\label{eq:lip2}
-L(\btjp) \le - L(\tj) - g_j^{\top}\be_j + \frac{M}{2}|\be_j|^2,
\end{equation}
where we denote $\be_j \coloneqq \btjp - \tj$.
Note that $\Gcal(\theta_j) = \tau^{-1}[\tj - \Pi(\tj - \tau g_j)]= \tau^{-1}(\tj - \btjp) = -\tau^{-1}\be_j$, whose magnitude is what we want to bound eventually.
Furthermore, due to the optimality of $\tjp$ in \eqref{eq:theta_iter} (which is convex in $\theta$), we know that
\begin{equation}
\label{eq:theta_opt}
0 \le \del[2]{G_j + \frac{\tjp-\tj}{\tau}}^{\top}(\btjp - \tjp) = \del[2]{G_j + \frac{e_j}{\tau}}^{\top}(\be_j - e_j).
\end{equation}
Adding \eqref{eq:lip1}, \eqref{eq:lip2}, and \eqref{eq:theta_opt} yields
\begin{equation}
\label{eq:est1}
L(\tjp) - L(\btjp) \le (g_j - G_j)^{\top}(e_j-\be_j) + \frac{e_j^{\top}(\be_j-e_j)}{\tau} + \frac{M}{2}|e_j|^2 + \frac{M}{2}|\be_j|^2.
\end{equation}
Repeating \eqref{eq:lip1}, \eqref{eq:lip2}, \eqref{eq:theta_opt}, and \eqref{eq:est1} with $\tjp$ and $\btjp$ replaced by $\btjp$ and $\tj$ respectively, and using the optimality of $\btjp$ in \eqref{eq:btheta_iter} with $g_j$, we obtain
\begin{equation}
\label{eq:est2}
L(\btjp) - L(\tj) \le - \del[2]{\frac{1}{\tau}- \frac{M}{2}} |\be_j|^2.
\end{equation}
Adding \eqref{eq:est1} and \eqref{eq:est2} yields
\begin{equation}
\label{eq:est3}
L(\tjp) - L(\tj) \le (g_j - G_j)^{\top}(e_j-\be_j) + \frac{e_j^{\top}(\be_j-e_j)}{\tau} + \frac{M}{2}|e_j|^2 - \del[2]{\frac{1}{\tau}-M}|\be_j|^2.
\end{equation}
Now due to Cauchy-Schwarz inequality, the definitions of $\tjp$ and $\btjp$ in \eqref{eq:theta_iter} and \eqref{eq:btheta_iter}, and that the projection $\Pi$ onto the convex set $\Theta$ is a non-expansive operator (i.e., $|\Pi(\theta)-\Pi(\hat{\theta})| \le |\theta - \hat{\theta}|$ for any $\theta,\hat{\theta}$), we can show that
\begin{align}
(g_j - G_j)^{\top}(e_j-\be_j)
=\ & (g_j - G_j)^{\top}(\tjp-\btjp) = (g_j - G_j)^{\top}(\Pi(\theta_j-\tau G_j)-\Pi(\theta_j-\tau g_j)) \nonumber \\
\le\ & |g_j-G_j| \; |\Pi(\tj-\tau G_j) - \Pi(\tj-\tau g_j)| \le \tau |g_j - G_j|^2.  \label{eq:cs}
\end{align}
Moreover, we have that
\begin{equation}
\label{eq:inner2sum}
\frac{e_j^{\top}(\be_{j} - e_{j})}{\tau}=\frac{1}{2\tau}\del[2]{|\be_j|^2 - |e_j|^2 - |e_j - \be_j|^2}.
\end{equation}
Substituting \eqref{eq:cs} and \eqref{eq:inner2sum} into \eqref{eq:est3}, we obtain
\begin{equation}
\label{eq:est4}
L(\tjp) - L(\tj) \le \tau |g_j-G_j|^2 - \del[2]{\frac{1}{2\tau}-M}|\be_j|^2  - \del[2]{\frac{1}{2\tau}-\frac{M}{2}}|e_j|^2 - \frac{1}{2\tau}|e_j-\be_j|^2.
\end{equation}
Taking expectation on both sides of \eqref{eq:est4} and discarding the last negative term, we obtain
\begin{equation}
\label{eq:est5}
\del[2]{\frac{1}{2}-\tau M}\tau \mathbb{E}[|\Gcal(\theta_j)|^2] = \del[2]{\frac{1}{2\tau}-M}\mathbb{E}[|\be_j|^2] \le L_j - L_{j+1} +  \tau \sigma^2 - \del[2]{\frac{1}{2\tau}-\frac{M}{2}}\mathbb{E}[|e_j|^2]
\end{equation}
where we used the fact $\mathbb{E}[|G_j-g_j|^2] \le \sigma^2$ and the notation $L_j \coloneqq \mathbb{E}[L(\theta_j)]$.
Now taking sum of \eqref{eq:est5} for $j=1,\dots,J$, dividing both sides by $(\frac{1}{2}-\tau M)\tau J$, and setting $\tau=\frac{1}{4M}$ (hence $\frac{1}{2}-\tau M = \frac{1}{4}$ and $\frac{1}{\tau}-\frac{M}{2}=\frac{7M}{2} > 0$), we know that
\begin{equation}
\min_{1\le j \le J} \mathbb{E}[|\Gcal(\theta_j)|^2] \le \frac{1}{J}\sum_{j=1}^J \mathbb{E}[|\Gcal(\theta_j)|^2] \le \frac{16M(L_1-L_{J+1})}{J} + 4\sigma^2 \le \frac{16M(L_1-L^*)}{J} + \frac{4E}{N} \le \varepsilon
\end{equation}
by choosing per-iteration sample complexity $N$ and iteration number $J$ as $N=J=[16M(L_1-L^*)+4E]\varepsilon^{-1}=O(\varepsilon^{-1})$ where $L^*\coloneqq \min_{\theta \in \Theta}L(\theta) \ge 0$ (and hence $L_{J+1}=\mathbb{E}[L(\theta_{J+1})]\ge J^*$).
This completes the proof.
\end{proof}

\section{Appendix: Problem Setting}
\label{app:prob_set}
The functions and parameters used in our experiments are summarized in table \ref{tab:prob_set}.

\begin{table}[ht]
\scriptsize
\newcommand{\tabincell}[2]{\begin{tabular}{@{}#1@{}}#2\end{tabular}}
\caption{Problem Settings for Tests 1-4 in subsection \ref{subsec:test}, where $\gamma^*$ denotes the true conductivity, $u^*$ denotes the true potential, and $f$ is the source function.} 
\label{tab:prob_set}
\centering 
\begin{tabular}{|c|c|c|c|} 
\hline 
Problem & \tabincell{c}{\\$\gamma^*$\\ \\} & $u^*$ & $f$ \\ \hline 
Test 1 & \tabincell{c}{$2(\exp(-|x-c_1|_{\Sigma_1}^2) + \exp(-|x-c_2|_{\Sigma_2}^2))$ \\ where $\Sigma_1=\mathrm{diag}(1.25,5,0,0,0)$, \\ $\Sigma_2=\mathrm{diag}(5,1.8,0,0,0)$, \\$c_1=(-0.5,0.5,0,0,0)$, \\$c_2 = (0.5,-0.5,0,0,0)$, \\and $|x|_\Sigma^2 \coloneqq x^{\top}\Sigma x$.} & $\cos(|x|^2)$ & \tabincell{c}{\\ $ 8\sum^{2}_{i,j=1}[\Sigma_j]_{ii}(x_i-[c_j]_i)x_i\sin(|x|^2)\exp(-|x-c_j|^2_{\Sigma_j})
$\\+$\gamma^*(2d\sin(|x|^2)+4|x|^2\cos(|x|^2))$\\where $[\Sigma]_{ij}$ and $[c]_j$ stand for \\ the $(i,j)$-th entry of the matrix $\Sigma$ \\and $j$th component of the vector $c$, respectively.\\ \\} \\ \hline
\tabincell{c}{Test 2\\ \& Test 3} & \tabincell{c}{\\ $ 0.5+1.5/(1+\delta(x))$ \\where $\delta(x)=\exp((|x-c|^2_{\Sigma}-0.6^2)/0.02)$ \\with $\Sigma=\mathrm{diag}(0.81,2,0.09,\dots,0.09)$ \\and $c=(0.1,0.3,0,\dots,0).$ \\ \\} & $|x|^2$ & $\frac{6(x-c)^{\top}\Sigma x}{\lambda(1/\delta(x)+2+\delta(x))}-2d\gamma^*$ \\\hline
Test 4(1) & \tabincell{c}{\\$0.5+3.5/(1+\delta_1(x))+1.5/(1+\delta_2(x))$ \\ where $\delta_1(x)=\exp((|x-c_1|^2_{\Sigma_1}-0.4^2)/0.02)$ \\and $\delta_2(x)=\exp((|x-c_2|^2_{\Sigma_2}-0.4^2)/0.02)$.\\ with $\Sigma_1=\mathrm{diag}(0.81,2,0.09,0.09,0.09)$, \\$\Sigma_2=\mathrm{diag}(2,0.81,0.09,0.09,0.09)$,\\ and $c_1=(-0.5,-0.5,0,0,0)$, \\$c_2=(0.5,0.5,0,0,0)$.\\  \\} & $|x|^2$
& \tabincell{c}{$\frac{14(x-c_1)^{\top}\Sigma_1 x}{\lambda(1/\delta_1(x)+2+\delta_1(x))}+\frac{6(x-c_2)^{\top}\Sigma_2 x}{\lambda(1/\delta_2(x)+2+\delta_2(x))}$\\ \\$-2d\gamma^*$} \\ \hline
Test 4(2) &\tabincell{c}{\\$0.5+1.5/(1+\delta_1(x))+1.5/(1+\delta_2(x)+\delta_3(x))$\\ where $\delta_1(x)=\exp((|x-c|^2_{\Sigma_2}-0.4^2)/0.02)$\\ $\delta_2(x)=\exp((|x_1+0.5|-0.15)/0.02)$\\ and $\delta_3(x)=\exp((|x_2|-0.6)/0.02)$ \\with $c=(0.55, 0, 0, 0, 0)$ \\and $\Sigma=\mathrm{diag}(1, 4, 0, 0, 0).$\\  \\}& $|x|^2$
&\tabincell{c}{\\$\frac{6(x-c)^T\Sigma x}{\lambda(1/\delta_1(x)+2+\delta_1(x))}+\frac{3(|x_1+0.5|\delta_2(x)+|x_2|\delta_3(x))}{\lambda(1+\delta_2(x)+\delta_3(x))^2}$\\ \\$-2d\gamma^*$\\ \\} \\ \hline
Test 4(3) & \tabincell{c}{\\ $0.5+\sum^{3}_{j=1}1.5/(1+\delta_{j1}(x)+\delta_{j2}(x))$ \\where $\delta_{j1}(x)=\exp((|x_1-c_j(1)|-r_{j}(1))/0.02)$, \\$\delta_{j2}(x)=\exp((|x_2-c_j(2)|-r_j(2))/0.02),$ \\for $j=1,2,3$ with $c_1=(-0.5, 0),$ \\$c_2=(-0.1, 0.6)$, $c_3= (-0.1, -0.6)$ \\and $r_1=(0.15, 0.8),$ $r_2=r_3=(0.55, 0.2).$ \\ \\} & $|x|^2$
&\tabincell{c}{$3\sum^{3}_{j=1}\frac{|x_1-c_j(1)|\delta_{j1}(x)+|x_2-c_j(2)|\delta_{j2}(x)}{\lambda(1+\delta_{j1}(x)+\delta_{j2}(x))^2}$\\ \\$-2d\gamma^*$}\\
\hline 
\end{tabular}
\end{table}

\section{Appendix: Recorded errors and running times in Tests 9 and 10}
\label{app:error_time_table}
The recorded errors and running times for {Test 9} and {Test 10} were present in table \ref{tab:error_time_layer_neuron} and table \ref{tab:error_time_Nr_Nb}, respectively.
\begin{table}[t]
\caption{Test 9 result on relative error of recovered conductivity $\gamma_\theta$ (top) and running time (bottom) using various combinations of $(K,d')$ for problem dimension $d=5$, where $K$ is the layer number and $d'$ is the per-layer neuron number.}
\label{tab:error_time_layer_neuron}
\centering
\begin{tabular}{|c|c|c|c|c|}
\hline
$d'$ & $K=5$ & $K=7$ & $K=9$ & $K=11$ \\ \hline \hline
5  & 0.060347& 0.040950  & 0.014905  & 0.019489      \\ \hline
10 & 0.053842& 0.029382  & 0.013325  & 0.011165      \\ \hline
20 & 0.017955& 0.016862  & 0.010916  & 0.011490      \\ \hline
40 & 0.012390& 0.010213  & 0.004422  & 0.005309      \\ \hline  \hline
5  & 1391.76(s)& 1432.90(s)& 1438.87(s)& 1458.12(s)     \\ \hline
10 & 1435.78(s)& 1468.95(s)& 1520.67(s)& 1568.49(s)     \\ \hline
20 & 1447.66(s)& 1522.84(s)& 1596.83(s)& 1614.04(s)      \\ \hline
40 & 1539.93(s)& 1623.20(s)& 1717.53(s)& 1809.30(s)    \\ \hline
\end{tabular}
\end{table}
\begin{table}[t]
\caption{Test 10 result on relative error of recovered conductivity $\gamma_\theta$ (top) and running time (bottom) with various combination of $(N_r,N_b)$ for problem dimension $d=5$, where $N_r$ is the number of sampled collocation points inside the region $\Omega$ and $N_b$ is the number of those on the boundary $\partial \Omega$.}
\label{tab:error_time_Nr_Nb}
\centering
\begin{tabular}{|c|c|c|c|c|}
\hline
$N_b$ & $N_r=25$K & $N_r=50$K & $N_r=100$K & $N_r=200$K \\ \hline \hline
$10\times 2d$ & 0.019023 & 0.020007  & 0.020898  & 0.012208      \\ \hline
$20\times 2d$ & 0.013999 & 0.012920  & 0.009681  & 0.010050      \\ \hline
$40\times 2d$ & 0.010668 & 0.012292  & 0.012053  & 0.010385      \\ \hline
$80\times 2d$ & 0.01061 & 0.009207  & 0.010200  & 0.007267      \\ \hline \hline
$10\times 2d$  & 528.88(s)&  554.64(s)&  556.08(s)&  552.86(s)      \\ \hline
$20\times 2d$  & 915.66(s)&  898.32(s)&  932.98(s)&  928.78(s)     \\ \hline
$40\times 2d$  & 1597.66(s)& 1577.03(s)& 1556.92(s)& 1590.87(s)      \\ \hline
$80\times 2d$ & 2798.17(s)& 2802.79(s)& 2795.28(s)& 2801.45(s)     \\ \hline
\end{tabular}
\end{table}

\section{Appendix: Comparison with classical numerical methods}
\label{sec:classical}
To further evaluate the proposed method, we provide an example that compares IWAN and a classical finite difference method (FDM) on a 2D problem in Test 2. We would like to acknowledge an anonymous reviewer for suggesting this valuable comparison. It is worth pointing out that the majority of existing numerical methods, such as FDM, require the knowledge of Dirichlet-to-Neumann (DtN) map for the EIT problem. However, to be consistent with the settings used in the present work, we conduct the comparison with FDM under the same setting of Test 2, where a DtN map is not available but only the boundary conditions of $u$ and $\gamma$ in \eqref{eq:eit} are given.

\begin{figure}[b]
\centering
\begin{subfigure}[b]{.3\textwidth}
\includegraphics[width=\textwidth]{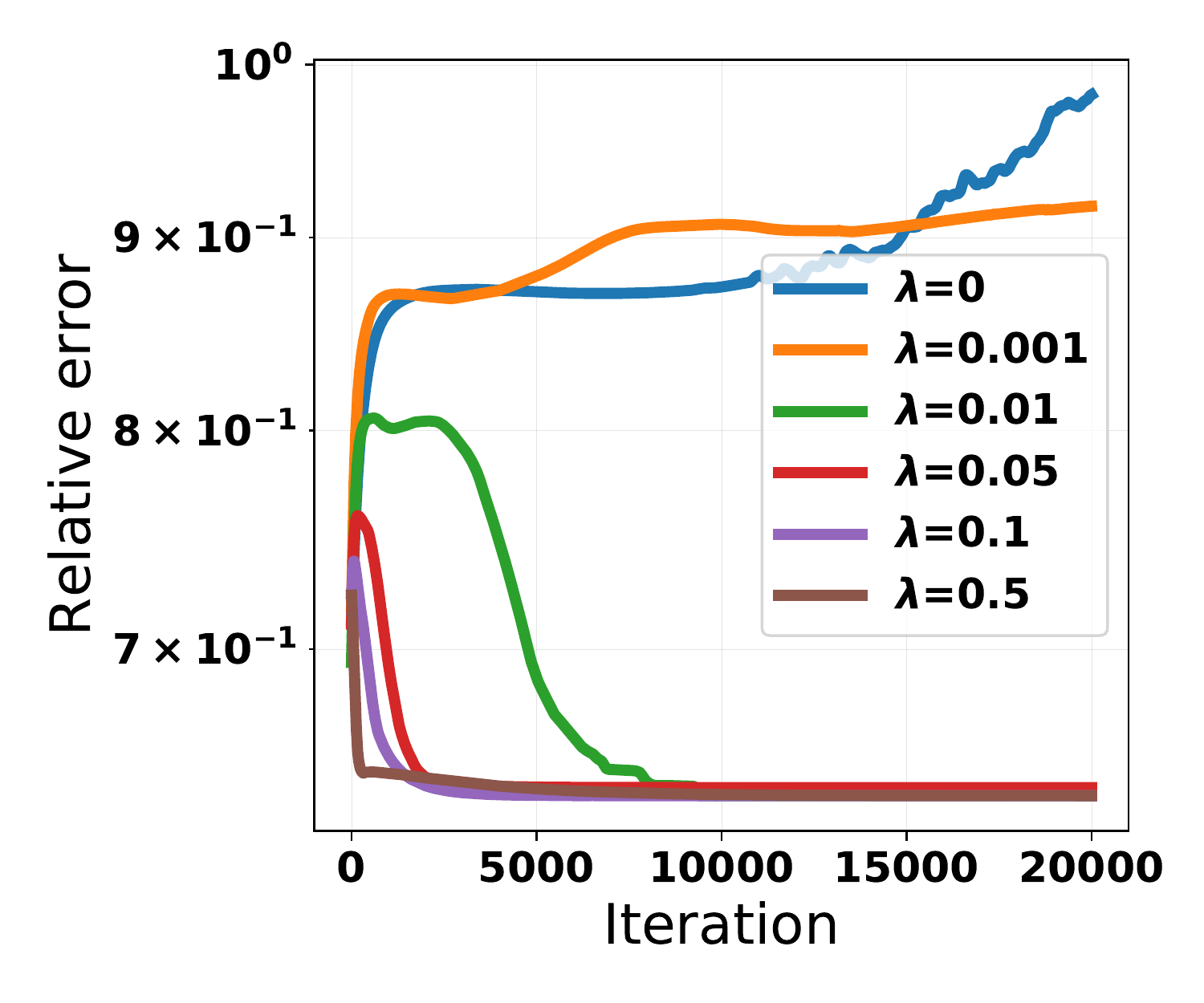}
\caption{Error vs.~iter by FDM}
\label{subfig:2d_lsq_err_iter_one}
\end{subfigure}
\begin{subfigure}[b]{.3\textwidth}
\includegraphics[width=\textwidth]{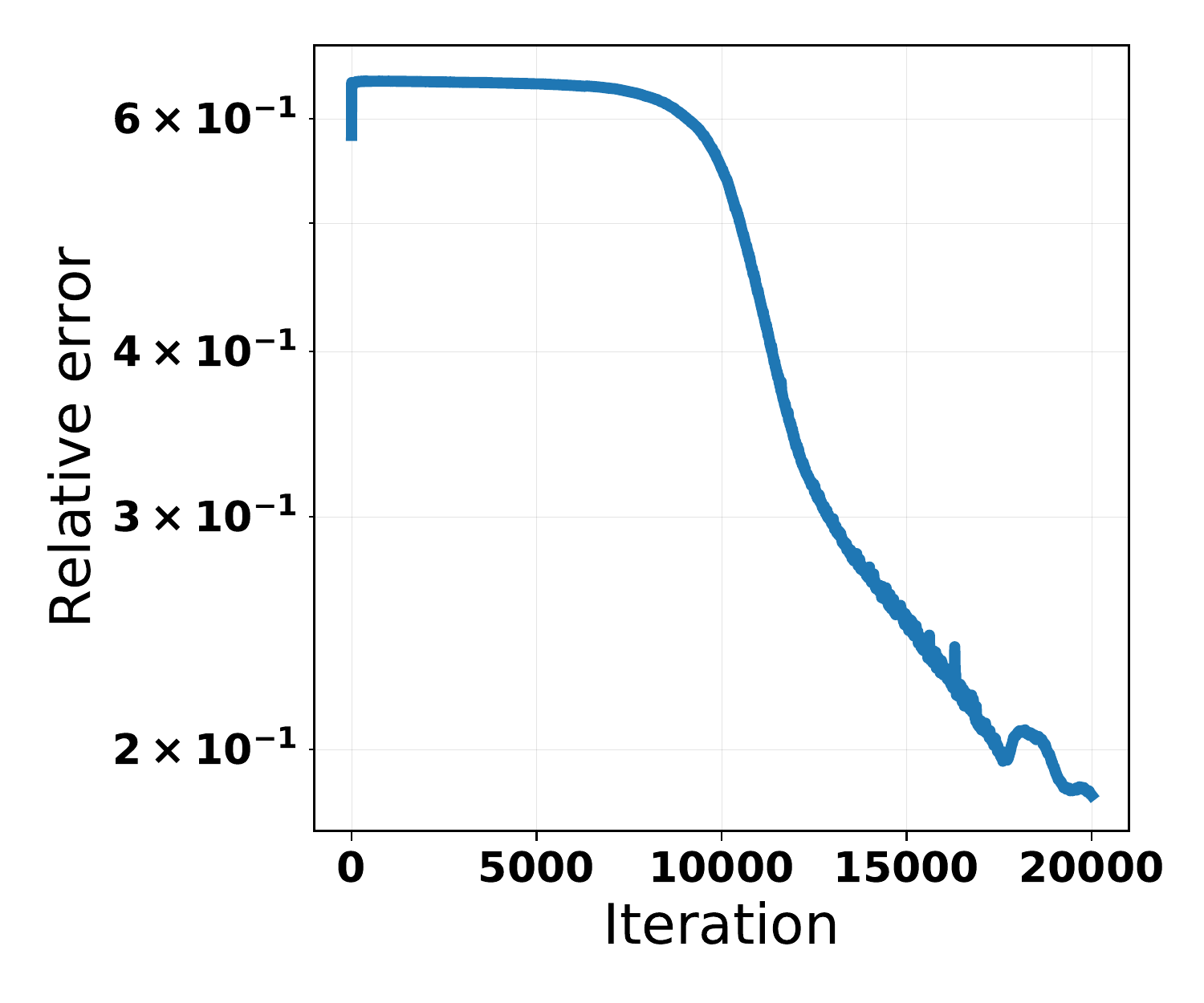}
\caption{Error vs.~iter by IWAN}
\label{subfig:2d_wan_err_iter}
\end{subfigure}
\begin{subfigure}[b]{.7\textwidth}
\includegraphics[width=\textwidth]{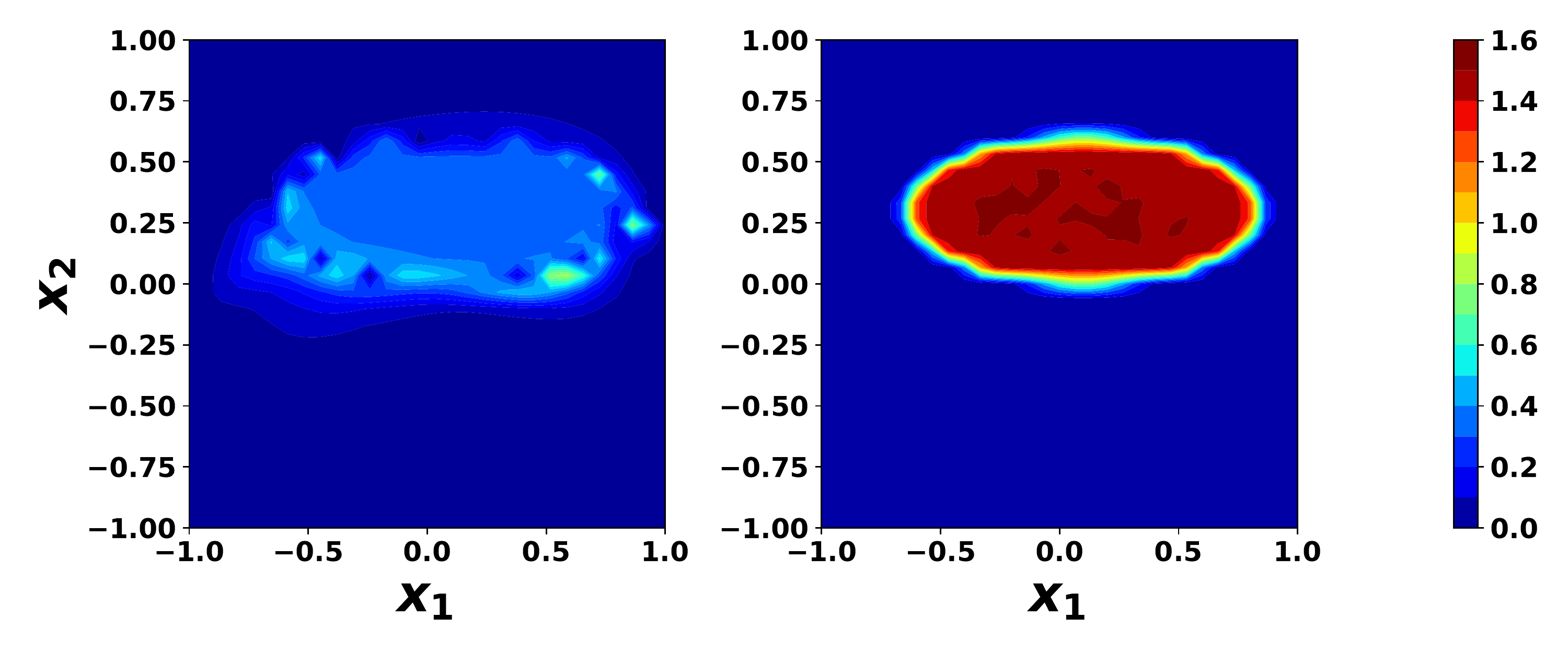}
\caption{$|\gamma_{\theta}-\gamma^*|$ by IWAN (left) and FDM (right)}
\label{subfig:2d_abs_error}
\end{subfigure}
\begin{subfigure}[b]{0.45\textwidth}
\centering
\includegraphics[width=.7\linewidth]{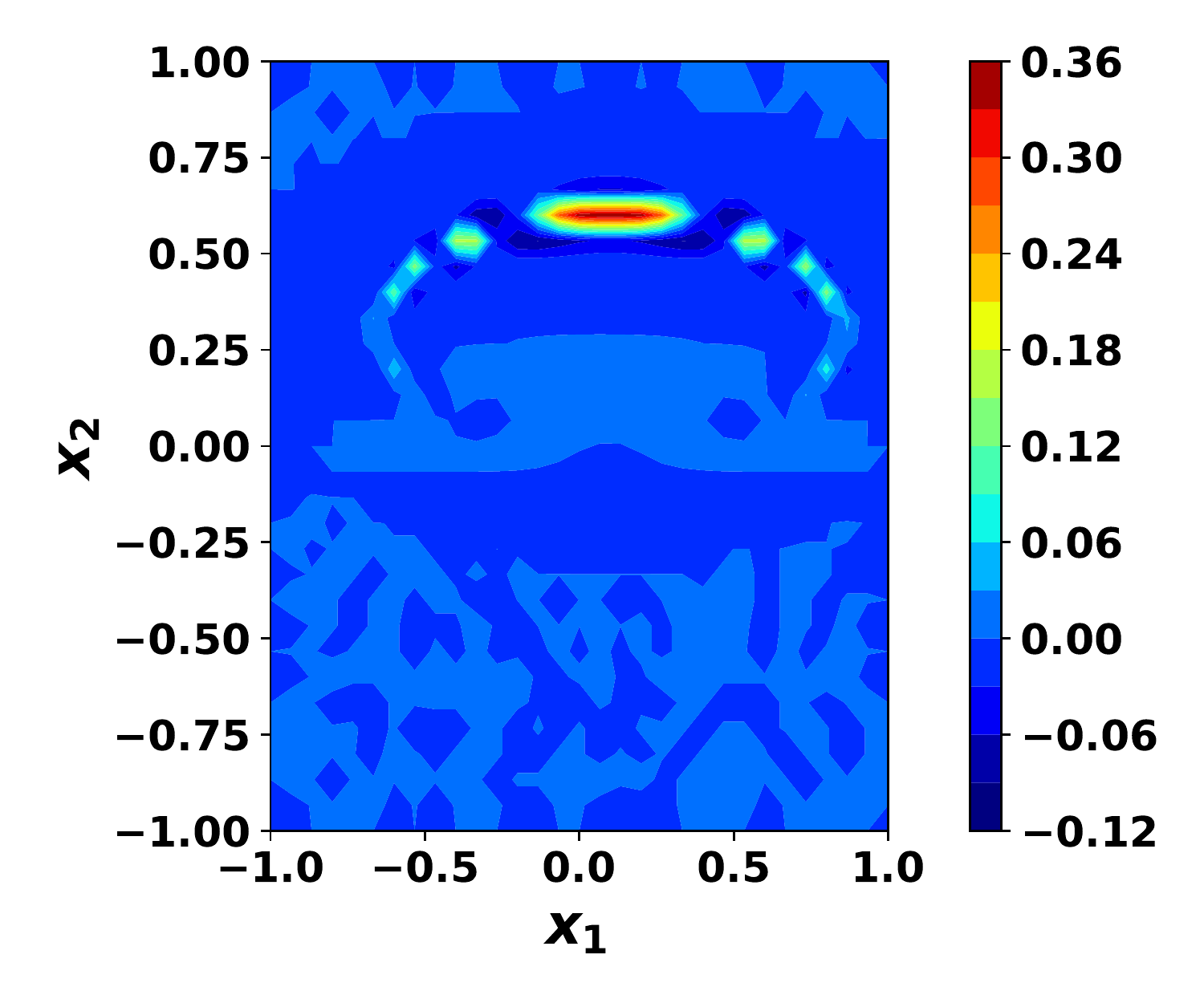}
\caption{$|-\nabla(\gamma\nabla u)-f|$ by FDM with $\lambda=0.1$}
\label{subfig:2d_abs_pde}
\end{subfigure}\qquad \qquad
\begin{subfigure}[b]{0.45\textwidth}
\centering
\includegraphics[width=.7\linewidth]{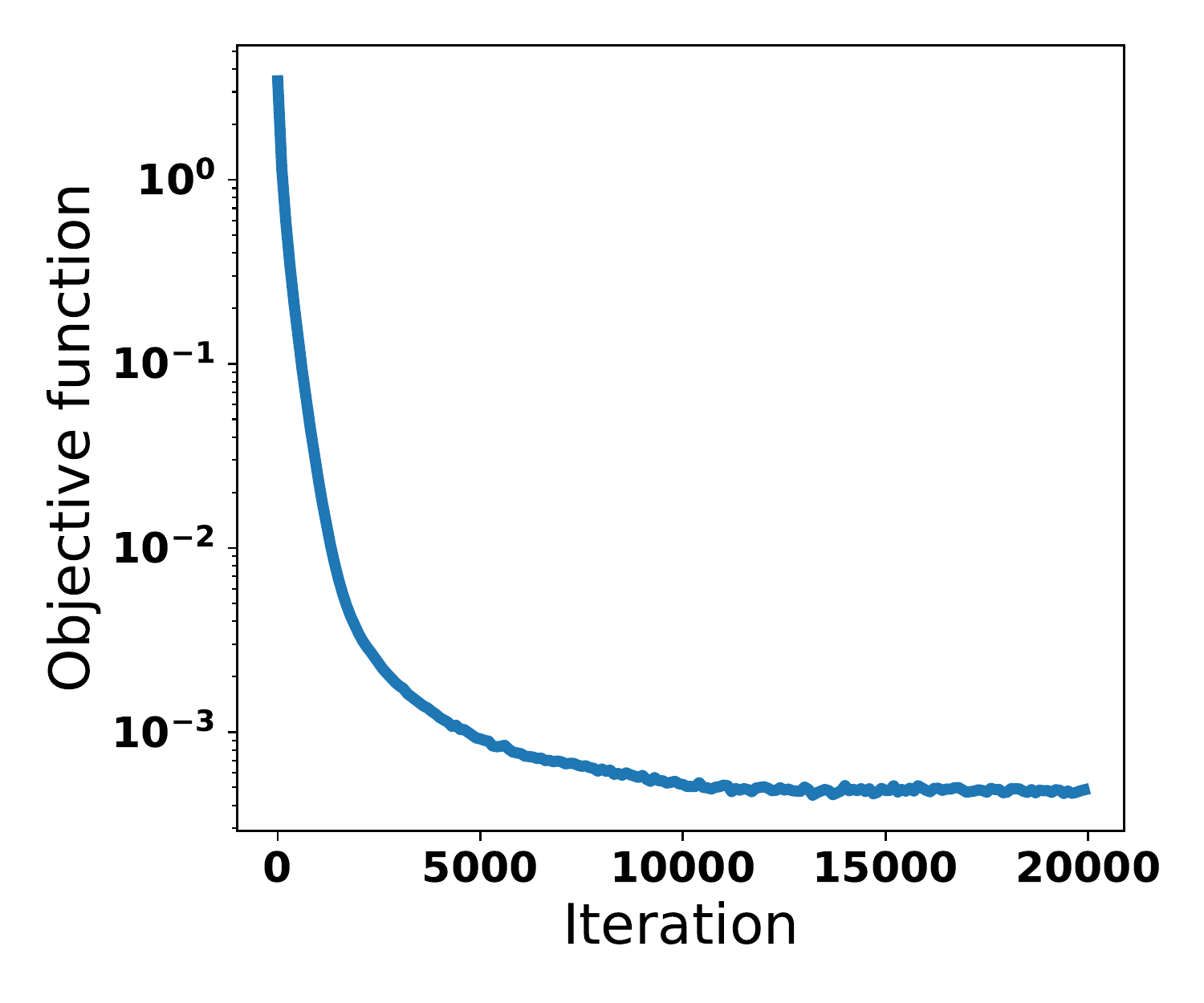}
\caption{Objective vs. iter by FDM with $\lambda=0.1$}
\label{subfig:2d_lsq_loss_iter}
\end{subfigure}
\caption{Test 2 on the comparison of IWAN and FDM to recover the conductivity $\gamma^*$ with in 2D case. (a) Relative error versus iteration number obtained by FDM with different $\lambda$. (b) Relative error versus iteration obtained by IWAN. (c) Pointwise absolute error $|\gamma-\gamma^*|$ with $\gamma$ obtained by IWAN (left) and FDM (right). (d)  $|-\nabla(\gamma\nabla u)-f|$ by FDM. (e) Objective function value versus iteration number by FDM.}
\label{fig:2d_compare}
\end{figure}

In FDM, we discretize the domain into $31\times 31$ mesh grids (about $900$ unknowns for each of $u$ and $\gamma$). We also experiment with higher resolution but it does not improve solution quality; see later for more explanations. We approximate the partial derivatives in the PDE by finite differences. As the problem is underdetermined, we use regularization and formulate as a minimization problem of $(u,\gamma)$, where the objective function is the sum of two terms: the mean square error of the PDE, and the TV regularization on $\gamma$. For comparison, we choose $4$ hidden layers with $15$ neurons per-layer to parameterize $u_{\theta}$ and $\gamma_{\theta}$ (each with $<800$ unknowns).
We set the number of collocation points to $N_r=1,000$ and $N_b=120$ (similar to the discretization resolution of FDM). For FDM, we test different regularization hyperparameter $\lambda$ (the weight of the TV term), and show the result in Figure \ref{subfig:2d_lsq_err_iter_one}. We observe that FDM achieved the best result when $\lambda=0.1$, which is used to generate the other images in Figure \ref{fig:2d_compare}.
The relative error obtained by IWAN is shown in Figure \ref{subfig:2d_wan_err_iter}. Figure \ref{subfig:2d_abs_error} shows the absolute error $|\gamma_{\theta}-\gamma|$ obtained by IWAN and FDM ($\gamma_{\theta}$ is the recovered conductivity by IWAN or FDM, and $\gamma$ is the ground truth), respectively, which demonstrates that IWAN can faithfully recover $\gamma_{\theta}$ but FDM cannot in the setting of Test 2. Figure \ref{subfig:2d_lsq_loss_iter} shows the objective function value versus iteration by FDM, which suggests that FDM has converged. We also show $|-\nabla(\gamma\nabla u)-f|$ obtained by FDM in Figure \ref{subfig:2d_abs_pde}, which shows that the $(u, \gamma)$ obtained by FDM indeed satisfies the PDE approximately. In addition, we increase the domain discretization resolution up to $101\times 101$ for FDM, but did not observe any noticeable improvement, and hence we omitted the results here.


\end{document}